\documentclass[11pt, a4paper]{article}
\usepackage[comma, authoryear, numbers]{natbib}
\usepackage{amsmath,amsfonts,amsthm,bm}
\usepackage[parfill]{parskip}
\usepackage{hyperref}
\usepackage{graphicx}
\usepackage[a4paper, margin=25mm]{geometry}
\usepackage{bm}
\usepackage{bbm}
\usepackage{breqn}
\usepackage[toc,page]{appendix}
\usepackage{tabularx}
\usepackage{amsthm}
\usepackage[dvipsnames]{xcolor}
\usepackage{multirow}
\usepackage{longtable}
\usepackage{booktabs}
\usepackage{bibentry}
\usepackage[multiple]{footmisc}
\usepackage{mathtools}
\usepackage{algpseudocode}
\usepackage{algorithmicx}
\usepackage{algorithm}
\usepackage{subcaption}
\usepackage{etoolbox}
\usepackage{dsfont}
\usepackage{breqn}
\usepackage{tabulary}
\usepackage{tocbibind}
\usepackage{comment}

\DeclareMathOperator*{\argmax}{argmax} 

\newcommand{\fa}[0]{\ \forall \ }

\newcommand{\T}[0]{\mathcal{T}}

\newcommand{\E}[0]{\mathbb{E}}
\newcommand{\R}{\mathbb{R}}
\renewcommand{\P}[0]{\mathbb{P}}
\newcommand{\m}[1]{\mathcal{#1}}

\renewcommand{\l}{\left}
\renewcommand{\r}{\right}

\usepackage{wrapfig}
\newtheorem{theorem}{Theorem}[section]

\newtheorem{lemma}[theorem]{Lemma}
\newtheorem{proposition}[theorem]{Proposition}
\allowdisplaybreaks

\title{\textbf{
    Parametric Distributionally Robust Optimisation Models for Budgeted Multi-period Newsvendor Problems
    }
}
\author{
    Ben Black\footnote{AnyVan Ltd, London, England, UK. Email:
    \href{mailto:benb17995@hotmail.co.uk}{benb17995@hotmail.co.uk}}
    \footnote{Corresponding author.}\ , 
    Trivikram Dokka\footnote{Advanced Analytics Group, Air Products Plc, United Kingdom. Email:
    \href{mailto:Trivikram.Dokka@yahoo.co.uk}{Trivikram.Dokka@yahoo.co.uk}}\ , 
    Christopher Kirkbride\footnote{Department of Management Science, Lancaster
    University Management School, United Kingdom. Email:
    \href{mailto:c.kirkbride@lancaster.ac.uk}{c.kirkbride@lancaster.ac.uk}.}
}
\date{}
\graphicspath{{images_tables/}}

\newcommand{\pderiv}[2]{\frac{\partial#1}{\partial #2}}

\usepackage{amsthm}
\usepackage[dvipsnames]{xcolor}

\begin{document}
\setcitestyle{nosort}


\maketitle
\vspace{-2em}
\begin{abstract}
    In this paper, we consider a static, multi-period newsvendor model under a budget constraint. In  the case where the true demand distribution is known, we develop a heuristic algorithm to solve the problem. By comparing this algorithm with off-the-shelf solvers, we show that it generates near-optimal solutions in a short time. We then consider a scenario in which limited information on the demand distribution is available. It is assumed, however, that the true demand distribution lies within some given family of distributions and that samples can be obtained from it. We consider the cases of normal and Poisson demands. For each case, we show that using maximum likelihood estimates in place of the true parameters can lead to poor estimates of the true cost associated with an order quantity. Hence, we make use of likelihood inference to develop confidence sets for the true parameters. These are used as ambiguity sets in a distributionally robust model, where we enforce that the worst-case distribution lies in the same family as the true distribution. We solve these models by discretising the ambiguity set and reformulating them as piecewise linear models. We show that these models quickly become large as the ambiguity set grows, resulting in long computation times. To overcome this, we propose a heuristic cutting surface algorithm that exploits theoretical properties of the objective function to reduce the size of the ambiguity set. We illustrate that our cutting surface algorithm solves orders of magnitude faster than the piecewise linear model, while generating very near-optimal solutions.

\end{abstract}


\section{Introduction}

The newsvendor problem~\citep{Newsboy1} is a classical problem in operational
research and operations management. In this
problem, we consider a retailer facing the decision of how much stock to order in order to try to meet demand as closely as possible. Based on the fact that newspapers are no longer saleable after the period is over, in the newsvendor problem any unsold stock is either lost or sold at a lower price. Furthermore, the newsvendor is typically penalised for missing demand. Even in its earliest form, the newsvendor problem models a situation in which demand is not known exactly. However, in early papers such as that of~\cite{Newsboy1}, it was typically assumed that the demand distribution was known exactly. The resulting problem is therefore stochastic, and usually nonlinear. 

In the years after its introduction, the newsvendor problem has been extended in many ways. Two such extensions are relevant to this paper. The first is the static multi-period newsvendor problem~\citep{CHEN2017}. In this version of the problem, the retailer must decide prior to the selling horizon on how much stock to have delivered at the start of each period. This corresponds to cases where the retailer and the supplier must agree on their stocking quantities beforehand, and applies to many industries. As discussed by~\cite{CHEN2017}, this type of ordering structure is relevant to fresh produce industries. In these industries, due to the short shelf-lives of products, orders cannot simply be delivered prior to the selling period and stored for its entirety. Some products can only be stored for a short time and hence must be delivered within the selling horizon to ensure that it does not perish before it is sold. Furthermore, ordering in advance allows the retailer to use one period's order for multiple subsequent periods, if this yields a higher profit. Depending on the costs associated with holding stock, this can be cheaper than making multiple separate orders. Our model makes an additional extension to the model of~\cite{CHEN2017}, through a budget constraint. This is not common in multi-period models, due to the added complexity this constraint brings, but it is more common in multi-product models~\citep{ALFARES2005465}. This constraint adds an extra level of realism, since it is very unlikely that any retailer can spend infinite amounts of money. In this paper, we develop an iterative algorithm for our model.

Another key extension of the newsvendor problem is the \textit{distributionally robust} (DR) problem. In this version of the problem, we do not assume that the distribution of demand is known. We assume that the true distribution lies in some set, referred to as an \textit{ambiguity set}. Then the problem is to find the ordering quantity with the highest worst-case expected profit over all distributions in this ambiguity set.  Early DR problems assumed that only some moments of the distribution, such as mean and variance, are known and fixed. The ambiguity set then contains all distributions whose moments are equal to these values. We refer to problems with moment-based ambiguity sets as \textit{distribution free} (DF) problems~\citep{oneperiodrobustsol}. The assumption that these moments are known is not always realistic, and typically they must be estimated from sample data~\citep{DRNV_Wasserstein}. It has been found that poor estimation of these parameters can lead to poor estimation of costs~\citep{ROSSI2014674} and bias in solutions~\citep{Siegel21}. Furthermore, DF models can lead to overly conservative solutions~\citep{Likelihood_RNVP}. Due to these issues, many recent DR models considered ambiguity sets containing distributions that lie within some pre-prescribed distance of a nominal distribution~\citep{zhao2015data}. 

Our model is different from the standard DR model in two key ways. Firstly, we only consider parametric distributions. To be more specific, we assume that the demand distribution lies in some parametric family, and then construct ambiguity sets that only contain distributions that also lie in this family. Parametric distributions have often been used in the newsvendor problem due to their ability to allow statistical estimation and analysis to be integrated with optimisation~\citep{LIYANAGE2005341, ROSSI2014674}. Secondly, due to the pitfalls of poor cost estimation, we do not simply assume that the estimates of these parameters are the truth. Instead, we use maximum likelihood theory to develop confidence sets for the true parameters, and use these as ambiguity sets. The concept of parametric \textit{distributionally robust optimisation} (DRO) was introduced by~\citep{black}, in the context of a resource planning problem under binomial demand distributions. In this paper, we further develop this framework by studying a multi-period newsvendor problem under normal and Poisson demands. This allows us to show how the methods of~\cite{black} can be applied to both continuous and discrete demand distributions with infinite support.

The model resulting from the parametric ambiguity set does not have any tractable reformulation. Hence, we discretise the ambiguity set and represent the inner objective using a set of constraints. Due to the nonlinearity of the objective function, these constraints are non-linear. For discrete distributions, these constraints are piecewise linear, and for continuous distributions we use piecewise linear approximations of them. One piecewise linear constraint is required for each distribution in the ambiguity set and each period considered by the model. Hence, the model can grow very large and become slow to solve. Therefore, inspired by~\cite{black}, we develop a heuristic \textit{cutting surface} (CS) algorithm in order to solve the model in a fast time. CS algorithms have often been used for solving non-parametric DRO problems to $\varepsilon$-optimality~\citep{mehrotra2014cutting}. Such CS algorithms were shown to become slow for large parametric problems by~\citep{black}. Hence, our algorithm is a heuristic CS algorithm that exploits knowledge of the distributional family to improve solution times. In general, a CS algorithm iteratively solves the model over subsets of the ambiguity set, each time adding a new distribution to the subset before solving again. This new distribution is chosen by finding the worst-case distribution for the most recent solution. To speed up this stage, we use theoretical properties of the objective function to develop a set of extreme distributions from which to select the worst-case. We will show that our heuristic CS algorithm performs very well and generates solutions in a short time.

In summary, the key contributions of our paper are as follows:
\begin{enumerate}
    \item We extend the multi-period model of~\cite{CHEN2017} by including a budget constraint. 
    \item We present an iterative solution algorithm for the multi-period model under a budget constraint. As far as we know, there is no existing algorithm for this problem.
    \item We show that the MLE approach provides poor estimates of the costs of a given order quantity, which can lead to predicting a profit for an order that would result in a loss. 
    \item We extend the concept of parametric DRO from earlier research~\citep{black} to the newsvendor literature. We consider both continuous and discrete random variables with infinite support. Only binomial random variables have been studied before. In addition, this paper shows how the work of~\cite{black} can be extended to models with nonlinear objective functions.
    \item We develop a fast heuristic version of the CS algorithm from the DRO literature. The general CS algorithm is $\varepsilon$-optimal, but can solve slowly for large ambiguity sets. Hence, our heuristic CS algorithm exploits properties of the distributional family to reduce the size of the ambiguity set considered, at little cost to solution quality. We perform extensive computational experiments to test the efficacy of our CS algorithm. 
\end{enumerate}
\section{Literature review}

Our research is based on three areas of the newsvendor literature: DR newsvendor problems, multi-period  and capacitated/budgeted newsvendor problems, and parametric newsvendor problems. We now review the literature corresponding to these three variations of the problem.

\subsection{Distributionally robust newsvendor problems}

The first example of the DR newsvendor problem comes from~\cite{oneperiodrobustsol}. Scarf formulated the first DF model under the assumption that only its mean and variance are known. He optimises the worst-case cost over all distributions with this mean and variance, proving a simple ordering rule to be optimal. The DF model has been the subject of many papers since the work of~\cite{oneperiodrobustsol}. For example, \cite{gallego_extensions} extended this work to the multi-product and random yield models. \cite{DF_balking} also extended the single-period DF model to the case where customers can balk, i.e.\ decide not to buy an item, after observing stock levels. Further extensions to the DF model come from~\cite{ALFARES2005465}, who extended the work of \cite{gallego_extensions} to the case where the model includes a shortage cost for missing demand. \cite{OUYANG20021} extended the DF model to incorporate uncertainty in parameters such as lost-sales and backorder rates. Later extensions of the model incorporate addition elements such as the cost of advertising~\citep{Lee11Advertising}, risk- and ambiguity-aversion~\citep{AverseNV}, and carbon emissions~\citep{LiuCarbon, BAI20162}. These authors were able to find closed-form solutions using derivatives and/or Karush-Kuhn-Tucker (KKT) conditions.

DF models are commonly used due to their tractability. However, the assumptions are sometimes unrealistic and can lead to issues with the resulting solutions. \cite{DRNV_Wasserstein} discuss these issues, citing that the main reason for not using a DF approach is that it is unrealistic to assume that any parameters of the distribution are known exactly. In addition, the moment-based ambiguity sets used in DF models can lead to overly conservative decisions~\citep{Likelihood_RNVP}. In practice, typically these parameters must be estimated from historical data, and~\cite{DRNV_Wasserstein} state that maximum likelihood estimation approaches can lead to suboptimal solutions. An example of this comes from~\cite{ROSSI2014674}, who used confidence interval analysis for the newsvendor problem under parametric distributions. Instead of singular solutions, their method provides a discrete set of order quantities that contains the true optimal order quantity with a given probability. They found that, even for large samples, the costs given by the maximum likelihood estimates were inaccurate with respect to the true cost. As far as we are aware, this is the only method so far that allows decision makers to hedge against uncertainty in parameters. However, it may not be practical, since the decision maker then faces the dilemma of which solution to choose. Since it provides a single order quantity to hedge against the worst costs, a DR solution is therefore more convenient for a risk-averse decision maker. 

In more recent literature, DR models with more complex ambiguity sets have been developed. Where the DF model's ambiguity set typically contains all distributions with some of their parameters fixed at nominal values, many models have used distance-based ambiguity sets. These sets contain all distributions that lie within some predefined distance from the nominal, according to some distance measure. These are common in the DRO literature, with an example distance measure being $\phi$-divergences~\citep{bayraksan2015data}. These ambiguity sets lead to problems requiring different treatment than DF problems, and are usually solved via tractable reformulations using duality.  One of the first examples of distance-based ambiguity sets for the newsvendor problem was given by~\cite{zhao2015data}. They used $\zeta$-structure probability metrics, a class of metrics containing the Wasserstein distance, to form their ambiguity sets. Their model was reformulated as a stochastic program, and solved by an iterative sampling algorithm. \cite{gao2017distributionally} studied a single-period newsvendor problem with Wasserstein and $\phi$-divergence ambiguity sets, which was solved using a convex programming reformulation. The Wasserstein distance is also studied by~\cite{EsfahaniWasserstein} and~\cite{DRNV_Wasserstein}, who also find that it leads to convex and linear programming reformulations. 

Our research differs from the above literature in three key ways. Firstly, we do not assume that any parameters of the distribution are known. In order to avoid the issues resulting from poor estimation, we instead assume that the demand distribution lies in a known parametric family, but that its parameters are completely unknown. Secondly, we enforce that the worst-case distribution selected by the DRO model also lies within this family, in order to ensure that it is a reasonable candidate for the true distribution. We enforce this via utilising the distribution's probability mass or density function in the objective function directly, and identifying the worst-case parameters instead of the worst-case distribution itself. In particular, we consider normal and Poisson demands. Hence, our model considers nonlinear objective functions, continuous distributions and discrete distributions with infinite support. This extends the work of~\cite{black}, where only binomial distributions and linear objective functions were studied. As was true for the resource planning problem of~\cite{black}, parametric ambiguity sets do not lead to tractable reformulations. Therefore, we solve using discretised ambiguity sets. This means that the inner objective can be represented by a finite set of constraints. Since the number of constraints required can be very large, we use theoretical properties of the cost function to develop a fast heuristic CS algorithm for our problem. We show that it performs very well, and that it generates solutions in a much shorter time than solving the full model. The final way in which our research differs is that we do not use distance measures to build our ambiguity sets. One of the key benefits of doing this, as stated by~\cite{DRNV_Wasserstein}, is that they allow the ambiguity sets to provide asymptotic probabilistic guarantees. In our methodology, we use maximum likelihood theory to develop confidence sets for the true parameters of the demand distribution. This allows us to provide ambiguity sets with the same probabilistic guarantees, without requiring the use of distance measures. 

\subsection{Multi-period and capacitated/budgeted newsvendor problems}

A natural extension of the classical newsvendor problem of~\cite{Newsboy1} is the multi-period newsvendor problem. This incorporates either the option to make orders at multiple periods in a planning horizon (dynamic problem), or to make orders for multiple periods in a selling horizon prior to the horizon concerned (static). The earliest multi-period problems are dynamic, where orders are made and delivered in each period. These problems are typically treated as dynamic programs, and can often be solved by base-stock policies~\citep{bellmaninv, Bouakiz92}. Under a base-stock policy, the order quantity for each period is chosen to bring the stock level up to a predefined level.  Much of the literature on multi-period newsvendor problems has considered the dynamic formulation. Some examples include~\cite{AHMED2007226}, who studied the case where the objective function is a coherent risk measure, \cite{Retsef07}, who considered that case where the demand distribution is unknown, and \cite{Altintas08}, who considered setting discounts for a multi-period newsvendor. Later examples considered service-dependent demand~\citep{Deng14}, non-stationary demand~\citep{KIM2015139}, and two-product models where the total demand is fixed~\citep{ZHANG2016143}. 

Rather than dynamic, our model is static. Static multi-period models are less common in the literature. \cite{MATSUYAMA2006170} considered a static problem with ordering cycles, but mainly focused on deriving theoretical properties of the optimal profit. Other papers on static multi-period problems typically do not include any constraints on the order quantities, because adding these can result in intractable models. \cite{CHEN2017} consider a situation in which a seller sets their price for selling to a multi-period, static newsvendor (buyer). They formulate the buyer's problem as a nonlinear program and solve using the KKT conditions. This solution would be complicated by additional constraints, as we will show later. \cite{Mehran19} also consider a price-setting scenario, this time through a distribution-free approach under price-dependent demand. This model is unconstrained, and solved using first-order conditions.  

Our model extends that of~\cite{CHEN2017} to the case where the buyer has a monetary budget. This adds a constraint to the model, meaning that standard KKT solutions are not easily applicable. To the best of our knowledge, there are no optimal algorithms in the literature for this problem. In fact, typically, multi-period models do not have capacity or budget constraints. However, such constraints are commonly found in the multi-product literature. In this literature, iterative algorithms are usually used~\citep{Lau1996TheNP, ABDELMALEK2005296, ALFARES2005465}. This is mainly due to the presence of non-negativity constraints in addition to the budget or capacity constraint. Adding Lagrange multipliers for each non-negativity constraint can make the KKT conditions difficult to solve. Furthermore, ignoring them can lead to negative order quantities~\cite{Lau1996TheNP}. Motivated by this, \cite{ABDELMALEK2005296} developed an iterative algorithm for the multi-product capacitated problem, based on Lagrangian relaxation. Such algorithms are now common ways to solve capacitated problems. To solve our multi-period problem, we adapt another iterative algorithm that was presented by~\cite{ALFARES2005465}. This algorithm is designed to solve a multi-product problem with a capacity constraint. It consists of first solving the KKT stationarity condition with no budget constraint and no non-negativity constraint, and then iteratively increasing the Lagrange multiplier from zero until the solution either becomes feasible or negative. If the solution becomes negative, then the corresponding order quantity is set to zero and it is no longer considered. The process then begins again. This allows us to solve the problem without ever explicitly incorporating the non-negativity constraints.

In summary, our research differs from the multi-period literature due to the presence of the budget constraint. Due to the additional complexity caused by this constraint, we use the KKT conditions of the problem to develop an iterative solution algorithm. Similar to the algorithm of~\cite{ALFARES2005465} for the multi-product capacitated problem, our algorithm iteratively increases the Lagrange multiplier for the budget constraint until the order quantities either become feasible or negative. If the latter case occurs first, one order is set to zero and removed from further consideration. However, compared to the multi-product model, the multi-period model has the additional complexity that the optimal unconstrained order quantities depend on one another. This makes it unclear which period's order should be set to zero when one becomes negative. Due to this, the proposed algorithm is not optimal in general for our problem. However, we will show that it performs comparably with algorithms from off-the-shelf solvers, and is even faster than these algorithms. More details can be found in Section~\ref{sec:fixed_dist}, where we present and test the algorithm.

\subsection{Parametric newsvendor problems}

The final area of the literature that we will review concerns parametric newsvendor problems. Parametric demand distributions allow for two main methodologies for the newsvendor problem. The first approach is the Bayesian approach. This consists of first assuming some prior distribution of demand, and then updating the prior distribution to obtain the posterior after some demand is realised. For example, \cite{HILL1997555} and~\cite{BensoussanBayes} used Bayesian approaches for estimating newsvendor demand under uniform and exponential priors, respectively, in the case of fully observable demand. The case of censored demand has also been addressed in a Bayesian fashion. For example, \cite{ChenUnobserved} considered a dynamic inventory problem with unobserved data, focussing on developing heuristics and bounds based on Bayesian updating. \cite{Merserau} uses Bayesian estimation to study the effect of inventory record inaccuracy under censored demand. As discussed by~\cite{ROSSI2014674}, the Bayesian approach can be inappropriate due to difficulty in selecting a prior and infinite iterations required to prove convergence. 

The second approach is the frequentist approach. This corresponds to authors directly estimating the parameters of the demand distribution from sample data. Early papers on this approach did not consider integrating the estimates into the model. For example, \cite{Nahmias1994739} studied different estimation methods for the parameters of a normal demand distribution in the situation where lost sales are not observed. \cite{Agrawal96} studied estimating negative binomial demand from sample data. Both of these studies assumed the order quantity was given. Later papers studied the problem of incorporating frequentist estimates into the optimisation models. \cite{LIYANAGE2005341} showed that simply using parameter estimates in place of true parameters can lead to suboptimal solutions. They then introduced the ``operational statistics'' framework, which integrates estimation and optimisation to reduce bias. There has also been research on which distributions should be used for newsvendor demand. \cite{GALLEGO2007} compared the profits obtained from normal, lognormal, gamma and negative binomial distributions in the case where the coefficient of variation was high. They find that the normal distribution can lead to large negative profits, and suggest using non-negative distributions instead. \cite{ROSSI2014674} propose a methodology that combines statistical confidence interval analysis with the newsvendor problem to reduce the effect of poor parameter estimation. They provide confidence intervals for the true parameters of the distributions and use this to provide confidence intervals for the true optimal cost. They argue that this avoids the pitfalls of simply using maximum likelihood estimates, namely poor estimation of optimal costs. \cite{Siegel21} also present a methodology to avoid said pitfalls, which is based on an asymptotic adjustment to the estimates in order to reduce bias.

Our methodology utilises maximum likelihood theory for parametric distributions in the newsvendor problem. Motivated by the problems that can occur from assuming that they are the true parameters, we study the effects of assuming that parameter estimates are the true values. We will show that this method performs similarly to how it did for~\cite{ROSSI2014674}, namely that it leads to suboptimality and poor cost estimates. In addition, we use multivariate confidence theory to develop confidence sets for the true parameters. These are then used as ambiguity sets in our models. Although our confidence sets are based on asymptotic theory, they are large for small samples and small for large samples. This means that the model actively accounts for small sample sizes through constructing a more conservative set.
\section{Fixed distribution model and solution}

In Section~\ref{sec:fd_model}, we present the model in the case where the demand distribution is known. Then, in Section~\ref{sec:fixed_dist}, we present our algorithm for solving this problem. This will be used to benchmark the DRO methods that are the main focus of the paper.

\subsection{Model}\label{sec:fd_model}

Adapting the model by~\cite{CHEN2017}, we assume
that there are $T$ periods for which an order needs to be made. We denote the periods by $t
\in \m{T} = \{1,\dots, T\}$. The decision vector is denoted $\bm{q} = (q_1, \dots, q_T)$, and gives the amount of stock to order for each period $t \in \m{T}$. The demand for period $t$ is denoted by $X_t$.
We also assume that the newsvendor has a total budget  
for expenditure of $W$, that limits the amount that can be ordered in total.
In practice, we do not have access to infinite resources and hence this is a
realistic addition to the model. We denote by $F$ the joint cumulative distribution function (CDF) of the demand vector $\bm{X}$.
Furthermore,  let
$I_t$ be the inventory available during period $t \in \{1,\dots,T\}$ and assume $I_0 = 0$ w.l.o.g. Finally,
let us define the costs for the model. Let $w_t$ be the cost of
ordering a unit of stock to sell on period $t$. We will assume that $w_t \ge w_{t+1} \fa t=1,\dots,T$. This is reasonable because, as shown by~\cite{CHEN2017}, any optimal set of prices for the supplier satisfies this condition. Let $h$ be the cost of holding one
unit of stock for one period (i.e.\ a \textit{holding cost}), and let $b$
be the cost of one unit of unmet demand (\textit{backorder cost}). Finally, let $p$ be the price a customer pays to purchase stock from the newsvendor. Then, in each period $t$ the following events occur:
\begin{enumerate}
    \item The order quantity $q_t$ arrives and the supplier pays $w_t q_t$ for the shipment.
    \item Demand $X_t$ is realised.
    \item Inventory at end of period is calculated as $I_t = I_{t-1} + q_t - X_t$.
    \item If $t \le T - 1$ then the newsvendor receives profit $p X_t - h I^+_t - b I^-_t$. Here $I^+_t = \max\{I_t, 0\}$ is the leftover product that must be held for the following period and $I^-_t = \max\{-I_t, 0\}$ is the unmet demand. This profit is the revenue from sales, minus the total holding costs and backorder costs. It assumes that all demand for period $t$ is met at some point during or after period $t$. If it cannot be met, this is not known until $t=T$, where the profit is reduced accordingly. If $t = T$ then the profit is $p\min\{X_T, I_{T-1} + q_T\} - hI^+_T - bI^-_T$.
    The first term in this expression represents that demand in period $T$ can only be met in period $T$, since there are no subsequent periods. Hence, it can only be met if the newsvendor has enough stock to meet it in period $T$.
\end{enumerate}

Given this, the expected cost for an ordering decision $\bm{q}$ is given by:
\begin{equation*}
    C_{F}(\bm{q}) = 
    \sum_{t=1}^T \l(h\E_{F}[I^+_t] + b\E_{F}[I^-_t] + w_t q_t\r)
    - p\l(\E_F\l[\sum_{t=1}^T X_t\r] - \E_F[I^-_T]\r),
\end{equation*}
which represents the total expected costs incurred from purchasing, holding and backordering, minus any revenue made. The final term is the expected revenue made from sales, where $\E_F\l[\sum_{t=1}^T X_t\r]$ is the expected total demand and $\E_F[I^-_T]$ is the expected number of missed sales, which corresponds to the unmet demand in period $T$. Then, rewriting slightly, the model for a fixed distribution $F$, referred to as  MPNVP, is given by:
\begin{align}
    \min_{\bm{q}}C_{F}(\bm{q}) &= p\E_{F}[I^-_T] + 
    \sum_{t=1}^T \l(h\E_{F}[I^+_t] + b\E_{F}[I^-_t] + w_t q_t
    - p \E_F[X_t]\r) \label{eq:MPNVP_obj}\\
    \text{s.t. } &\sum_{t=1}^T w_t q_t \le W,\label{eq:MPNVP_1} \\
    &I_t = I_{t-1} + q_t - X_t \fa t = 1, \dots, T,\label{eq:MPNVP_2} \\
    &q_t \ge 0 \fa t = 1,\dots,T.\label{eq:MPNVP_3}
\end{align}
Here, constraint~\eqref{eq:MPNVP_1} is the budget constraint and~\eqref{eq:MPNVP_2} is the inventory constraint.

\subsection{Solution procedure}\label{sec:fixed_dist}

\subsubsection{A solution of the KKT stationarity condition}

In order to develop a solution algorithm for this problem, we will now assume that $F$ represents a continuous distribution. This allows us to use the KKT conditions to derive an algorithm. The resulting algorithm iteratively increases the Lagrange multiplier for the budget constraint, each time generating a new solution by evaluating the inverse CDF of the demand distribution. Since the inverse CDF of a discrete distribution will always provide integer order quantities, the algorithm we develop can also be applied to discrete distributions. We will show that it performs well for both the Poisson and normal distributions. In particular, we will use the following result:
\begin{theorem}\label{thm:KKT_sol}
    For a given Lagrange multiplier $\nu$, suppose that the following two conditions hold:
    \begin{enumerate}
        \item $\nu \le \frac{b}{w_t - w_{t+1}} - 1 \fa t \in \{1,\dots,T-1\}$,
        \item $\nu  \le \frac{b+p}{w_T} - 1$.
    \end{enumerate}
    Then, the solution:
    \begin{equation}\label{eq:KKT_sol}
    \small
    \l.\begin{aligned}
            q^*_1 &= \tilde{F}^{-1}_1\l(\frac{b - (1+\nu)(w_1 - w_2)}{h+b}\r)\\
            q^*_t &= \tilde{F}^{-1}_t\l(\frac{b - (1+\nu)(w_t - w_{t+1})}{h+b}\r) - \tilde{F}^{-1}_{t-1}\l(\frac{b - (1+\nu)(w_{t-1} - w_{t})}{h+b}\r),  (2\le t \le T-1) \\
            q^*_T &=  \tilde{F}^{-1}_T\l(\frac{b - (1+\nu)w_T + p }{h+b+p}\r) - \tilde{F}^{-1}_{T-1}\l(\frac{b - (1+\nu)(w_{T-1} - w_{T})}{{h+b}}\r).
        \end{aligned}
        \r\}
    \end{equation}
    satisfies the stationarity KKT condition of MPNVP.
\end{theorem}
This theorem is proved in Appendix~A. Even though this solution satisfies the stationarity condition of MPNVP, it may not be feasible. We also require that $q_t \ge 0 \fa t$ and $\sum_{t=1}^T w_t q_t \le W$. Adding Lagrange multipliers for every non-negativity constraint results in an intractable set of KKT conditions, as discussed by~\cite{ALFARES2005465, ABDELMALEK2005296} and~\cite{Lau1996TheNP}, who each studied multi-product newsvendor models with budget constraints. Due to this, along with the fact that relaxing non-negativity constraints can lead to negative order quantities, the common approach among these papers has been to develop iterative algorithms that do not require that the non-negativity constraints are enforced directly. 

\subsubsection{An iterative solution algorithm: FD}

We now present our algorithm for solving MPNVP. We will refer to this algorithm as \textit{fixed distribution} (FD). Details on deriving FD can be found in Appendix~B.1. Let us denote by $\m{T}^0$ the set of days that have had their orders set to zero and been removed from consideration. Let $\bm{q}^*(\nu)$ be the solution obtained from evaluating~\eqref{eq:KKT_sol} with Lagrange multiplier $\nu$. Similarly, define $\tilde{\bm{q}}^*(\nu)$ as the solution obtained from evaluating~\eqref{eq:KKT_sol_adj} with Lagrange multiplier $\nu$.
\begin{equation}\label{eq:KKT_sol_adj}
\l.
\begin{aligned}
    \tilde{q}^*_1 &= \begin{cases}
        \tilde{F}^{-1}_1\l(\frac{b - (1+\nu)(w_1 - w_2)}{h+b}\r) &\text{ if } 1 \notin \m{T}^0\\
        0 &\text{ if } 1 \in \m{T}^0
    \end{cases}\\
    \tilde{q}^*_t &= \begin{cases}
         \tilde{F}^{-1}_t\l(\frac{b - (1+\nu)(w_t - w_{t+1})}{h+b}\r) - \sum_{k=1}^{t-1} \tilde{q}^*_k &\text{ if } t \notin \m{T}^0\\
        0 &\text{ if } t \in \m{T}^0
    \end{cases}
    \fa t \in \{2,\dots,T-1\}\\
    \tilde{q}^*_T &= \begin{cases}
        \tilde{F}^{-1}_T\l(\frac{b - (1+\nu)w_T + p }{h+b+p}\r) - \sum_{k=1}^{T-1}\tilde{q}^*_k  &\text{ if } T \notin \m{T}^0\\
        0 &\text{ if } T \in \m{T}^0
    \end{cases}.
\end{aligned}
\r\}
\end{equation}
This solution adjusts the orders to account for the fact the days in $\m{T}^0$ have been set to zero, to ensure the stationarity condition is met. Then, FD can be described as follows: 
\begin{enumerate}
    \item Set $\m{T}^0 = \{\}$.
    \item Take $\bm{q} = (\max\{\tilde{q}^*_1(0), 0\},\dots,\max\{\tilde{q}^*_T(0), 0\})$, i.e.\ $\tilde{\bm{q}}^*(0)$ but setting negative values to zero.
    \item If $\sum_{t=1}^T w_t q_t \le W$ then go to step 7. If $\sum_{t=1}^T w_t q_t > W$ then if $\m{T}^0 = \{\}$  go to step 4 and otherwise go to step 5. 
    \item The budget constraint is binding. Compute upper bounds on $\nu$, i.e.\ 
    \begin{equation*}
        \nu^{\text{UBs}} = \l\{\frac{b}{w_t - w_{t+1}} - 1 \fa t \in \{1,\dots,T-1\}\r\} \cup \l\{\frac{b+p}{w_T} - 1\r\}
    \end{equation*}
    and set $\nu^{\text{UB}} = \min\{\nu^{\text{UBs}}\}$. Then $\nu \le \nu^{\text{UB}}$ ensures that the inverse CDFs are all defined at $\nu$. Go to step 5.
    \item Starting from $\nu=0$, increase $\nu$ until the first occurrence of one of the following:
    \begin{enumerate}
        \item The resulting $\tilde{\bm{q}}^*(\nu)$ has $\tilde{q}^*_k(\nu) < 0$ for some $k$. In this case define $\nu^* = \nu$ and go to step 6.
        \item The resulting $\tilde{\bm{q}}^*(\nu)$ has $\sum_{t=1}^T w_t \tilde{q}^*_t(\nu) \le W$. In this case go to step 7.
    \end{enumerate}
    \item Select which day to add to $\m{T}^0$:
    \begin{enumerate}
        \item For each day $t \in \m{T} \setminus \m{T}^0$, calculate a lower bound on the cost resulting from setting $q_t = 0$ and adding $t$ to $\m{T}^0$. Specifically, calculate $\bm{\tilde{q}}^*(\nu^*)$ with $\m{T}^0 = \m{T}^0 \cup \{t\}$.
        \item Select $k$ to be the day with the lowest lower bound resulting from adding it to $\m{T}^0$.
        \item Set $\m{T}^0 = \m{T}^0 \cup \{k\}$ and go to step 2.
    \end{enumerate}  
    \item Return $\tilde{\bm{q}}^*(\nu)$.
\end{enumerate}
The algorithm uses a line search on $[0, \nu^{\text{UB}}]$ to find the optimal value of $\nu$ in Step 5. More details on step 5 can be found in Appendix~B.2. The reason why FD is only a heuristic is step 6. The idea behind step 6 is that, unlike in the multi-product model, the order quantity that becomes negative may not be the one that should be set to zero. This is because, in the multi-period model, setting a negative order to zero reduces subsequent orders.

In order to select which order quantity should be set to zero, we estimate the effect of each choice by finding a lower bound on the resulting cost. For each $t \in \m{T} \setminus \m{T}^0$, we calculate this lower bound by temporarily adding $t$ to $\m{T}^0$ and evaluating the cost of the resulting $\bm{\tilde{q}}^*(\nu)$. The reason for evaluating the next decision at $\nu^*$ is as follows. Since step 5(a) caused us to proceed to step 6, we know that a negative order occurs before the budget constraint is met. Since period $t$'s order was negative, if this order is set to zero then the budget constraint will still not be met (otherwise 5(b) would have occurred first). Hence, the remaining orders must be reduced further, i.e.\ we must have $\nu \ge \nu^*$ if  $\bm{\tilde{q}}^*(\nu)$ is to meet the budget constraint. Since $\nu$ has to be increased further to enforce the budget constraint, the cost at optimality must be higher than the cost of $\bm{\tilde{q}}^*(\nu)$. Hence, the cost of this solution gives a lower bound on the cost of the solution where $t$ has been added to $\m{T}^0$. Given all lower bounds, we then select $t$ with the lowest lower bound and add it to $\m{T}^0$.
 
\subsubsection{Testing FD}\label{sec:FD_results}

In order to establish the performance of FD, we tested it on 1536 problem instances. We obtained these by considering problems with $T \in \{2,\dots,5\}$, $p, h, b \in \{1, 2\}$, $W \in \{10, 25, 50\}$, $w \in \{(2T, 2(T-1),\dots, 2), (1, \frac{1}{2},\dots,\frac{1}{T})\}$. In addition, we selected 2 mean vectors randomly for each $T$ and set the standard deviation vector to be $\bm{\sigma} = \frac{\bm{\mu}}{4}$ in each case. This only affects the normal experiments. Furthermore, the tolerance for the budget constraint was either 0 or $10^{-6}$. The tuning parameter $\tau$ for the line search used by FD was tested for $\tau \in \{0.25, 0.5\}$. This generated 768 instances for normal and 768 instances for Poisson demands. Since our algorithm is a heuristic, we compare it with 3 off-the-shelf algorithms that might be applied to solve MPNVP. The first two algorithms are Sequential Least Squares Quadratic Programming (SLQP)~\citep{SLSQP} and Trust Constraint (TC)~\cite{trust_constr} from the Python library Scipy~\citep{2020SciPy-NMeth}. For more details on these algorithms, see Appendix~B.3. The final algorithm is Piecewise Linear Approximation (PLA), an algorithm that uses Gurobi~\citep{gurobi} to solve a piecewise linear approximation of MPNVP.

\begin{figure}[htbp!]
    \centering
    \begin{subfigure}[t]{0.45\textwidth}
        \centering
        \includegraphics[width=0.95\textwidth]{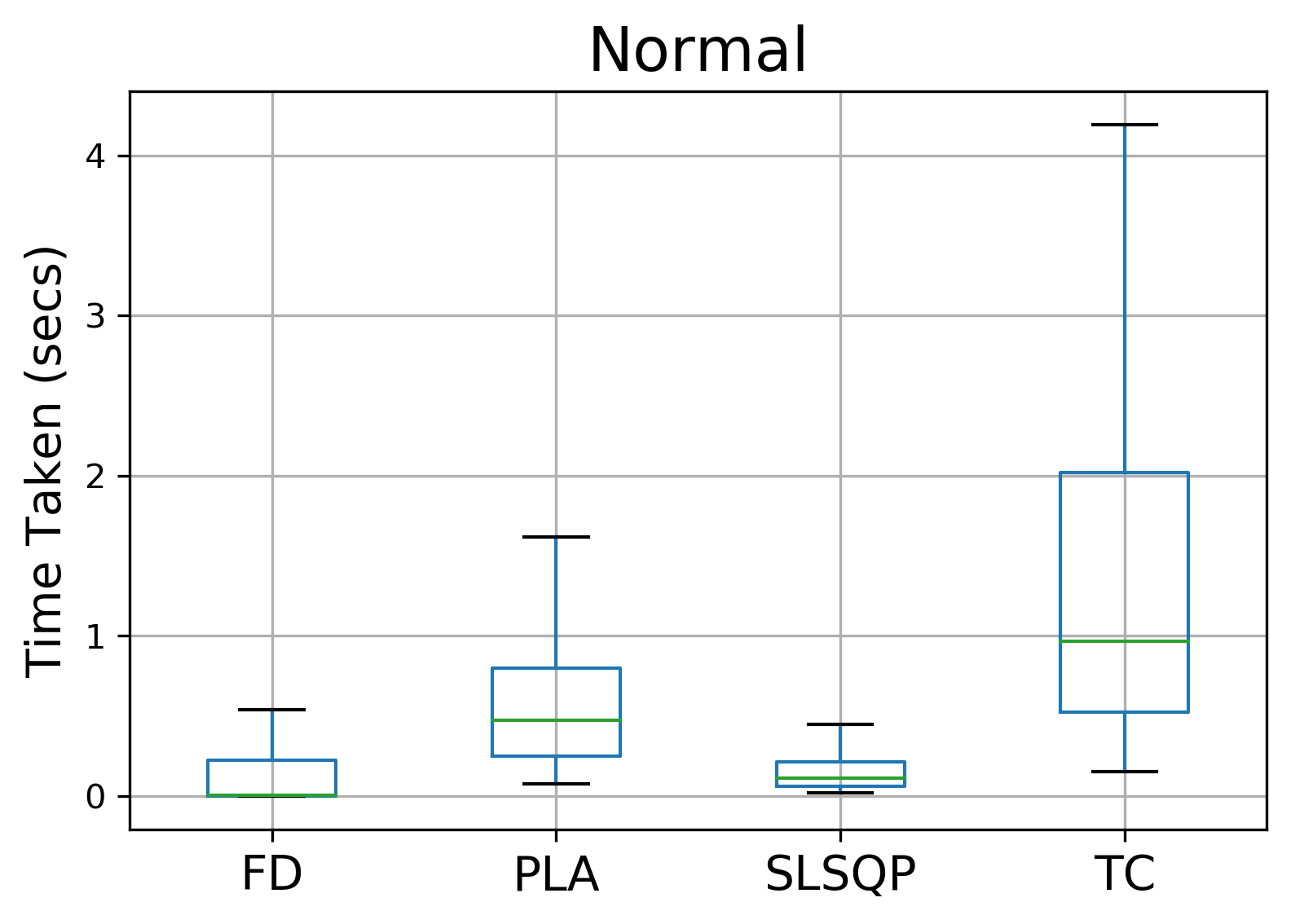}
        \caption{}
        \label{fig:FD_normal_times}
    \end{subfigure}
    \begin{subfigure}[t]{0.45\textwidth}
        \centering
        \includegraphics[width=0.95\textwidth]{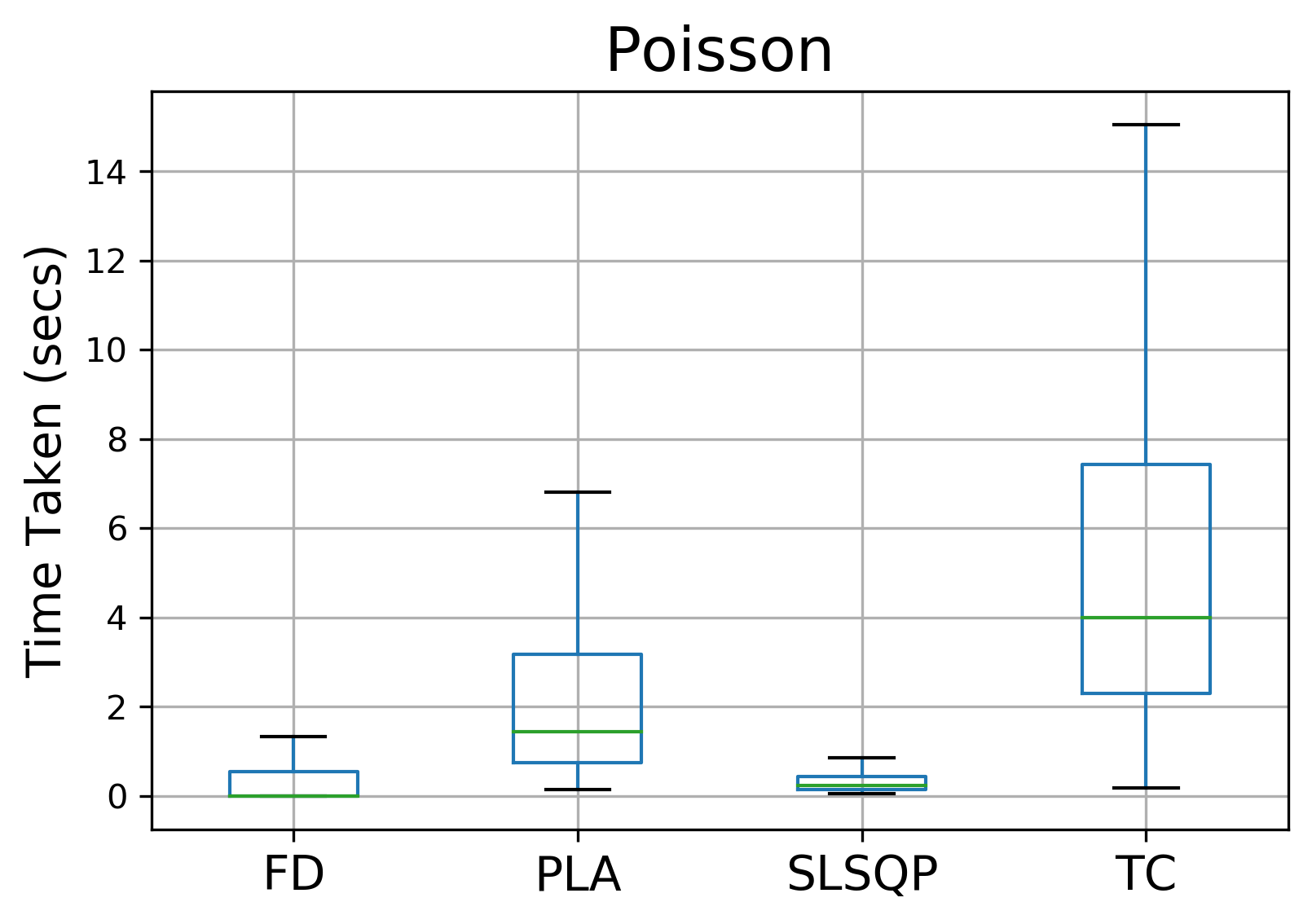}
        \caption{}
        \label{fig:FD_poisson_times}
    \end{subfigure}
    \caption{Boxplots summarising times taken by the 4 algorithms}
    \label{fig:FD_times}
\end{figure}

We now summarise the results of our experiments. Figure~\ref{fig:FD_times} shows the times taken by each algorithm to solve MPNVP under normal and Poisson demands, with outliers removed. Points that are more than 1.5 times the IQR above the 75th percentile or below the 25th percentile are treated as outliers. Note that this is the case for all boxplots in this paper, and it is the standard metric for defining outliers in boxplots. The plots show that FD and SLSQP are the fastest algorithms for both normal and Poisson demands. On average, for normal demands, FD took 0.14 seconds and SLSQP took 0.15 seconds. For Poisson demands, FD took 0.23 seconds whereas SLSQP took 0.32 seconds, on average. Hence, our results suggest that FD is the fastest algorithm on average.

\begin{figure}[htbp!]
    \centering
    \begin{subfigure}[t]{0.45\textwidth}
        \centering
        \includegraphics[width=0.95\textwidth]{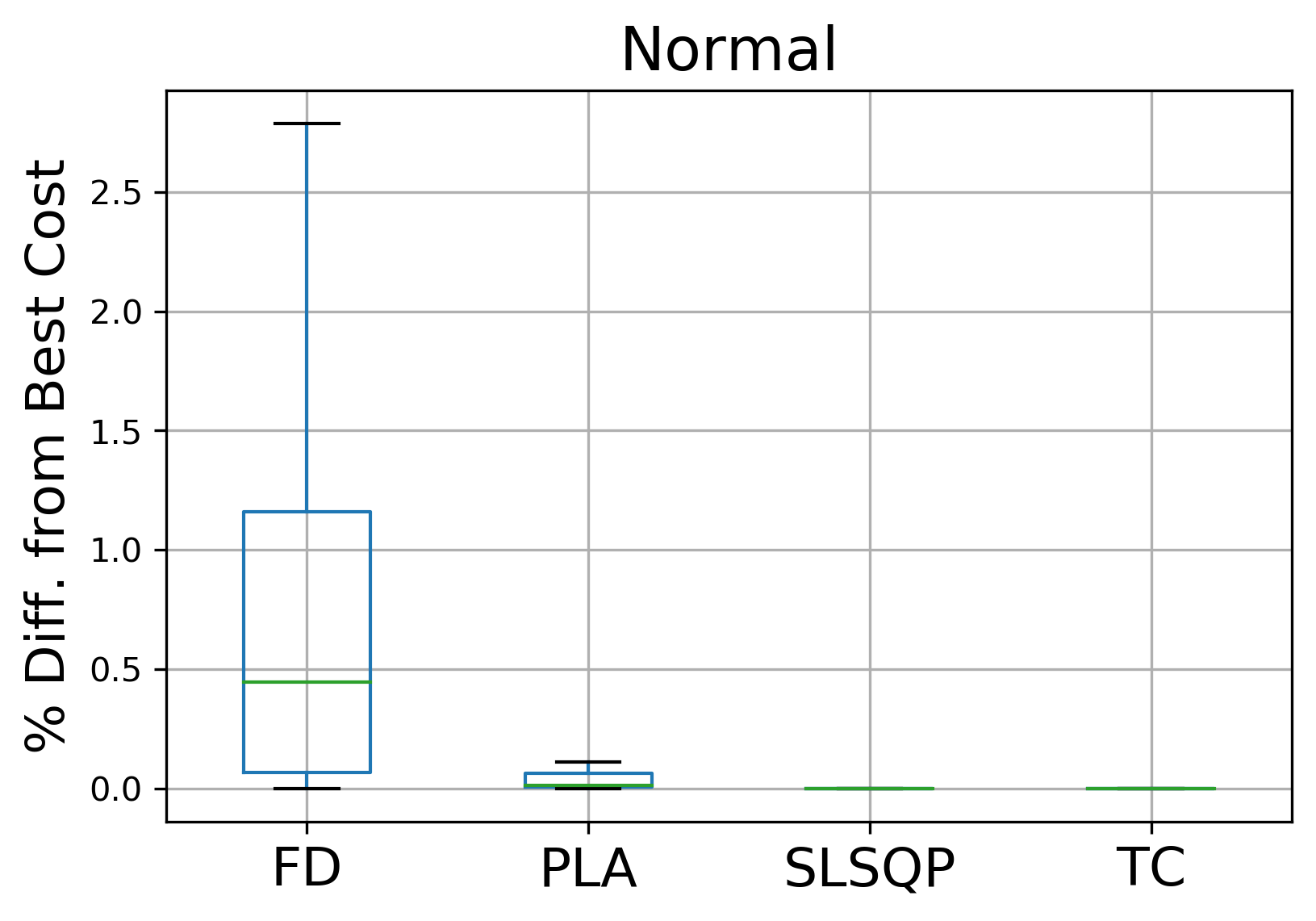}
        \caption{}
        \label{fig:FD_normal_gaps}
    \end{subfigure}
    \begin{subfigure}[t]{0.45\textwidth}
        \centering
        \includegraphics[width=0.95\textwidth]{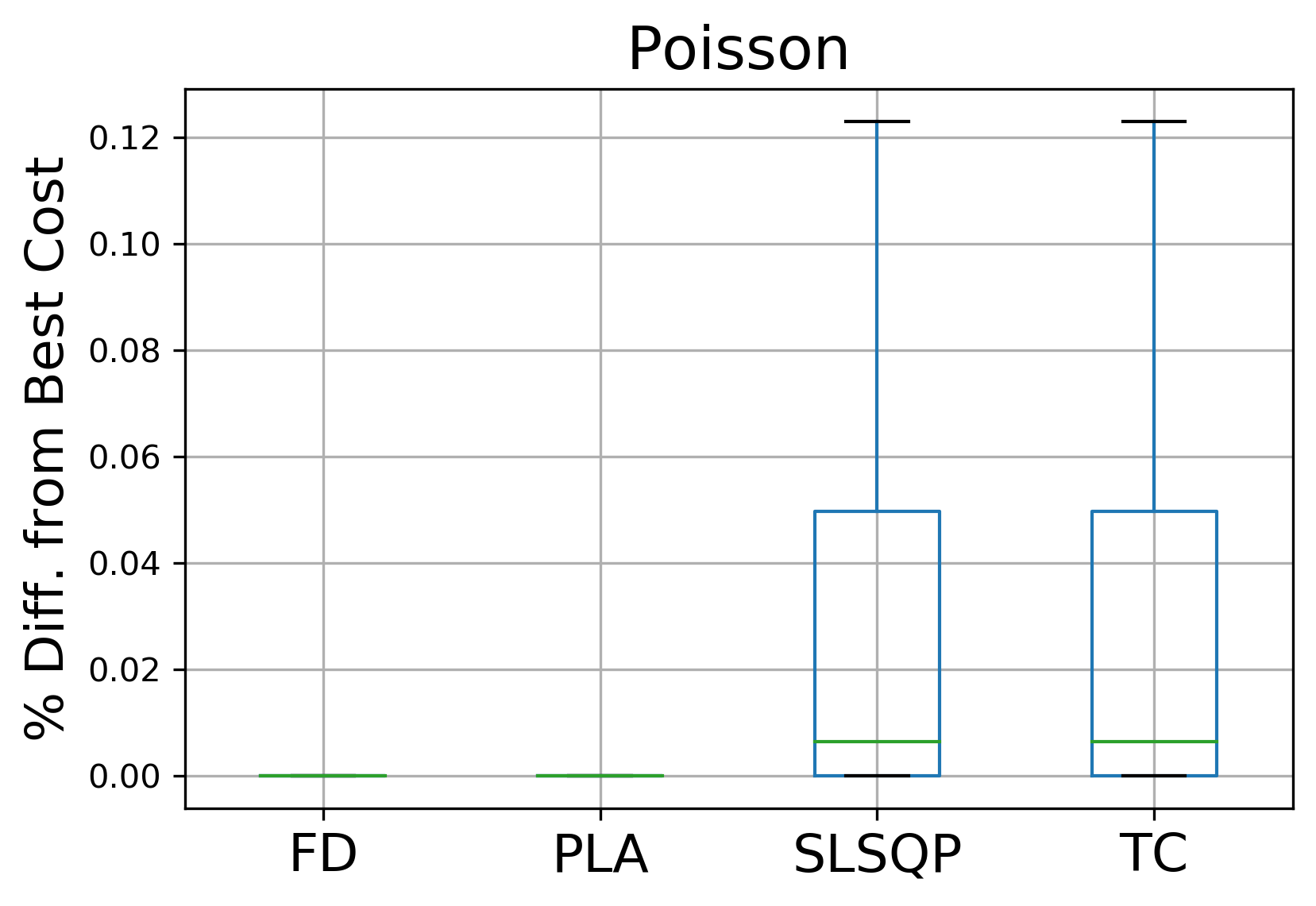}
        \caption{}
        \label{fig:FD_poisson_gaps}
    \end{subfigure}
    \caption{Boxplots summarising percentage optimality gaps of the 4 algorithms}
    \label{fig:optgaps_alf}
\end{figure}

To understand the performance of the algorithms, we define a percentage gap from the best solutions as follows:
$$100 \times \frac{C_{F}(\bm{q}^a) - \min_{y \in \m{Y}} C_F(\bm{q}^y)}{\min_{y \in \m{Y}} C_F(\bm{q}^y)},$$
where $a, y \in \m{Y} = \{\text{FD}, \text{PLA}, \text{SLSQP}, \text{TC}\}$ represents an algorithm, and $\bm{q}^y$ is the solution returned by algorithm $y$. Given this definition, we summarise the gaps in Figure~\ref{fig:optgaps_alf}. This plot shows that FD is typically within 2.5\% of the best solution for normal demands. SLSQP and TC perform the best here. However, for Poisson demands, FD and PLA are the best algorithms in terms of objective value. FD had the best objective value in 84\% of Poisson instances. On average over all instances, PLA is the best performing algorithm. It is important to mention that, while SLSQP performs well for Poisson demands, in 678 out of 768 instances with Poisson demands, its solution was not integer. There is no way to enforce integer variables within this algorithm. Hence, FD and PLA have this as an advantage over SLSQP: if the distribution is discrete, their solutions are integer.

The final result we would like to mention is the exceedance of the budget. One key difference between SLSQP and TC is that TC allows us to specify a tolerance for exceeding the budget, where SLSQP does not. We recorded the number of times each algorithm exceeded the budget by more than the tolerance, and by how much. We find that SLSQP exceeded the tolerance in 206 instances, with a maximum exceedance of $6.58 \times 10^{-6}$. FD was the only algorithm that never exceeded the tolerance.  



From these results, we can conclude four main points. Firstly, FD was the fastest algorithm of those that we tested. Secondly, FD was typically close to the best solution on average, and performed better than SLSQP and TC for Poisson demands. Thirdly, SLSQP and TC do not return integer solutions for the Poisson distribution in general. There is also no way to enforce this. Despite FD being designed for continuous order quantities, when the demand is discrete, its solution will always be integer. This is due to its solution being found from evaluating the inverse CDF of the demand distribution. This is an important point, as when demand is discrete the order quantity should also be. Finally, FD is the only algorithm that can be guaranteed not to exceed the budget constraint by more than some user-specified tolerance. 
\section{DRO model with normal demands}

In this section, we describe and solve the DRO model under normal demands. We first detail the model and ambiguity sets used in Section~\ref{sec:normal_model}. Following this, in Section~\ref{sec:normal_theory} we develop a set of extreme distributions based on theoretical results concerning the objective function. In Section~\ref{sec:normal_CS}, we use these results to develop our CS algorithm. Finally, in Section~\ref{sec:normal_results}, we perform extensive computational experiments to test the MLE approach and the CS algorithm. 

\subsection{Formulation and ambiguity sets}\label{sec:normal_model}

The first parametric DRO model that we study is the one resulting from independent normal demands. Suppose we have a set $\Omega$ of pairs $\bm{\omega} = (\bm{\mu}, \bm{\sigma})$
such that each $\bm{\omega}$ gives a unique multivariate normal distribution. This means there is a unique CDF $F$ for each $\bm{\omega}$ and so we can replace $F$ with $\bm{\omega}$.
Then the DRO model is:
\begin{align}
    \min_{\bm{q}}&\max_{\bm{\omega}\in \Omega}\label{eq:normal_first} C_{\bm{\omega}}(\bm{q})\\
    \text{s.t. }& \sum_{t=1}^T w_t q_t \le W\\
    & \bm{q} \in \R^T_+\label{eq:normal_last} 
\end{align}
Lemma~\ref{lemma:normal_obj} allows us to rewrite the objective function in a more convenient form.
\begin{lemma}\label{lemma:normal_obj}
    When $X_t \sim \m{N}(\mu_t, \sigma^2_t)$ for $t=1,\dots,T$ and the $X_t$'s are independent, we can rewrite $C_{\bm{\omega}}(\bm{q})$ as:
    \begin{equation}\label{eq:normal_obj}
        C_{\bm{\omega}}(\bm{q}) \hiderel{=}
         \sum_{t=1}^T\l\{(a_t \Phi(\beta_t) - h)\sum_{k=1}^t(\mu_k - q_k)
        + a_t \phi(\beta_t)\sqrt{\sum_{k=1}^t \sigma_k^2} +
        q_t w_t- p \mu_t\r\},
    \end{equation}
    with $\beta_t = \frac{ \sum_{k=1}^t (\mu_k - q_k)}{\sqrt{\sum_{k=1}^t
\sigma_k^2}}$.
\end{lemma}
The proof of Lemma~\ref{lemma:normal_obj} can be found in Appendix~C.1. Given this formulation, we next decide on the form that our ambiguity set $\Omega$ will take. As is common in the DRO literature, we will use maximum likelihood estimation to derive an approximate $100(1-\alpha)\%$ confidence set for the true parameter vector, which we will denote by $\bm{\omega}^0 = (\bm{\mu}^0, \bm{\sigma}^0)$.
Suppose for each $t = 1,\dots,T$ that we have access to a sample $\bm{x}_t = (x^1_t,\dots,x^N_t)$ of size $N$ from the true demand distribution $\m{N}(\mu^0_t, (\sigma^0_t)^2)$. We can use these samples to create MLEs of the true parameters, $\hat{\bm{\mu}} = (\hat{\mu}_1,\dots,\hat{\mu}_T)$ and $\hat{\bm{\sigma}} = (\hat{\sigma}_1,\dots,\hat{\sigma}_T)$. We also define $\hat{\bm{\omega}}$ as the MLE of the true parameter vector $\bm{\omega}^0$, i.e.\ $\hat{\bm{\omega}} = (\hat{\bm{\mu}}, \hat{\bm{\sigma}})$. It is a well-known result that the MLEs for normal parameters are given by:
\begin{equation*}
    \hat{\mu}_t = \frac{1}{N}\sum_{n=1}^N x^n_t, \quad
    \hat{\sigma}_t = \sqrt{\frac{1}{N} \sum_{n=1}^N (x^n_t - \hat{\mu}_t)^2}. \label{eq:MLEs}
\end{equation*}
Given these estimates, we can construct an approximate confidence set using Proposition~\ref{prop:normal_conf_set}. This result is proved in Appendix~C.2.
\begin{proposition}\label{prop:normal_conf_set}
    Suppose that $X_t \sim \m{N}(\mu^0_t, (\sigma^0_t)^2)$ for $t=1,\dots,T$, where $\bm{\omega}^0 = (\bm{\mu}^0, \bm{\sigma}^0)$ is unknown and the $X_t$'s are independent. Furthermore, suppose that $\hat{\mu}_t$ and $\hat{\sigma}_t$ are MLEs obtained from $N$ samples from the distribution of $X_t$, for $t=1,\dots,T$. Then, an approximate $100(1-\alpha)\%$ confidence set for $\bm{\omega}^0$ is given by:
    \begin{equation}\label{eq:normal_conf}
    \Omega_{\alpha} = \l\{(\bm{\mu}, \bm{\sigma}) \in \R^{2T}_{+}: \sum_{t=1}^T\l(\frac{N}{\hat{\sigma}^2_t}\l(\hat{\mu}_t - \mu_t\r)^2 + \frac{2N}{\hat{\sigma}^2_t}\l(\hat{\sigma}_t - \sigma_t\r)^2\r) \le \chi^2_{2T, 1-\alpha} \r\}.
    \end{equation}
\end{proposition}

In order for our model to be tractable, we assume a discrete ambiguity set. This means that the inner objective can be represented by a finite number of constraints. As described by~\cite{black}, parametric ambiguity sets do not yield tractable reformulations in any other way. Hence, in this paper, we will work with discrete subsets of the set in~\eqref{eq:normal_conf}. Since it is difficult to discretise this set directly, we will first construct a baseline ambiguity set $\Omega_{\text{base}}$ such that $\Omega_\alpha \subseteq \Omega_{\text{base}}$ and discretise this set. Then, we build a discretisation of $\Omega_\alpha$ by taking all elements of the discretisation of $\Omega_{\text{base}}$ that also lie in $\Omega_\alpha$. A logical way in which to construct $\Omega_{\text{base}}$ is as follows. The definition of~\eqref{eq:normal_conf} implies that any $(\bm{\mu}, \bm{\sigma}) \in \Omega_\alpha$ satisfies:
\begin{align*}
     \mu_t  &\in {\mu}^{\text{int}}_t = \l[\hat{\mu_t} - \hat{\sigma}_t\sqrt{\frac{\chi^2_{2T, 1-\alpha}}{N}}, \hat{\mu_t} + \hat{\sigma}_t\sqrt{\frac{\chi^2_{2T, 1-\alpha}}{N}}\r],  \fa t \in \{1,\dots,T\}\\
    \sigma_t &\in {\sigma}^{\text{int}}_t = \l[\hat{\sigma_t} - \hat{\sigma}_t\sqrt{\frac{\chi^2_{2T, 1-\alpha}}{2N}}, \hat{\sigma_t} + \hat{\sigma}_t\sqrt{\frac{\chi^2_{2T, 1-\alpha}}{2N}}\r] \fa t \in \{1,\dots,T\}.
\end{align*}
Therefore, defining $\Omega_{\text{base}}$ using~\eqref{eq:Omega}, we have $\Omega_\alpha \subseteq \Omega_{\text{base}}$.
\begin{equation}\label{eq:Omega}
    \Omega_{\text{base}} = \l\{(\bm{\mu}, \bm{\sigma}): \mu_t \in \mu^I_t \land \sigma_t \in \sigma^I_t \fa t = 1,\dots,T\r\}. 
\end{equation}
Let us refer to the upper and lower boundaries of the individual intervals as $\mu^l_t, \mu^u_t$ and $\sigma^l_t, \sigma^u_t$ for each $t = 1,\dots,T$. Then, we can create discretisations of these intervals using:
\begin{equation*}
    \tilde{\mu}^{\text{int}}_t = \l\{\mu^l_t + m\frac{\mu^u_t - \mu^l_t}{M - 1}:  m = 0,\dots, M - 1 \r\},
    \tilde{\sigma}^{\text{int}}_t = \l\{\sigma^l_t + m\frac{\sigma^u_t - \sigma^l_t}{M - 1}: m = 0,\dots,M-1
    \r\},
\end{equation*}
Then we can construct a discretisation of the baseline ambiguity set via:
\begin{equation*}
    \Omega^{\prime}_{\text{base}} = \l\{(\bm{\mu}, \bm{\sigma}): \mu_t \in \tilde{\mu}^{\text{int}}_t \land \sigma_t \in \tilde{\sigma}^{\text{int}}_t \fa t = 1,\dots,T\r\}.
\end{equation*}
Then, a discretisation of $\Omega_\alpha$ is given by $\Omega^{\prime}_{\alpha} = \Omega^{\prime}_{\text{base}} \cap \Omega_\alpha$. This is the set that we will use in our experiments. For the discrete ambiguity set $\Omega^{\prime}_{\alpha}$, we can reformulate the normal DRO model~\eqref{eq:normal_first}-\eqref{eq:normal_last} as:
\begin{equation}\label{eq:normal_reform}
    \min\l\{\theta: \theta \le C_{\bm{\omega}}(\bm{q}) \fa \bm{\omega} \in \Omega^{\prime}_{\alpha}, \sum_{t=1}^T w_t q_t \le W, \bm{q} \in \mathbb{R}^T_+\r\}.
\end{equation}
The first constraint here is implemented in Gurobi using piecewise linear approximations. For the piecewise linear approximation reformulation of this model, see Appendix~C.3.

\subsection{Extreme and dominated distributions}\label{sec:normal_theory}

Depending on the size of $\Omega^{\prime}_{\alpha}$, the piecewise linear model can become very slow. This is due to the fact that one piecewise linear approximation is required for every $\bm{\omega} \in \Omega^{\prime}_{\alpha}$ and every $t \in \{1,\dots,T\}$. The resulting model can have a significant number of constraints. In order to generate solutions in a shorter time, we will develop a CS algorithm. The idea of this algorithm is to iteratively solve the model on a small subset of $\Omega^{\prime}_{\alpha}$, and after each iteration generate a new distribution to add to the previous subset. Since finding the worst-case parameter for a given $\bm{q}$ can be cumbersome when $\Omega^{\prime}_{\alpha}$ is large, we will create a set $\Omega^{\text{ext}}$ of extreme parameters that are close to the worst case. In order to create this set, we derive the following properties of the objective function. These results are proved in Appendix~C.4.
\begin{theorem}\label{thm:normal_extreme}
    When $X_t \sim \m{N}(\mu_t, \sigma^2_t)$ for $t=1,\dots,T$ and the $X_t$'s are independent, $C_{\bm{\omega}}(\bm{q})$ is:
    \begin{enumerate}
        \item Convex in $\mu_t$ for all $t \in \{1,\dots,T\}$,
        \item Increasing in $\sigma_t$ for all $t \in \{1,\dots,T\}$.
    \end{enumerate}
\end{theorem}
Given these results, the question then
arises of how to combine them in order to find the most extreme $\bm{\omega}$
values. Let us define the sets of all distinct $\bm{\mu}$ and $\bm{\sigma}$
vectors as $\m{M}$ and $\m{S}$:
\begin{align*}
    \m{M} &= \l\{\bm{\mu} \in \R^T_{+}\ | \ \exists \ \bm{\sigma} \in \R^T_{+}: (\bm{\mu}, \bm{\sigma}) \in \Omega^{\prime}_{\alpha}\r\}\\
    \m{S} &= \l\{\bm{\sigma} \in \R^T_{+} \ | \ \exists \ \bm{\mu} \in \R^T_{+}: (\bm{\mu}, \bm{\sigma}) \in \Omega^{\prime}_{\alpha}\r\}.
\end{align*}
From $\m{M}$ and $\m{S}$ we can extract the most extreme $\bm{\mu}$'s and the most
extreme $\bm{\sigma}$'s individually using Theorem~\ref{thm:normal_extreme}. We first consider characterising the worst-case $\bm{\sigma}$. Since the cost function is increasing in each $\sigma_t$, it may seem reasonable that the worst-case $\bm{\omega}$ will have at least one $\sigma_t$ set at its maximum. However, if one $\sigma_t$ is at its maximum, then the others may need to be closer to $\hat{\sigma}_t$ in order for $\bm{\sigma}$ to remain sufficiently close to $\hat{\bm{\sigma}}$. In other words, maximising one $\sigma_t$ may mean reducing the others, and this can  yield a lower expected cost than if multiple standard deviations were allowed to be moderately large. In order to account for this fact, we consider $\bm{\sigma}$'s with the largest sum. This allows us to avoid prioritising any one day in particular and instead ensure that there is maximal total variability over all days. We now discuss the worst-case $\bm{\mu}$. Since the cost function is convex in each $\mu_t$, we will consider $\bm{\mu}$ such that at least one $\mu_t$ is equal to its maximum or minimum value. This is similar to the idea of maximising one $\sigma_t$, and therefore has the same problems. However, it is possible that some $\mu_t$'s should be maximised and others minimised. Therefore, maximising or minimising $\sum_{t=1}^T \mu_t$ is not logical here, since it implies that all should be maximised or all should be minimised. 

In addition, it may be the case that any $(\bm{\mu}, \bm{\sigma}) \in \Omega^{\prime}_{\alpha}$ that has $\sum_{t=1}^T \sigma_t$ maximised does not have any $\mu_t$ at its maximum or minimum. Similarly, any $(\bm{\mu}, \bm{\sigma})$ such that some $\mu_t$ is maximised or minimised may not have $\sum_{t=1}^T \sigma_t$ being maximised. Therefore, in order to construct a set of extreme parameters, we use a two-step procedure. For each
$\bm{\mu} \in \m{M}$ we can find the $\bm{\sigma}$'s with the largest sum over all $\bm{\omega} \in \Omega^{\prime}_{\alpha}$ whose mean is equal to $\bm{\mu}$. This corresponds to finding the most extreme distributions for this fixed $\bm{\mu}$. Hence, we define the set:
\begin{equation*} \label{eq:dominated_sigmas}
    \Omega^{\text{ext}}_1 = \l\{(\bm{\mu}, \bm{\sigma}) \in \Omega^{\prime}_{\alpha}: \sum_{t=1}^T \sigma_t = \max_{(\bm{\mu}, \tilde{\bm{\sigma}}) \in \Omega^{\prime}_{\alpha}} \sum_{t=1}^T \tilde{\sigma}_t\r\}.
\end{equation*}
Given this set, we now eliminate pairs from $\Omega^{\text{ext}}_1$ where
$\bm{\mu}$ is dominated when the corresponding $\bm{\sigma}$ is fixed. We can
characterise this by finding the maximum and minimum for each $\mu_t$
conditional on $\bm{\sigma}$. Then, we can remove any pairs from
$\Omega^{\text{ext}}_1$ that do not have any $\mu_t$ equal to their maximum or
minimum. Therefore, we define:
\begin{equation*}
    \mu^{\max}_t(\bm{\sigma}) = \max_{(\bm{\mu}, \bm{\sigma}) \in \Omega^{\text{ext}}_1} \mu_t, \ \mu^{\min}_t(\bm{\sigma}) = \min_{(\bm{\mu}, \bm{\sigma}) \in \Omega^{\text{ext}}_1} \mu_t,
\end{equation*}
\begin{equation*} \label{eq:dominated_mus}
    \Omega^{\text{ext}}_2 = \l\{(\bm{\mu}, \bm{\sigma}) \in \Omega^{\text{ext}}_1 \ | \
    \exists t \in \T: \mu_t \in \{\mu^{\min}_t(\bm{\sigma}), \mu^{\max}_t(\bm{\sigma})\}\r\}.
\end{equation*}
We can write this in a simpler form as a single set using:
\begin{equation}\label{eq:extreme_distributions}
    \Omega^{\text{ext}} = \l\{(\bm{\mu}, \bm{\sigma}) \in \Omega^{\prime}_{\alpha} \ \Bigg |\ \l(\sum_{t=1}^T \sigma_t = \max_{(\bm{\mu}, \tilde{\bm{\sigma}}) \in \Omega^{\prime}_{\alpha}} \sum_{t=1}^T \tilde{\sigma}_t\r) \land \l(\exists t \in \T: \mu_t \in \{\mu^{\min}_t(\bm{\sigma}), \mu^{\max}_t(\bm{\sigma})\}\r)\r\}
\end{equation}
Note that this set may not contain the true worst-case $\bm{\omega} \in \Omega^{\prime}_{\alpha}$, since it may be the case that selecting a less extreme mean/standard deviation in favour of a more extreme standard deviation/mean may lead to a worse distribution. However, in Section~\ref{sec:normal_results}, we will show that the worst-case in $\Omega^{\text{ext}}$ is typically very close to the true worst-case.

The size of $\Omega^{\prime}_{\alpha}$ can be too large for the model to be solved by PLA in a reasonable time. In addition to constructing $\Omega^{\text{ext}}$, Theorem~\ref{thm:normal_extreme} allows us to remove some dominated distributions that can never be the worst case. Consider two parameter values $\bm{\omega}^1$ and $\bm{\omega}^2$ such that $\bm{\mu}^1 = \bm{\mu}^2$, there is at least one $t$ such that $(\sigma^1)_t > (\sigma^2)_t$, and there is no $t$ such that $(\sigma^1)_t < (\sigma^2)_t$. Since, for a fixed $\mu$, the objective function is increasing in each $\sigma_t$, this means that $\bm{\omega}^2$ is \textit{dominated} by $\bm{\omega}^1$. Hence, $\bm{\omega}^2$ can never be a worse parameter than $\bm{\omega}^1$. Since this parameter can never be the worst-case, it can be removed from $\Omega^{\prime}_{\alpha}$ and the solution to the model will not change. In our experiments, in order to reduce the computational burden for PLA, we will remove all dominated distributions prior to solving. 


\subsection{Cutting surface algorithm}\label{sec:normal_CS}

We now present our CS algorithm. The idea of the algorithm is as follows. We start with some initial singleton ambiguity set $\Omega^1$. Then, for $k=1,\dots,k^{\max}$, we solve the piecewise linear approximation of the DRO model in~\eqref{eq:normal_reform} over $\Omega^{k}$ to get a solution $\bm{q}^k$. Then, we find the worst-case parameter $\bm{\omega}^k$ for the solution $\bm{q}^k$ over the set $\Omega^{\text{ext}}$ defined in~\eqref{eq:extreme_distributions}. We then set $\Omega^{k+1} = \Omega^{k} \cup \{\bm{\omega}^k\}$ and solve again. This process repeats until stopping criteria are met. The CS algorithm for this problem is as follows. 
\begin{enumerate}
    \item Initialise by choosing initial set of one distribution $\Omega^1 =
    \{\bm{\omega}^0\} \subseteq \Omega^{\prime}_{\alpha}$, optimality tolerance $\epsilon$, $k^{\max}$, and the
    gap $\varepsilon$ between successive $z$ points for the piecewise linear
    approximation.
    \item Construct set of extreme distributions $\Omega^{\text{ext}}$ using~\eqref{eq:extreme_distributions}.
    \item Set $k=1$.
    \item While $k \le k^{\max}$:
    \begin{enumerate}
        \item Solve piecewise linear approximation of the DRO model~\eqref{eq:normal_reform} using
        the set $\Omega^k$, to get solution $(\bm{q}^k, \tilde{\bm{\omega}})$.
        \item Use true objective function~\eqref{eq:normal_obj} to find objective value $\theta^k = C_{\tilde{\bm{\omega}}}(\bm{q}^k)$.
        \item Find worst-case cost $C^k = \max_{\bm{\omega} \in
        \Omega^{\text{ext}}} C_{\bm{\omega}}(\bm{q}^k)$ and use brute force computation to find worst-case
        parameter $\bm{\omega}^k \in \argmax_{\bm{\omega} \in
        \Omega^{\text{ext}}} C_{\bm{\omega}}(\bm{q}^k)$.
        \item If $C^k \le \theta^k + \frac{\epsilon}{2}$ or $\bm{\omega}^k \in \Omega^k$ then set $k=k^{\max} + 1$. 
        \item Else set $\Omega^{k+1} = \Omega^k \cup \{\bm{\omega}^k\}$ and $k =
        k+1$.
    \end{enumerate}
    \item Return solution $(\bm{q}^k, \bm{\omega}^k)$. 
\end{enumerate}
Note that, in step 4, CS uses the true objective function to evaluate $\bm{q}^k$ and to select a worst-case distribution for $\bm{q}^k$. Therefore, CS may give better solutions than solving the full model using a piecewise linear approximation. Furthermore, this algorithm typically ends in only a few iterations. This means that the DRO model is only ever solved over a small set of distributions, meaning that the piecewise linear models used solve very quickly. 

\subsection{Experiments on confidence-based ambiguity sets}\label{sec:normal_results}

In this section, we describe the design and results of our numerical experiments assessing the performance of our methods. In Section~\ref{sec:normal_expdes}, we describe the experimental design used. Following this, in Section~\ref{sec:MLE_performance} we study the order quantities and costs resulting from assuming that $\bm{\omega}^0 = \hat{\bm{\omega}}$. Finally, in Section~\ref{sec:normal_CSvsPLA} we assess CS as a heuristic for solving the full DRO model. These experiments were
run in parallel on a computing cluster (STORM) which has 486 CPU cores. The
solver used in all cases was the Gurobi Python package, gurobipy~\citep{gurobi}.
The version of gurobipy used was 9.0.1. The node used on STORM was the Tukey
node, which runs the Linux Ubuntu 16.04.7 operating system, Python version
2.7.12, and 46 AMD Opteron 6238 CPUs.

\subsubsection{Experimental design}\label{sec:normal_expdes}

In this section, we describe the parameter values used in our experiments. As an ambiguity set, we will use the discretisation $\Omega^{\prime}_{\alpha}$ of the confidence set~\eqref{eq:normal_conf} as described in Section~\ref{sec:normal_model}. The parameters used in our experiments are as follows. The number of periods was $T \in \{2,3,4\}$. The gap for the piecewise linear approximation in Appendix~C.3 was $\varepsilon \in \{0.1, 0.25, 0.50\}$ and the number of points used in the discretisation of univariate confidence intervals $M \in \{3,5,10\}$. The number of samples taken from the true distribution was $N\in \{10, 25, 50\}$. The model costs were $p, b, h \in \{100,200\}$ and $w_t = 100(T - t + 1)$ for $t=1,\dots,T$. The budget was $W = 4000$ when $T = 2,3$ and $W = 8000$ for $T = 4$. We randomly selected 3 values of $\bm{\mu}^0$ from the set $\{1,\dots,20\}^T$ and 3 values of $\bm{\sigma}^0$ from the set $\{1,\dots,10\}^T$. These values are chosen randomly for two reasons. Firstly, selecting only 3 keeps the number of experiments manageable. Secondly, selecting the values randomly means that we are not biasing the results by selecting them favourably.

Note that we only select means and standard deviations such that $3\sigma^0_t \le \mu^0_t \fa t$, to ensure that the demand distributions are unlikely to generate negative values. We also do not consider instances where $T = 4$ and $M=10$. This is because in these cases the base ambiguity sets are too large to be feasible. The reason why we use $W = 8000$ for $T=4$ is that the optimal solution if $W = 4000$ is always $(0,0,0,40)$. This is not interesting, and increasing the budget gives a wider variety of solutions. The above inputs yield 864 instances, not including those with $T=4$ and $M = 10$.

\subsubsection{Performance of MLE approaches}\label{sec:MLE_performance}

In this section, we assess the performance of the MLE approaches as a heuristic for solving the model under the true distribution. In particular, we study how close to optimality the solutions resulting from the MLE distributions are, and how well the MLE distribution's cost function approximates the true cost function. As a reminder, we use the term \textit{MLE cost function} to mean~\eqref{eq:normal_obj} with $\bm{\omega} = \hat{\bm{\omega}}$, i.e. $C_{\hat{\bm{\omega}}}$. Similarly, we use \textit{true cost function} to mean~\eqref{eq:normal_obj} with $\bm{\omega} = \bm{\omega}^0$, i.e. $C_{\bm{\omega}^0}$.
We ran the MLE approaches on the 864 instances generated by the inputs from Section~\ref{sec:normal_expdes}, solving each instance with both FD and PLA. 

\begin{figure}[htbp!]
    \centering
    \begin{subfigure}[t]{0.45\textwidth}
        \centering
        \includegraphics[width=0.95\textwidth]{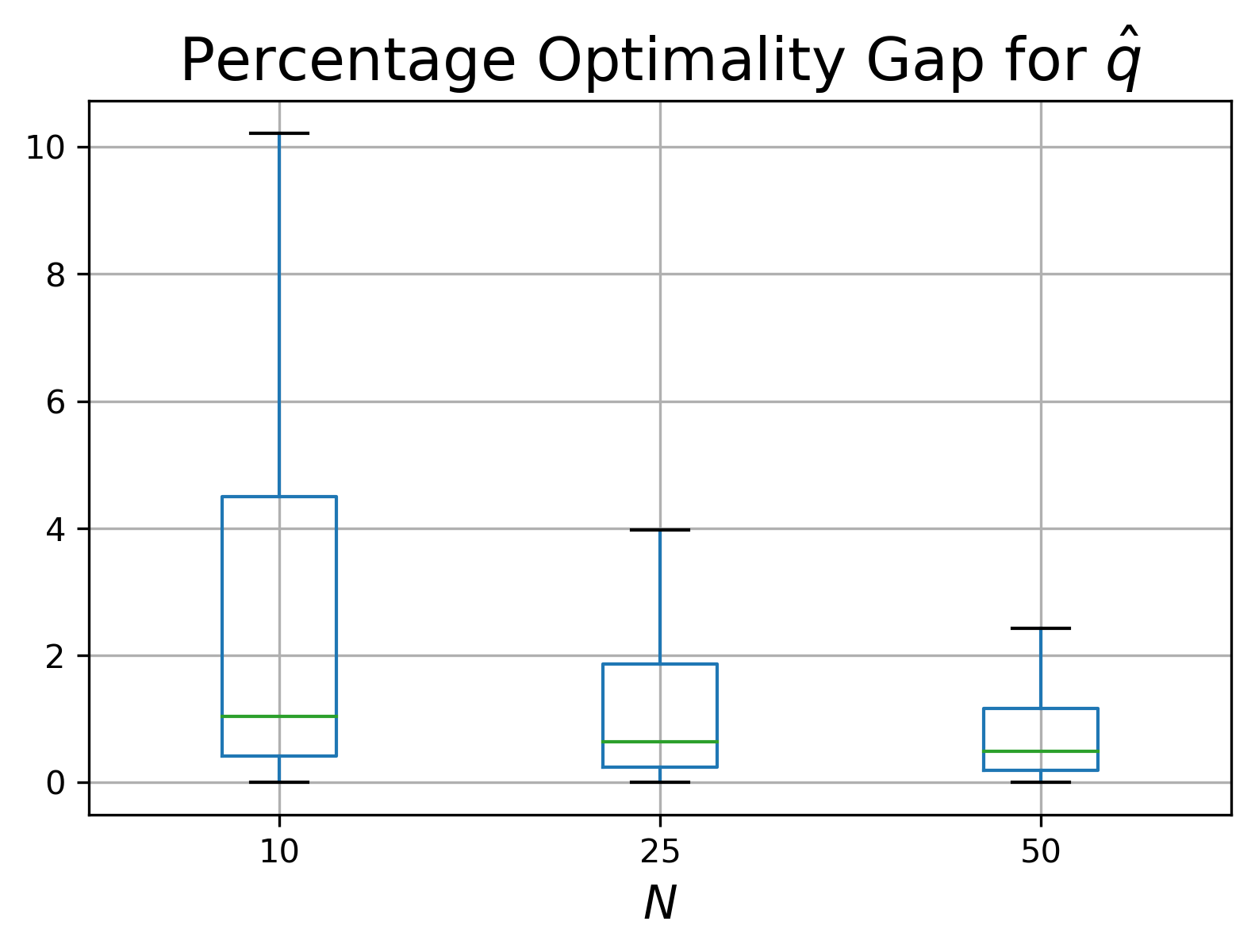}
        \caption{}
        \label{fig:FD_optgap}
    \end{subfigure}
    \begin{subfigure}[t]{0.45\textwidth}
        \centering
        \includegraphics[width=0.95\textwidth]{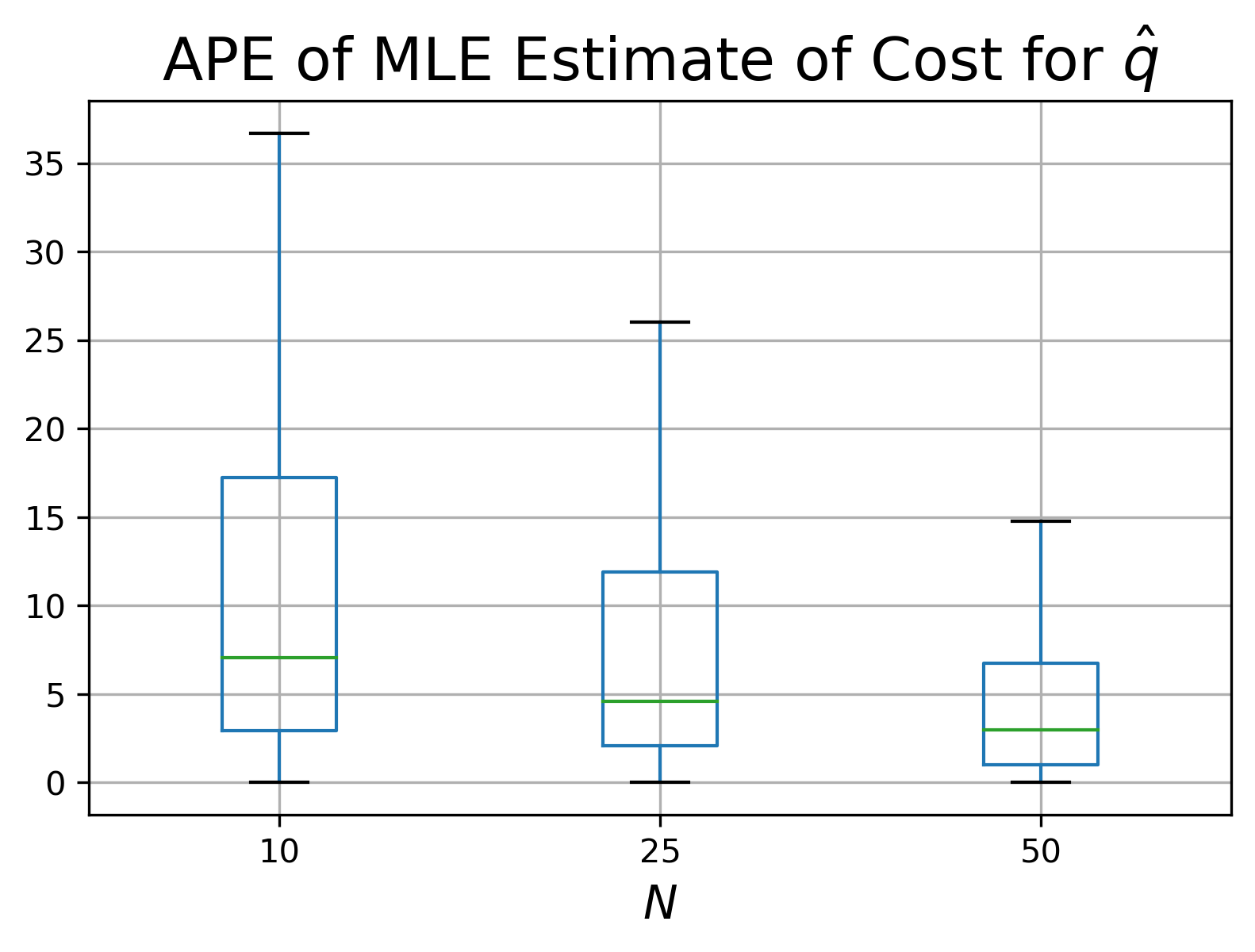}
        \caption{}
        \label{fig:FD_APE}
    \end{subfigure}
    \caption{Boxplots summarising performance of FD in comparison with true optimal solution. }
    \label{fig:FD_performance}
\end{figure}

We present boxplots summarising the performance of the MLE approach in Figure~\ref{fig:FD_performance}. Note that, as in Section~\ref{sec:FD_results}, boxplots are presented with outliers (points further than 1.5IQR above the 75th percentile or below the 25th percentile) removed. Let $\hat{\bm{q}}$ be the solution obtained from applying FD assuming $\bm{\omega} = \bm{\hat{\omega}}$. Figure~\ref{fig:FD_optgap} shows the absolute percentage optimality gap of $\hat{\bm{q}}$ when compared with the optimal $\bm{q}$ under $\bm{\omega}^0$, given by:
$$100 \times \l\lvert\frac{\min_{\bm{q}} C_{\bm{\omega}^0}(\bm{q})
- C_{\bm{\omega}^0}(\hat{\bm{q}})}{\min_{\bm{q}} C_{\bm{\omega}^0}(\bm{q})}\r\rvert.$$
This plot shows that $\hat{\bm{q}}$ is quite close to optimality under $\bm{\omega}^0$. For the smallest samples, i.e.\ $N=10$, it typically was no further than 10\% from optimal. The optimality gap reduces as $N$ increases. We also study the absolute percentage error (APE) of $C_{\hat{\bm{\omega}}}(\hat{\bm{q}})$ as an estimate of $C_{\bm{\omega}^0}(\hat{\bm{q}})$, calculated as:
$$100 \times \l\lvert\frac{C_{\bm{\omega}^0}(\hat{\bm{q}})
- C_{\hat{\bm{\omega}}}(\hat{\bm{q}})}{C_{\bm{\omega}^0}(\hat{\bm{q}})}\r\rvert.$$
Figure~\ref{fig:FD_APE} shows that the MLE distribution can result in poor estimates of the cost of $\hat{\bm{q}}$. The plot shows that for $N=10$, the APE can reach over 35\%. The accuracy of the estimate increases with $N$, but even for $N=50$ we see errors of 15\%.

In addition, these errors can be even more extreme than shown here, due to outliers being removed from the plot. In one instance, the APE was as high as 15,600\%. This occurred in a 2-period instance with $N=25$, and the corresponding costs were $C_{\hat{\bm{\omega}}}(\hat{\bm{q}}) = -386.06$ and $C_{\bm{\omega}^0}(\hat{\bm{q}}) = -2.46$. In addition to this, we found that in 21 instances the MLE solution gave an order quantity that was predicted to make a profit but ultimately resulted in a loss. In the most extreme example of this, we found $C_{\hat{\bm{\omega}}}(\hat{\bm{q}}) = -479.12$ and $C_{\bm{\omega}^0}(\hat{\bm{q}}) = 52.42$. The lowest cost under $\bm{\omega}^0$ in this instance was $-5.42$, meaning that it was possible to make a profit, but that using the MLE approach would have resulted in a loss. Interestingly, this happened for all values of $N$ and was not restricted to $N=10$.  In general, we found that the MLE approach was likely to underestimate the cost of its chosen $\bm{q}$. This was the case in 61\% of instances.

We can also confirm that these results are not simply due to the performance of FD as a heuristic for MPNVP. We found that FD and PLA's solutions were typically within -2\% and 1.5\% of one another in terms of cost. We also found that FD sometimes gave better solutions than PLA. In addition, we found similar results about the performance of the MLE approach when the model was solved using PLA. Namely, the MLE cost function predicted a profit for PLA's solution when there would in fact be a cost in 16 instances, and it underestimated the cost of PLA's decision in 59\% of instances.


From the results shown here, we can make three key conclusions. Firstly, using the MLE distribution will typically lead to solutions that are reasonably close to optimal under the true distribution. Secondly, we can conclude that the MLE distribution can lead to poor estimates of associated costs associated with a given $\bm{q}$. In particular, solving using the MLEs can lead to the suggestion that the resulting decision will lead to a profit when it will actually lead to a cost. Finally, we can conclude that using the MLE distribution is likely to lead to an underestimation of the cost associated with the selected decision. We can confirm that the DRO objective in~\eqref{eq:normal_first} never predicted a profit when the solution would actually result in a loss. For the aforementioned instance where $\hat{\bm{q}}$ had a predicted cost of $-479.12$ and an actual cost of 52.4, the cost that would be associated with this decision under the DRO model was 156.84.


\subsubsection{Performance of DRO algorithms}\label{sec:normal_CSvsPLA}

Each algorithm was given a maximum time of four hours in which to carry out all pre-computation and to build and solve the resulting models. By pre-computation, we mean computing the sets $\Omega^{\prime}_{\text{base}}$ and $\Omega^{\prime}_{\alpha}$. Since this precomputation was required to generate the inputs for both algorithms and is not technically part of either algorithm, we do not include this in the times that we now present. Instead, we present \textit{run times}, i.e.\ the amount of time the algorithms ran for after all precomputation was complete, until it either finished running or was timed out. We use the term {timed out} to mean that the algorithm's run time exceeded the time that was remaining after all precomputation was complete. For example, if an algorithm timed out with a run time of 1 minute, this meant that it spent 3 hours and 59 minutes in pre-computation.  

Of the {864} instances that we ran, 6 instances resulted in ambiguity sets with only 1 parameter. Hence, they were skipped. This means that we have results for 858 instances. Of these 858 instances, PLA timed out in {1}. CS did not time out in any instance. PLA took 7 minutes and 26 seconds on average, whereas CS took just 1.2 seconds. The 25\% and 75\% quantiles for PLA's time taken are 0.7 seconds and 52.6 seconds. For CS, these values are 0.72 seconds and 3.3 seconds. In addition, PLA ran for up to a maximum of 3 hours and 53 minutes, whereas CS always finished running in less than 8 seconds. Note that PLA never ran for more than 4 hours because the time limit accounts for pre-computation time and run time. In the above example, the pre-computation took 7 minutes and PLA ran for 3 hours and 53 minutes. It is worth noting that PLA's average time is strongly affected by some large run times. PLA took no more than 2 minutes to finish running in 79\% of instances. However, it took over 3 hours in approximately 1\% of instances. 

\begin{figure}[htb!]
    \begin{subfigure}[t]{0.45\textwidth}
        \centering
        \includegraphics[width=0.95\textwidth]{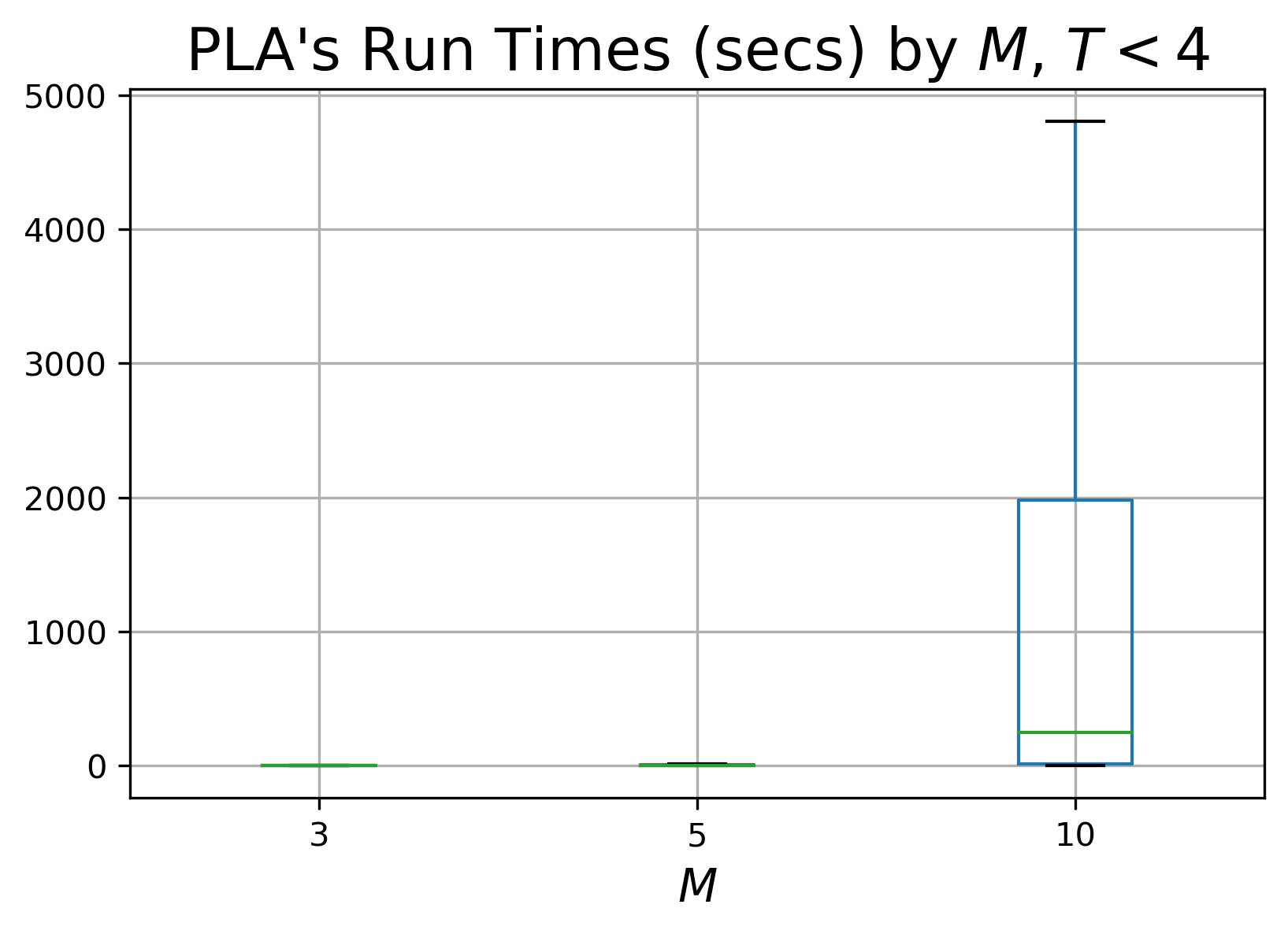}
        \caption{}
        \label{fig:normal_PLA_tt_by_M2}
    \end{subfigure}
    \begin{subfigure}[t]{0.45\textwidth}
        \centering
        \includegraphics[width=0.95\textwidth]{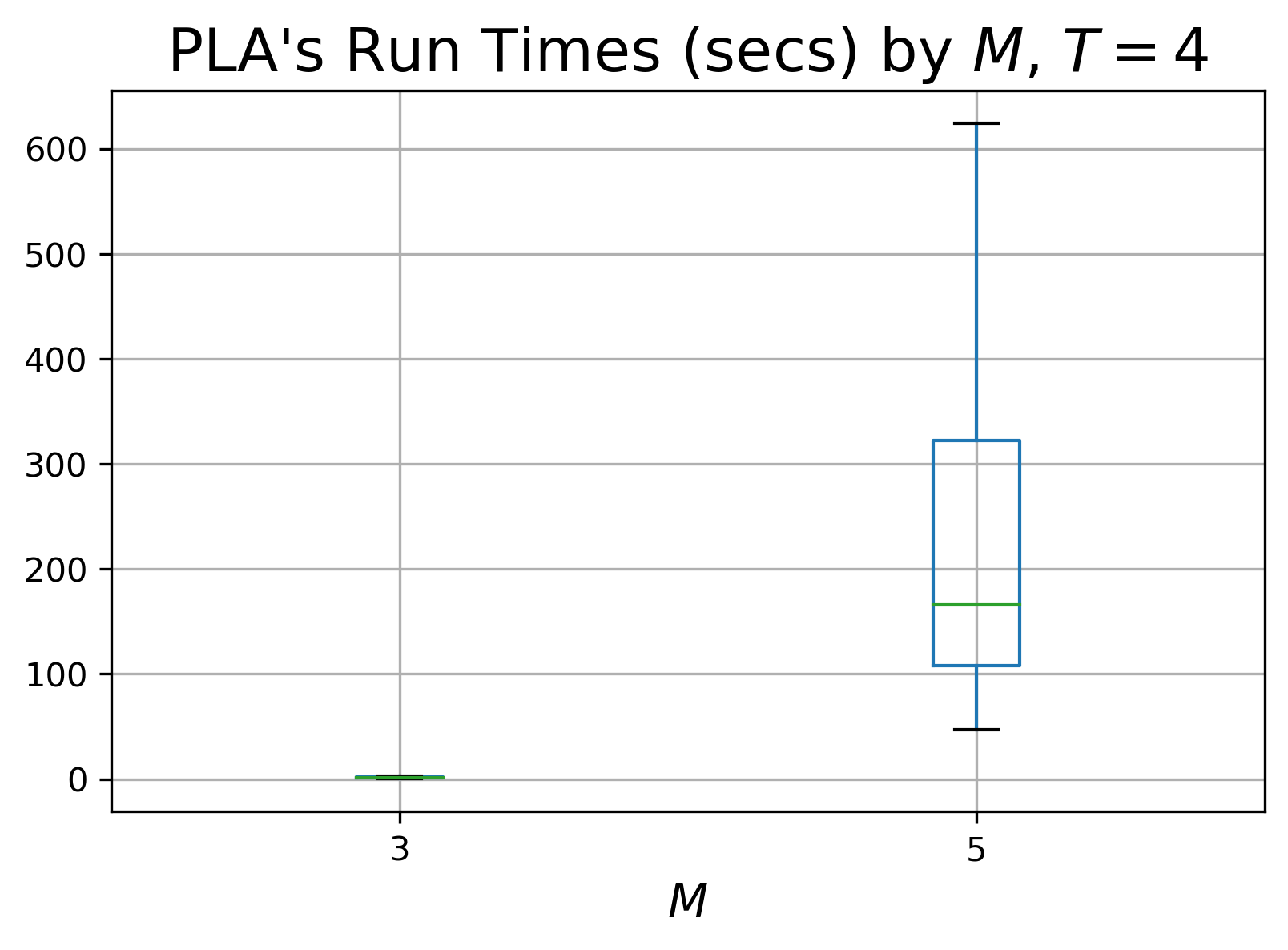}
        \caption{}
        \label{fig:normal_PLA_tt_by_M}
    \end{subfigure}
    \caption{Boxplots of run times of PLA, grouped by $M$.}
    \label{fig:normal_timestaken_by_M_PLA}
\end{figure}

\begin{figure}[htb!]
    \begin{subfigure}[t]{0.45\textwidth}
        \centering
        \includegraphics[width=0.95\textwidth]{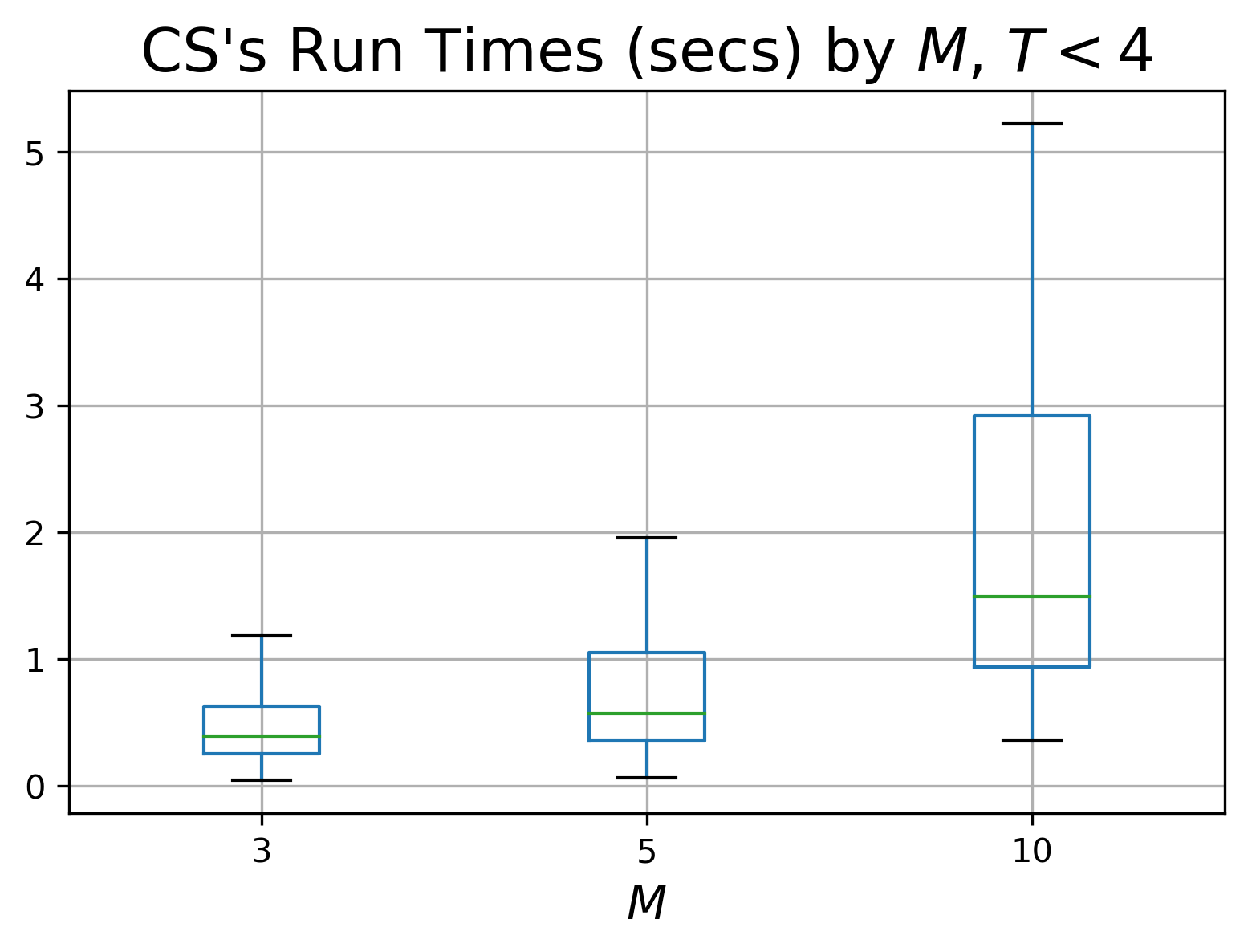}
        \caption{}
        \label{fig:normal_CS_tt_by_M_2}
    \end{subfigure}
    \begin{subfigure}[t]{0.45\textwidth}
        \centering
        \includegraphics[width=0.95\textwidth]{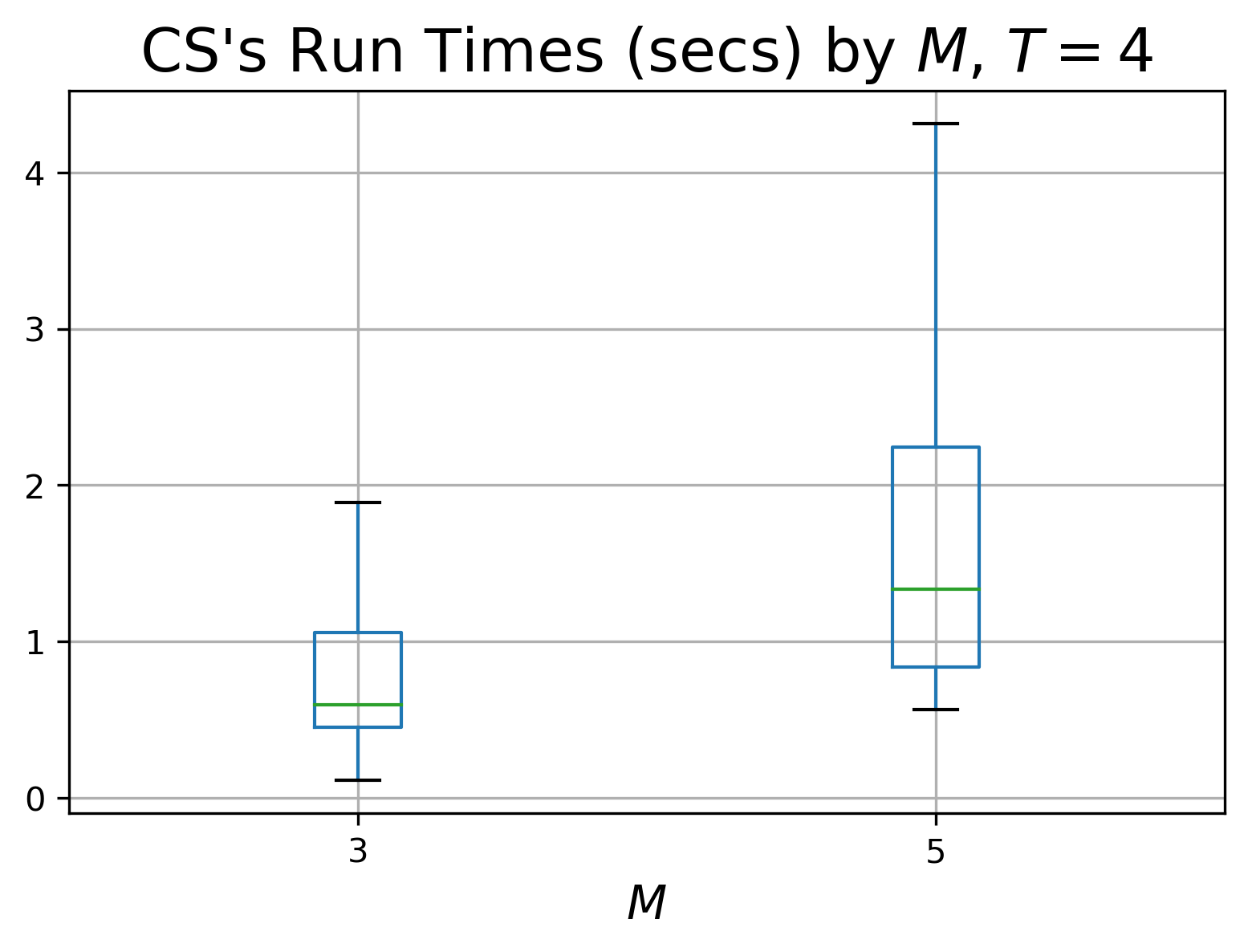}
        \caption{}
        \label{fig:normal_CS_tt_by_M}
    \end{subfigure}
    \caption{Boxplots of run times of CS, grouped by $M$.}
    \label{fig:normal_timestaken_by_M_CS}
\end{figure}

We can now look at the factors affecting the run times. From our analysis, the main factor affecting PLA's run times is {$M$}. For large $M$, the base ambiguity set $\Omega^{\prime}_{\text{base}}$ becomes very large. Therefore, $\Omega^{\prime}_{\alpha}$ is also large and even constructing it takes quite some time. To highlight the effect of $M$ on PLA and CS's run times, we present boxplots grouped by $M$ in Figures~\ref{fig:normal_timestaken_by_M_PLA} and~\ref{fig:normal_timestaken_by_M_CS}. Since $M=10$ was not run for $T=4$, we separate these plots into two. Figure~\ref{fig:normal_PLA_tt_by_M} shows that PLA begins to run very slowly for $M=10$ when $T < 4$. In fact, PLA's timeout occurred when $T=3$ and $M=10$. Figure~\ref{fig:normal_PLA_tt_by_M2} shows that PLA is also slow for $T=4$ and $M=5$, but it did not time out in any of these instances. The effect of $M$ on CS is not so drastic. As shown in Figure~\ref{fig:normal_CS_tt_by_M_2}, for $T<4$, CS's run times do increase with $M$, but not by a large amount. For $T=4$, Figure~\ref{fig:normal_CS_tt_by_M} shows that CS does take longer to run for $M = 5$ than for $M=3$, but this difference is much less noticeable. We do not present run times by $T$, since these plots would simply mirror the effect of $M$.

We now assess the performance of CS as a heuristic for solving the full DRO model. This is done via comparing its objective values and returned $\bm{\omega}$ values with those from PLA. We first present the number of times that CS and PLA obtained the worst-case parameters for their chosen $\bm{q}$. We present these values for PLA as well as CS since both algorithms use a reduced ambiguity set. CS uses $\Omega^{\text{ext}}$ and PLA uses $\Omega^{\prime}_{\alpha}$ with dominated parameters removed. Hence, studying these results allows us to confirm that the reduced set used by PLA contains the true worst-case parameter for each solution. We present the optimality counts in Table~\ref{tab:normal_worstcase}. Note that these counts are given only over the instances where the algorithm concerned finished running. Since PLA timed out in one instance, it finished running in 1 less instance than CS.

CS chose the worst-case parameters for its selected $\bm{q}$ in 85\% of instances.  PLA chose the worst-case parameters in every instance. This is expected, since despite the fact that the approximate objective can lead to the wrong worst-case parameter, we use the true objective function to generate the cost of PLA's solution after it finishes running. The main purpose of this is to confirm that our removal of dominated distributions did not remove the worst-case parameter in any instance. CS failed to select the worst-case parameters in 126 instances.

 \begin{table}[htbp!]
    \centering
    \begin{tabular}{lrlll}
\toprule
 $T$ &  CS No. & CS Worst-case Count (\%) & PLA No. & PLA Worst-case Count(\%) \\
\midrule
 All &     858 &            732 (85.31\%) &     857 &            857 (100.0\%) \\
   2 &     320 &            270 (84.38\%) &     320 &            320 (100.0\%) \\
   3 &     323 &            257 (79.57\%) &     322 &            322 (100.0\%) \\
   4 &     215 &            205 (95.35\%) &     215 &            215 (100.0\%) \\
\bottomrule
\end{tabular}

    \caption{Number of instances in which algorithms selected true worst case parameters.}   
    \label{tab:normal_worstcase}
\end{table}

Table~\ref{tab:normal_worstcase} only shows the number of times CS chose the worst-case parameters. In order to gain a better understanding of its performance, we can study two quantities:
\begin{enumerate}
    \item Percentage $\bm{\omega}$-gap: The percentage difference between the worst-case cost for $\bm{q}^{\text{CS}}$ and the cost obtained under the $\bm{\omega}$ returned by CS. Defined as:
    $$100\times \frac{\max_{\bm{\omega} \in \Omega^{\prime}_{\alpha}} C_{\bm{\omega}}(\bm{q}^{\text{CS}}) - C_{\bm{\omega}^{\text{CS}}}(\bm{q}^{\text{CS}})}{\max_{\bm{\omega} \in \Omega^{\prime}_{\alpha}}C_{\bm{\omega}}(\bm{q}^{\text{CS}})}.$$
    \item Percentage $\bm{q}$-gap: The percentage difference between the worst-case cost attained by $\bm{q}^{\text{CS}}$ and that attained by $\bm{q}^{\text{PLA}}$. Defined as:
    $$100 \times \frac{\max_{\bm{\omega} \in \Omega^{\prime}_{\alpha}} C_{\bm{\omega}}(\bm{q}^{\text{PLA}}) - \max_{\bm{\omega} \in \Omega^{\prime}_{\alpha}} C_{\bm{\omega}}(\bm{q}^{\text{CS}})}{\max_{\bm{\omega} \in \Omega^{\prime}_{\alpha}} C_{\bm{\omega}}(\bm{q}^{\text{PLA}})}.$$
\end{enumerate}

\begin{figure}[htbp!]
    \begin{subfigure}[t]{0.45\textwidth}
        \centering
        \includegraphics[width=0.95\textwidth]{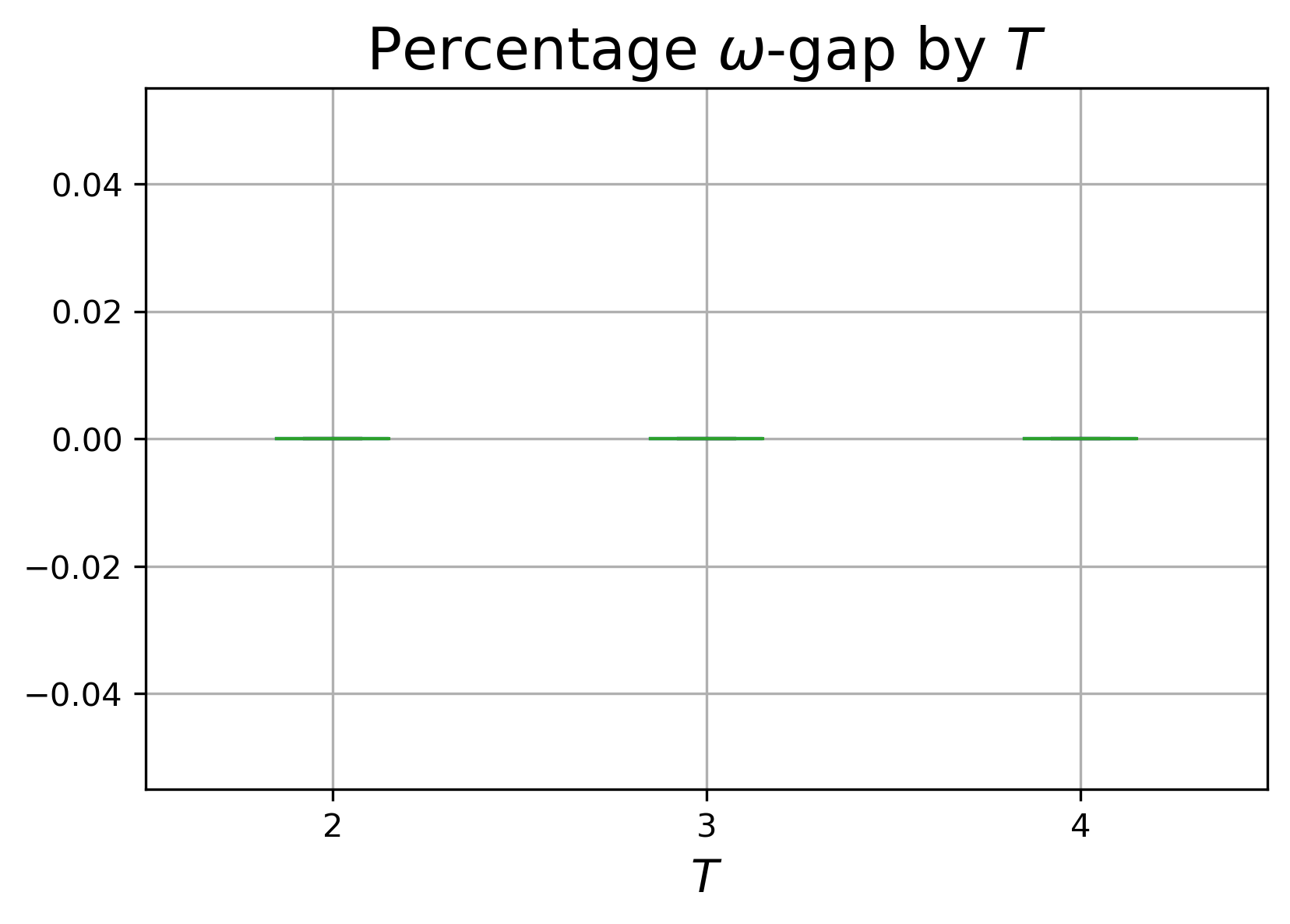}
        \caption{}
        \label{fig:normal_omega_gap}
    \end{subfigure}
    \begin{subfigure}[t]{0.45\textwidth}
        \centering
        \includegraphics[width=0.95\textwidth]{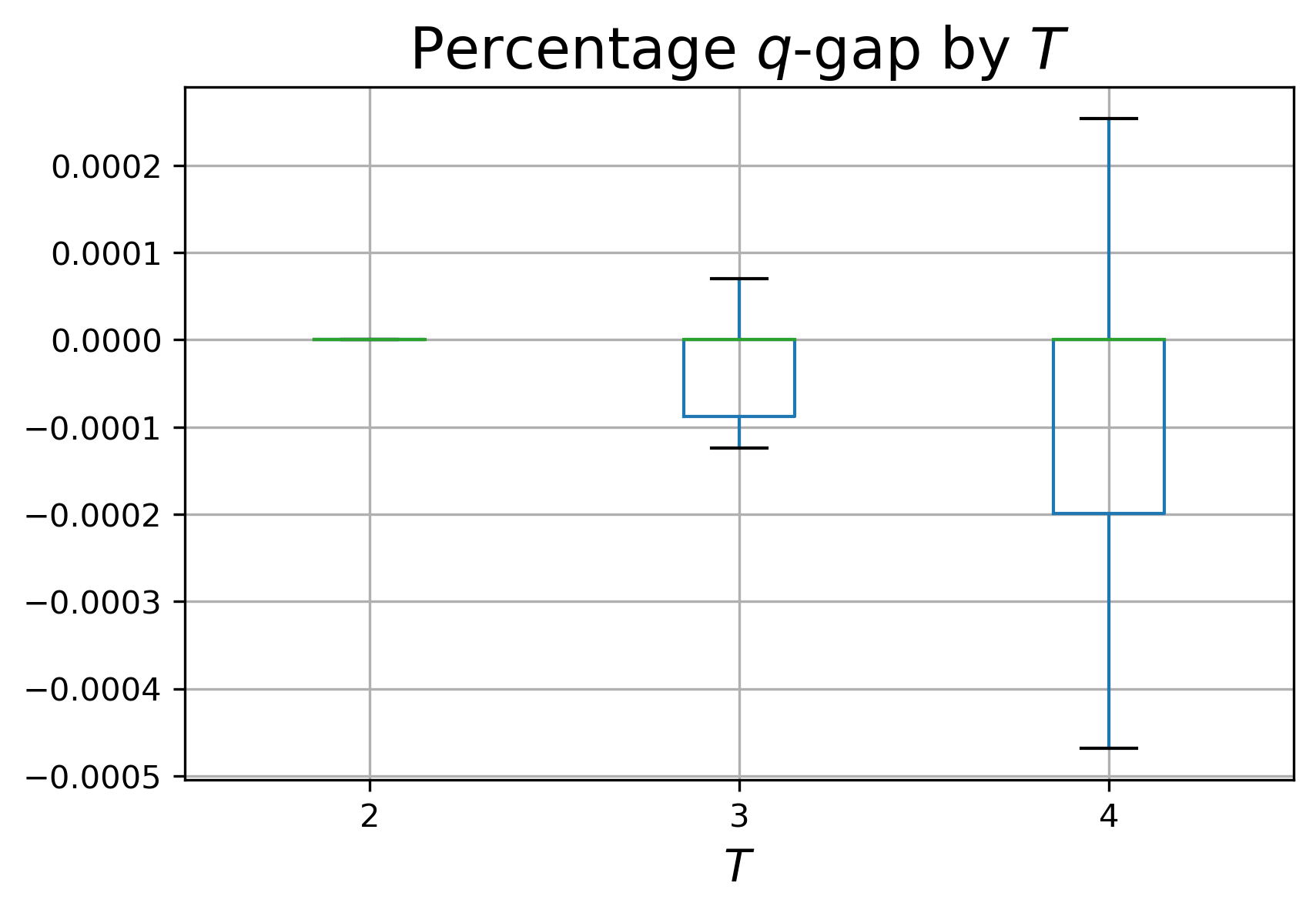}
        \caption{}
        \label{fig:normal_q_gap}
    \end{subfigure}
    \caption{Summary of percentage gaps}
    \label{fig:normal_gaps}
\end{figure}

Note that positive values (lower expected cost) for the $\bm{\omega}$-gap implies that CS gave a suboptimal $\bm{\omega}$. Negative values (higher worst-case expected cost) for the $\bm{q}$-gap imply that CS gave a suboptimal $\bm{q}$. Figure~\ref{fig:normal_gaps} summarises the percentage gaps using boxplots. The first thing to notice is that the $\bm{\omega}$-gaps were so commonly zero that everything non-zero was considered as an outlier. On average, the percentage $\bm{\omega}$-gap was 0.045\%. We see from Figure~\ref{fig:normal_q_gap} that CS's chosen $\bm{q}$ had a worst-case cost that was very close to PLA's, in general. The percentage optimality gaps typically do not exceed 0.003\%. The overall average percentage $\bm{q}$-gap was $-0.06\%$. CS even performed better than PLA in 153 instances. Since CS can be viewed as a heuristic for solving the full model that PLA solves, this may be unexpected. This occurs due to one main reason. CS's solution can have a higher worst-case cost than PLA's solution under the approximate cost function used by PLA, even when the worst-case cost of CS's solution is lower than that of PLA's under the true cost function.

From these results we can conclude 4 main points. Firstly, CS was much faster than PLA and scales better. It never took more than 8 seconds to finish running, whereas PLA took up to 3 hours and 53 minutes. Secondly, PLA performed better in selecting the worst-case parameters for a fixed $\bm{q}$. This is expected, as CS applies cardinality reduction to the ambiguity set. CS's returned parameters were very close to the worst-case in terms of cost, however. Thirdly, CS was very close to optimal with respect to its chosen $\bm{q}$. On average, it returned solutions that are only $0.045\%$ away from optimal. Finally, CS performed better than PLA in some instances, due to the inaccuracy of the piecewise linear approximation.

\section{DRO model with Poisson demands}

In this section, we present and solve the DRO model when the demands are Poisson random variables. In Section~\ref{sec:poisson_model}, we give the model and describe the ambiguity sets used. Then, in Section~\ref{sec:poisson_extreme}, we describe the set of extreme parameters used by CS. Finally, in Section~\ref{sec:poisson_results}, we test the MLE approach and our CS algorithm using computational experiments.

\subsection{Formulation and ambiguity sets}\label{sec:poisson_model}

When the demands are independent Poisson random variables, i.e.\ $X_t \sim \text{Pois}(\lambda^0_t)$ for unknown $\bm{\lambda}^0$, then we have that $\sum_{k=1}^t X_k \sim \text{Pois}\l(\sum_{k=1}^t \lambda^0_k\r)$ for $t=1,\dots,T$. We will now consider integer decision variables, i.e.\ $\bm{q} \in \mathbb{N}^T_0$, since the demand is discrete. Furthermore, we can use Lemma~\ref{lemma:pois_obj} to rewrite the objective function in a more convenient fashion. This lemma is proved in Appendix~D.1.
\begin{lemma}\label{lemma:pois_obj}
    Suppose that the demands are independent Poisson random variables, i.e.\ $X_t \sim \text{Pois}(\lambda_t) \fa t = 1,\dots,T$. Then, for $\bm{q} \in \mathbb{N}^T_0$, the objective function $C_{\bm{\lambda}}(\bm{q})$ can be written as:
    \begin{equation*}
        C_{\bm{\lambda}}(\bm{q}) = \sum_{t=1}^T\l(a_t Q_t \tilde{F}_t(Q_t) - \Lambda_t a_t \tilde{F}_t(Q_t-1) + (b+p\mathds{1}\{t=T\})(\Lambda_t - Q_t) + w_t q_t - p \lambda_t\r),
    \end{equation*}
    with $a_t = h + b + p\mathds{1}\{t = T\}$ and $\Lambda_t = \sum_{k=1}^t \lambda_k$.
\end{lemma}
Therefore, the DRO model for Poisson demands is given by:
\begin{align}
    \min_{\bm{q}}&\max_{\bm{\lambda} \in \Omega} C_{\bm{\lambda}}(\bm{q}) = \sum_{t=1}^T\l(a_t Q_t \tilde{F}_t(Q_t) - \Lambda_t a_t \tilde{F}_t(Q_t-1) + (a_t - h)(\Lambda_t - Q_t) + w_t q_t - p \lambda_t\r)\label{eq:poisson_obj}\\
    \text{s.t. } & \sum_{t=1}^T w_t q_t \le W,\\
    & q_t \in \mathbb{N}_0  \fa t \in \{1,\dots,T\}.\label{eq:poisson_last}
\end{align}
Similarly to the case for normally distributed demands, we will construct confidence sets for use as ambiguity sets. Suppose that we draw $N$ samples $\bm{x}^n = (x^n_1,\dots,x^n_T)$ from $\bm{X}$. Then, the MLEs are given by $\hat{\lambda}_t = \frac{1}{N} \sum_{n=1}^N x^n_t$  for $t=1,\dots,T$.
We can construct an approximate $100(1-\alpha)\%$ confidence set for $\bm{\lambda}^0$ using Proposition~\ref{prop:pois_conf_set}, which is proved in Appendix~D.2.
\begin{proposition}\label{prop:pois_conf_set}
    Suppose that $X_t \sim \text{Pois}(\lambda^0_t)$ for $t \in \{1,\dots,T\}$ and that the $X_t$ are independent. Then, an approximate $100(1-\alpha)\%$ confidence set for $\bm{\lambda}^0$ is given by:
    \begin{equation}\label{eq:pois_conf_set}
    \Omega_{\alpha} = \l\{\bm{\lambda} \in \R^T_{+}: \sum_{t=1}^T \frac{N}{\hat{\lambda}_t}(\hat{\lambda}_t - \lambda_t)^2 \le \chi^2_{T, 1-\alpha}\r\}.
\end{equation}
\end{proposition}
As before, we will construct a discretisation of $\Omega_\alpha$. In the Poisson case, we have that $\bm{\lambda} \in \Omega_\alpha$ implies that:
\begin{align*}
    \lambda_t \in \lambda^{\text{int}}_t &= \l[\hat{\lambda}_t -  \sqrt{\frac{\chi^2_{T, 1-\alpha}\hat{\lambda}_t}{N}} \le \lambda_t \le \hat{\lambda}_t +  \sqrt{\frac{\chi^2_{T, 1-\alpha}\hat{\lambda}_t}{N}}\r] \quad (t=1,\dots,T).
\end{align*}
Therefore, if we define $\Omega_{\text{base}} = \l\{\bm{\lambda} \in \R^T_{+}: \lambda_t \in {\lambda}^{\text{int}}_t \fa t = 1,\dots,T \r\}$, then we have $\Omega_{\alpha} \subseteq \Omega_{\text{base}}$.
As before, we can discretise the intervals $\lambda^{\text{int}}_t$ using:
\begin{equation*}
    \tilde{\lambda}^{\text{int}}_t = \l\{\lambda^l_t + m \frac{\lambda^u_t - \lambda^l_t}{M-1}, m=0,\dots,M-1\r\},
\end{equation*}
where $\lambda^l_t, \lambda^u_t$ are the lower and upper bounds of $\lambda^{\text{int}}_t$.  Hence, we can calculate a discretisation of $\Omega_{\text{base}}$ by taking $\Omega^{\prime}_{\text{base}} = \{\bm{\lambda} \in \R^{T}_{+}: \lambda_t \in \tilde{\lambda}^{\text{int}}_t \fa t \in \{1,\dots, T\}\}$.
Finally, we can take a discretisation of the confidence set $\Omega_{\alpha}$ by taking  $\Omega^{\prime}_{\alpha} = \Omega^{\prime}_{\text{base}} \cap \Omega_\alpha$. This is the set used in our experiments. For the discrete ambiguity set $\Omega^{\prime}_{\alpha}$, we can reformulate~\eqref{eq:poisson_obj}-\eqref{eq:poisson_last} as:
\begin{equation*}
    \min\l\{\theta: \theta \le C_{\bm{\lambda}}(\bm{q}) \fa \bm{\lambda} \in \Omega^{\prime}_{\alpha}, \sum_{t=1}^T w_t q_t \le W, \bm{q} \in \mathbb{N}^T_0\r\}.
\end{equation*}
This model is implemented in Gurobi using piecewise linear constraints. For  more details on how the model is implemented, see Appendix~D.3.

\subsection{Extreme distributions}\label{sec:poisson_extreme}

As for the case of normally distributed demands, we will develop a CS algorithm that produces solutions in a comparatively short time. In order to construct the set $\Omega^{\text{ext}}$ of extreme parameters, we will use the result in Theorem~\ref{thm:pois_convex}.
\begin{theorem}\label{thm:pois_convex}
    Suppose that $X_t \sim \text{Pois}(\lambda^0_t)$ for $t \in \{1,\dots,T\}$ and that the $X_t$ are independent. Then $C_{\bm{\lambda}}(\bm{q})$ is convex in $\lambda_t$ for each $t=1,\dots,T$.
\end{theorem}
This result is proved in Appendix~D.4. This result means that the worst-case $\lambda_t$ will either be the smallest in the ambiguity set or the largest (when all other $\lambda_j$ are fixed). Therefore, we can construct a set of extreme distributions for a given $\Omega^{\prime}_\alpha$ as follows:
\begin{equation}\label{eq:poisson_extreme}
    \Omega^{\text{ext}} = \l\{\bm{\lambda} \in \Omega^{\prime}_{\alpha} \bigg| \exists t \in \{1,\dots,T\}: \lambda_t \in \l\{ \min_{\bm{\nu} \in \Omega^{\prime}_{\alpha}} \nu_t, \max_{\bm{\nu} \in \Omega^{\prime}_{\alpha}} \nu_t\r\}\r\}.
\end{equation}
The cutting surface algorithm is therefore the same as the one in Section~\ref{sec:normal_CS}, but using~\eqref{eq:poisson_extreme} as the set of extreme distributions and $C_{\bm{\lambda}}$ instead of $C_{\bm{\omega}}$.

\subsection{Experimental design and results}\label{sec:poisson_results}

To test the CS algorithm on Poisson demands, we run experiments on the same inputs as for the normal distribution. In every case, we take the value of $\bm{\lambda}^0$ to be equal to the sampled value of $\bm{\mu}_0$ in the normal experiments. This allows to see the effect on our results if the only aspect allowed to change is the distribution used. Note, however, that we now run the tests with $T = 4$ and $M = 10$, since there is now only one parameter per day and so the Poisson ambiguity sets are much smaller. Since the objective function is piecewise constant between integer points, we can simply use an $\varepsilon$ of 1 in all instances. Hence, we run our experiments for an additional 6 randomly selected values of $\bm{\lambda}^0$, in order to generate more instances. These parameters yield 972 instances. Note that, as in the normal case, all boxplots have outliers (points further than 1.5IQR above the 75th percentile or below the 25th percentile) removed. These experiments were also run on STORM, this time on the Wald node. The Wald node runs 70 Intel Xeon E5-2699 CPUs, and also runs Python 2.7.12 and Gurobipy 9.0.1.

\subsubsection{Performance of MLE approaches}

\begin{figure}[htbp!]
    \centering
    \begin{subfigure}[t]{0.45\textwidth}
        \centering
        \includegraphics[width=0.95\textwidth]{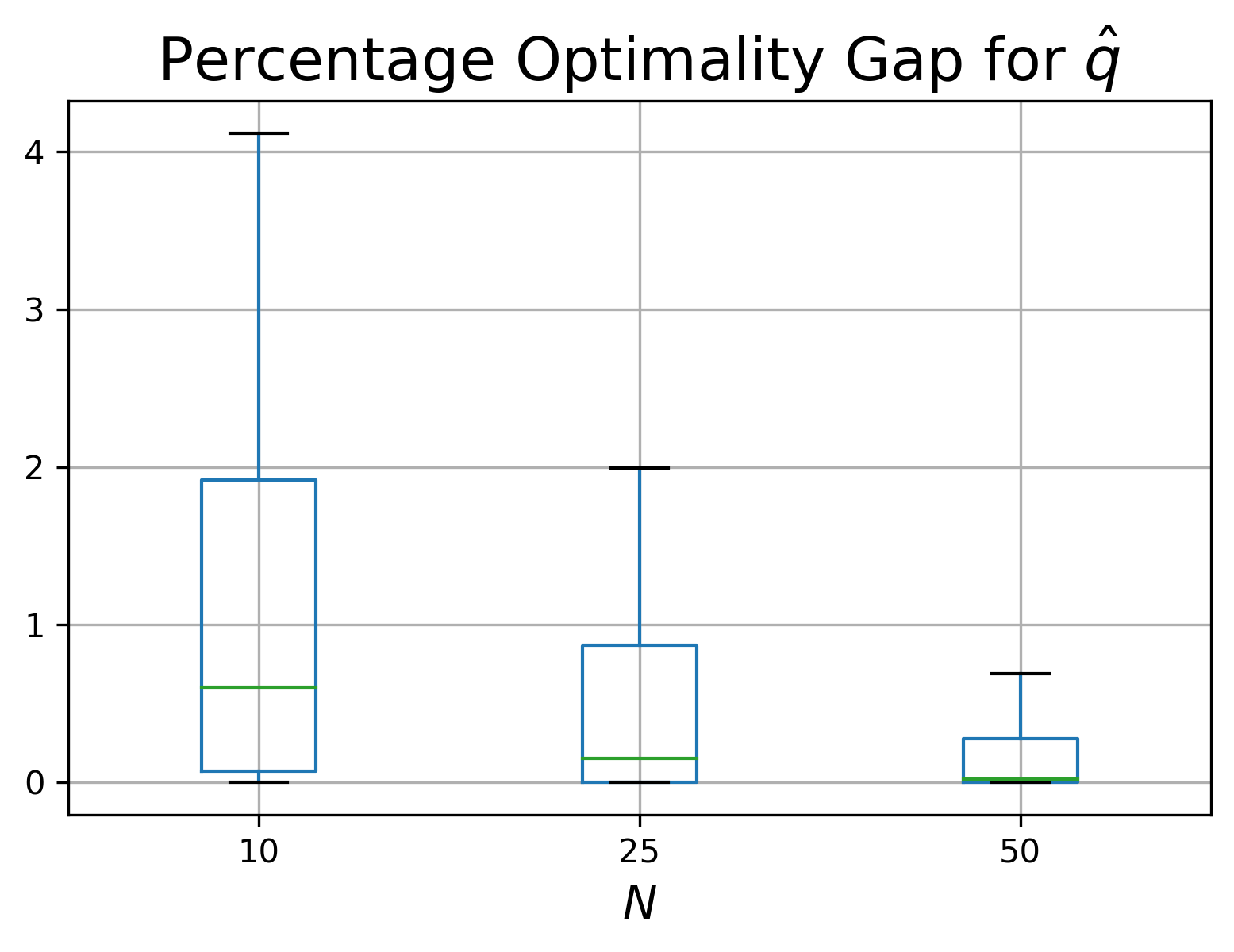}
        \caption{}
        \label{fig:FD_optgap_poisson}
    \end{subfigure}
    \begin{subfigure}[t]{0.45\textwidth}
        \centering
        \includegraphics[width=0.95\textwidth]{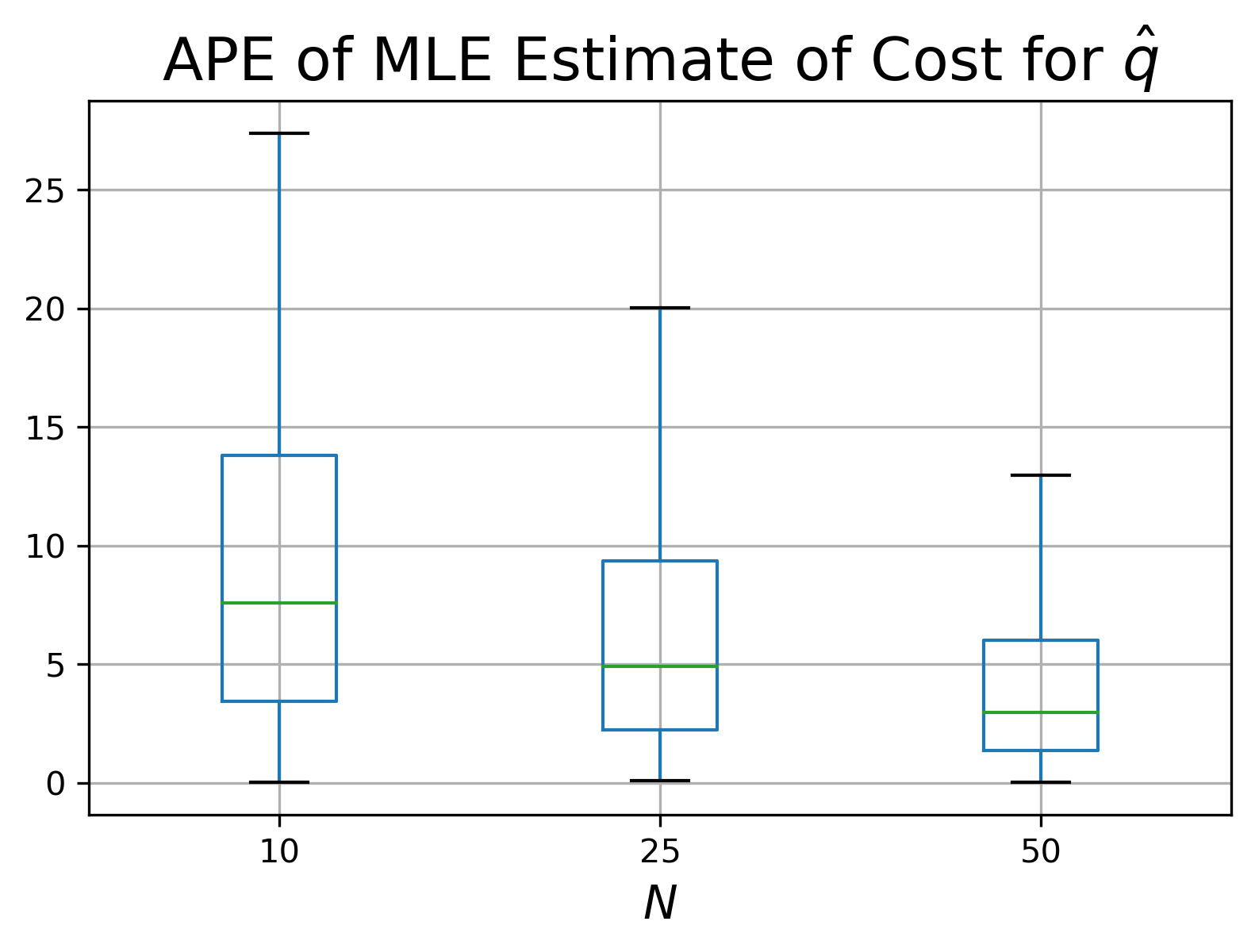}
        \caption{}
        \label{fig:FD_APE_poisson}
    \end{subfigure}
    \caption{Boxplots summarising performance of MLE approach. Figure 7(a) shows the optimality gap of $\hat{\bm{q}}$ as a solution to the model under $\bm{\lambda}^0$. Figure 7(b) shows the APE of $C_{\hat{\bm{\lambda}}}(\hat{\bm{q}})$ as an estimate of $C_{\bm{\lambda}^0}(\hat{\bm{q}})$.}
    \label{fig:FD_performance_poisson}
\end{figure}

Firstly, as for the normal distribution, we assess the effect of solving the newsvendor problem using the MLE distribution. We present the APE of the MLE's cost estimate along with the optimality gap of the MLE solution in Figure~\ref{fig:FD_performance_poisson}. This figure indicates that the MLE approach is close to optimal for Poisson demands, with solutions typically not being further than 4\% from optimal. However, the cost estimates are still not accurate. For $N=10$, we can see some APEs above 25\%. Even for $N=50$, we see APEs of almost 15\%. In addition, this plot does not show outliers. There were 19 instances where the APE exceeded 100\%, and 2 where it exceeded 10000\%.

In addition, we again find instances where the MLE predicts a profit for a decision that would actually attain a cost. This happened in 3 instances for the Poisson distribution, as opposed to 24 for the normal distribution. We summarise these 3 instances in Table~\ref{tab:FD_posneg}. This table shows the values of $\bm{\lambda}_0$ and $\hat{\bm{\lambda}}$ as well as the different costs associated with the MLE solution. This table shows that, in every instance, the DRO cost was positive. This suggests that the DRO objective would always show the newsvendor that the solution could lead to a cost. We also found that the MLE distribution predicted a cost when it would make a profit, although this only happened in two instances.

\begin{table}[htbp!]
    \centering
    \begin{tabular}{llrrr}
\toprule
$\lambda^0$ & $\hat{\lambda}$ &  FD Prediction &  True Cost &   DRO Cost \\
\midrule
   (10, 17) &    (8.8, 15.72) &      -28.62960 &  118.14383 &  177.11568 \\
    (9, 16) &    (8.12, 15.2) &      -25.11351 &   52.26802 &  130.41717 \\
   (15, 16) &  (14.24, 16.24) &      -16.40814 &  110.62268 &  264.08408 \\
\bottomrule
\end{tabular}

    \caption{Summary of instances where FD predicted profit but actually made a cost}
    \label{tab:FD_posneg}
\end{table}

Similarly to the normal case, we can confirm that the poor accuracy is not due to FD being a heuristic for the MPNVP model.  In the Poisson case, the two algorithms are even closer to one another than in the normal case.  We found that FD returned the same solution as PLA in 921 of the 972 instances considered. In addition, the solutions from the two algorithms had the same objective values in the other 51 instances. Therefore, the two algorithms have exactly the same performance for Poisson demands. Hence, the issues presented here are likely only due to the effects of treating the MLE as the true parameter.


In conclusion, the MLE approach performs better for Poisson demands than for normal ones, but it still has the same issues. In particular, it provides poor estimates of the cost associated with an order quantity. It is even capable of predicting a profit when a cost would be incurred, and vice versa. In addition, we confirmed that the DRO solution does not do this. In fact, in these instances the DRO cost is sometimes closer to reality than the MLE cost.

\subsubsection{Performance of DRO algorithms}

We now present the results for the 952 instances where we solved the DRO model with Poisson demands. The other 20 instances only had singleton ambiguity sets. We will summarise the run times of PLA and CS, along with CS's performance as a heuristic for the full model. Both algorithms were given a maximum time of 4 hours to build and solve the models. Firstly, we found that PLA timed out in 27 instances and CS did not time out in any instance. PLA took an average of 11 minutes and 12 seconds to either return a solution or time out, and its maximum time was 4 hours. CS took an average of 0.75 seconds, and its maximum time was 5.9 seconds. 
The 25th and 75th percentiles of PLA were 0.45 seconds and 16.84 seconds, respectively. For CS, these values are 0.22 and 0.61 seconds. Hence, again it is clear that PLA's average is affected by a small number of very large run times.

\begin{figure}[htbp!]
    \begin{subfigure}[t]{0.45\textwidth}
        \centering
        \includegraphics[width=0.97\textwidth]{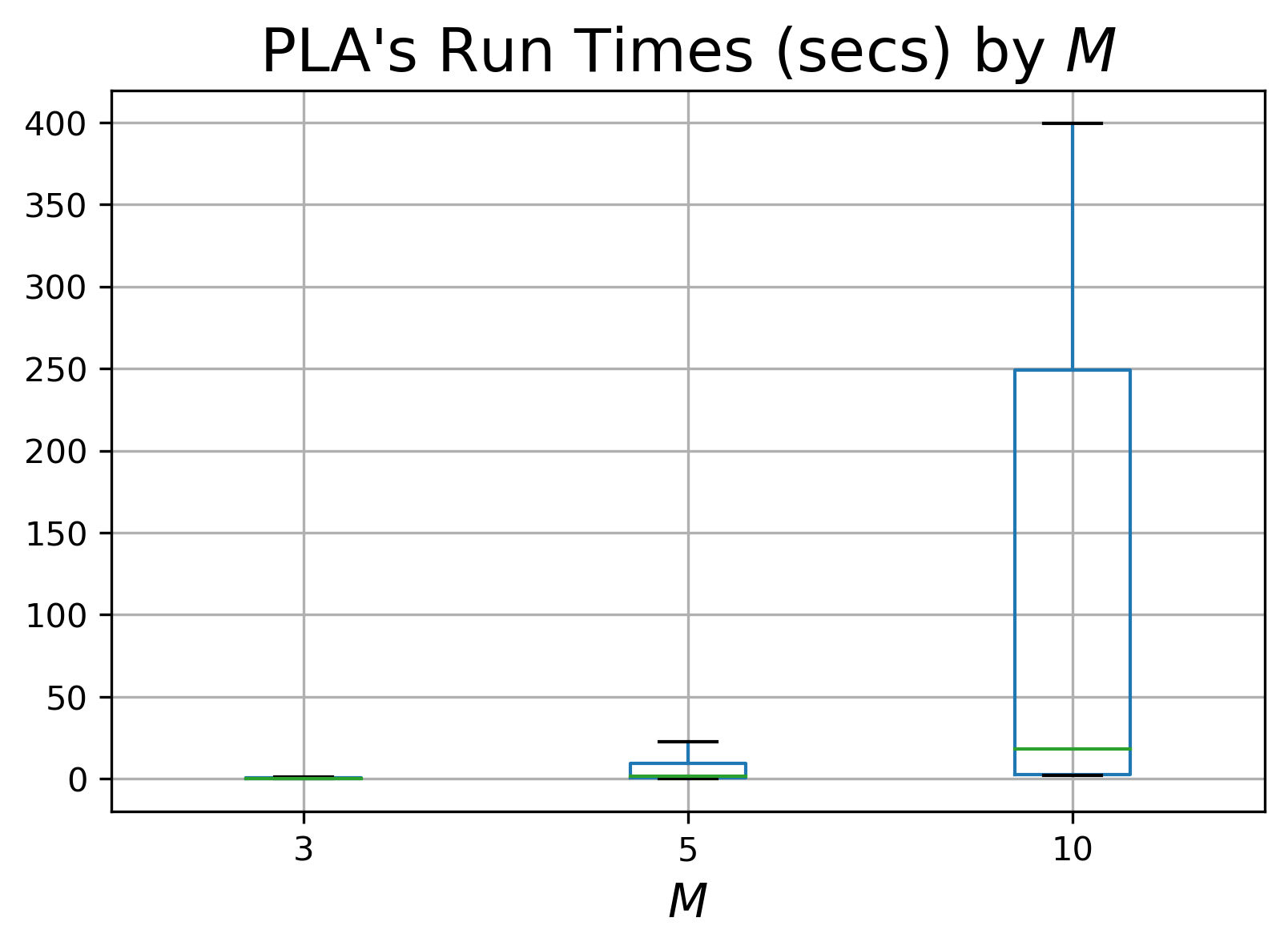}
        \caption{}
        \label{fig:poisson_PLA_tt_by_M}
    \end{subfigure}
    \begin{subfigure}[t]{0.45\textwidth}
        \centering
        \includegraphics[width=0.94\textwidth]{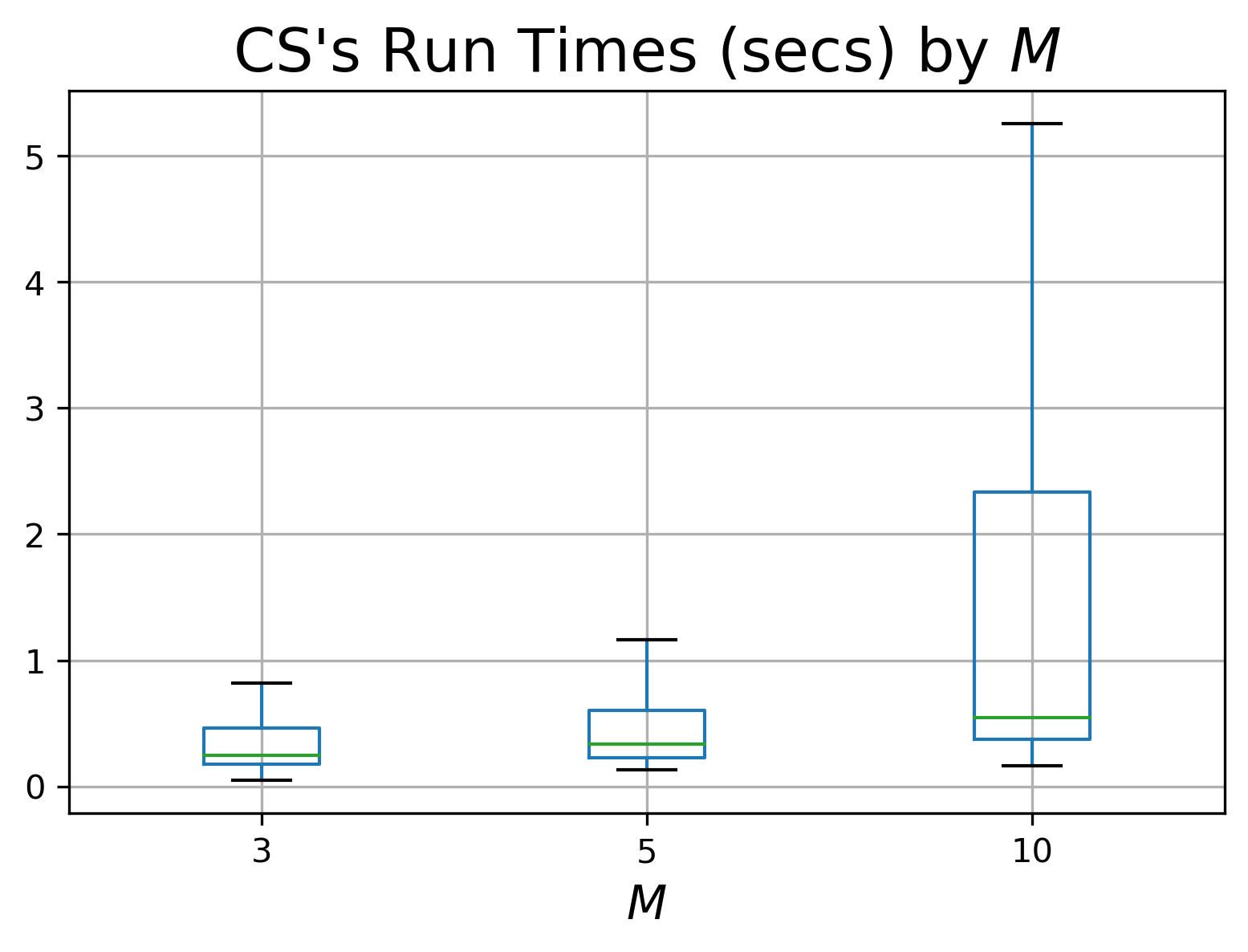}
        \caption{}
        \label{fig:poisson_PLA_tt_by_M2}
    \end{subfigure}
    \caption{Boxplots of run times, grouped by $M$}
    \label{fig:poisson_timestaken_by_M}
\end{figure}

The main variable affecting the solution times was, again, $M$. We present the run times by $M$ in Figure~\ref{fig:poisson_timestaken_by_M}. We can see PLA is greatly affected by $M$; it begins to take as long as 7 minutes for $M=10$. CS is also effected by $M$, but much less drastically. The difference is most noticeable for $M=10$, and this is the value at which both algorithms start to solve more slowly. However, CS did not take more than 6 seconds in any instance. The run times by $T$ reflect the results suggested by Figure~\ref{fig:poisson_timestaken_by_M}. This is because the reason why $M$ affects times so drastically is that it affects the size of the ambiguity set. 

\begin{figure}[htbp!]
    \begin{subfigure}[t]{0.45\textwidth}
        \centering
        \includegraphics[width=0.90\textwidth]{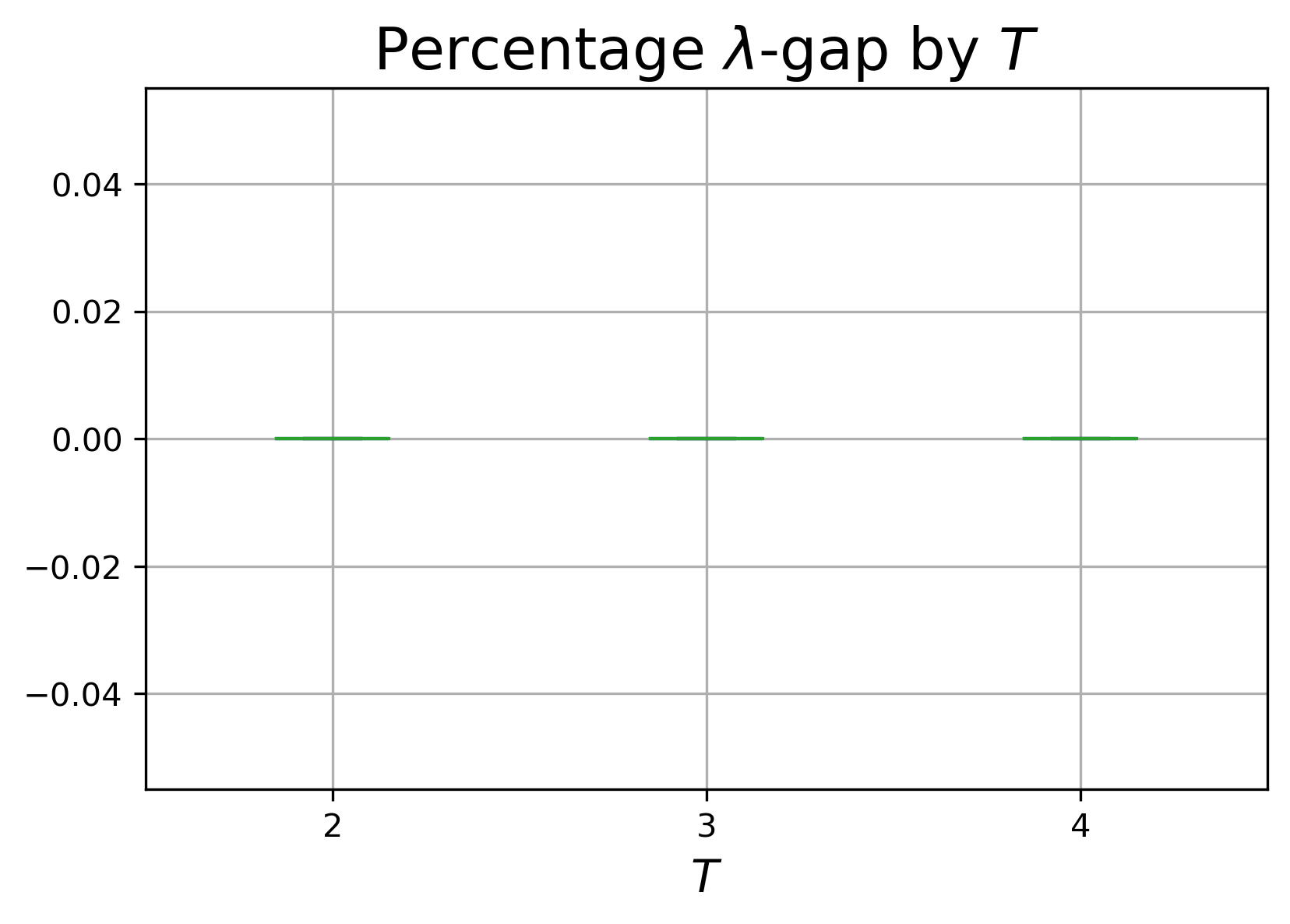}
        \caption{}
        \label{fig:poisson_omega_gap}
    \end{subfigure}
    \begin{subfigure}[t]{0.45\textwidth}
        \centering
        \includegraphics[width=0.95\textwidth]{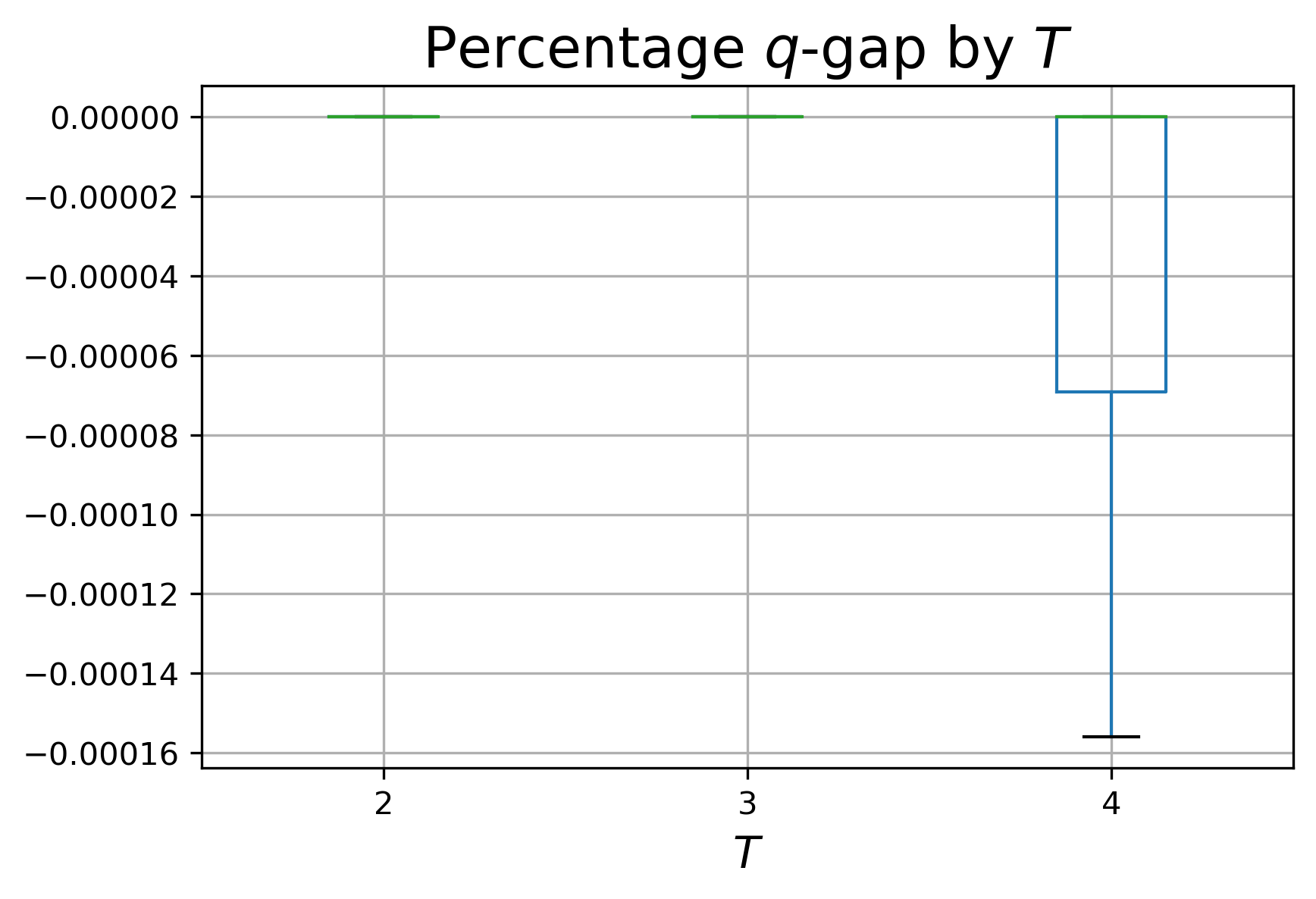}
        \caption{}
        \label{fig:poisson_q_gap}
    \end{subfigure}
    \caption{Summary of percentage gaps}
    \label{fig:poisson_gaps}
\end{figure}

We now present the optimality of CS with respect to $\bm{\lambda}$ and $\bm{q}$. Firstly, we present the number of times each algorithm selected the worst-case parameter for its chosen $\bm{q}$ in Table~\ref{tab:poisson_worstcase}. This table shows that CS selected the worst-case $\bm{\lambda}$ for its $\bm{q}$ in 94.2\% of instances. We will show that its other choices were very close to the worst-case, however. PLA selected the worst-case in every instance. The percentage $\bm{\lambda}$- and $\bm{q}$-gaps are defined as they were for the normal experiments in Section~\ref{sec:normal_CSvsPLA}, and we present summaries of these gaps in Figure~\ref{fig:poisson_gaps}. Firstly, Figure~\ref{fig:poisson_omega_gap} reflects the fact that CS's chosen $\bm{\lambda}$ was very close to the true worst-case. The overall average $\bm{\lambda}$-gap was 0.09\%. The fact that CS did not always select the true worst-case is due to the fact that $\Omega^{\text{ext}}$ is constructed by combining univariate results. Despite the fact that the expected cost can always be increased by setting one $\lambda_t$ to its minimum or maximum, this may not be possible without changing another $\lambda_j$ given~\eqref{eq:pois_conf_set}. The cost may improve as a result of these two changes. Hence, the worst-case sometimes did not have any parameters at their maximum or minimum.

 \begin{table}[htbp!]
    \centering
    \begin{tabular}{lrlll}
\toprule
 $T$ &  CS No. & CS Worst-case Count (\%) & PLA No. & PLA Worst-case Count(\%) \\
\midrule
 All &     952 &            880 (92.44\%) &     925 &            925 (100.0\%) \\
   2 &     311 &            311 (100.0\%) &     311 &            311 (100.0\%) \\
   3 &     319 &            311 (97.49\%) &     319 &            319 (100.0\%) \\
   4 &     322 &            258 (80.12\%) &     295 &            295 (100.0\%) \\
\bottomrule
\end{tabular}

    \caption{Number of instances in which algorithms selected true worst case parameters.}   
    \label{tab:poisson_worstcase}
\end{table}

Secondly, Figure~\ref{fig:poisson_q_gap} shows that CS typically selected a $\bm{q}$ with a worst-case cost of less than $1.6\times 10^{-4}\%$ higher than PLA's. In fact, CS's $\bm{q}$ had the same worst-case expected cost as PLA's in 796 (86\%) of the 925 instances in which PLA did not time out. The overall average $\bm{q}$-gap was $-0.02\%$,  and the minimum value of this gap was $-3.14\%$. CS did not encounter much suboptimality until $T=4$. In these instances, its suboptimality was likely due to either its selection of the incorrect worst-case parameter, solver-specific issues like Gurobi's pre-solve model reductions, or different starting points for its MIP algorithms.

From the results shown here, we can conclude two main points. Firstly, CS again solved much faster than PLA. Secondly, it was very near optimal if not optimal, in general.  It was capable of selecting the worst-case parameter for its chosen $\bm{q}$ in the large majority of instances, and this $\bm{q}$ was often optimal.
\section{Conclusions and further work}

In this paper, we considered a static, multi-period newsvendor problem under a budget constraint. We first developed an algorithm to solve the problem under a known demand distribution. This was done by adapting an algorithm from the literature that was used to solve a multi-product budgeted newsvendor problem. It entails first finding a solution to the stationarity KKT condition, and then optimising the Lagrange multiplier until the solution becomes feasible. We tested this algorithm against three benchmarks, and found that it is faster than all 3 of them and provides comparable solutions.

Following this, we studied parametric distributionally robust optimisation models for this newsvendor problem. We assumed that the newsvendor knows the family of distributions in which the true demand distribution lies. We showed how maximum likelihood estimation can be used to construct a confidence-based ambiguity set for the true parameters, and then how to build and solve the resulting distributionally robust model for normal and Poisson demands. Due to the size of the discrete ambiguity sets used, we developed a decomposition algorithm named CS for each family of distributions. We did so by obtaining theoretical results about the behaviour of the expected cost as a function of the distribution's parameters, in order to characterise the worst-case distribution. We showed that the CS algorithm gave very close to optimal solutions in a matter of seconds, for both families. We also assessed the performance of the method that solved the fixed distribution problem assuming the maximum likelihood estimate is the true parameter. We found that, while the solutions this yielded were not far from optimal, the cost estimates from this model were far from accurate. In fact, we found that there were cases where the maximum likelihood estimate suggested a profit would be made, when in reality a cost would be incurred. We also confirmed that, in these instances, the distributionally robust model would in fact make the newsvendor aware that their order could result in a cost. Hence, using distributionally robust models avoids these issues. 

There are a number of areas for future work that would prove interesting. The most immediate one is to extend our work to the case of dependent demands. It is common in real-life scenarios that demand for a product over different days will be correlated. However, this creates a much more complex model. We would like to study how this could be done in future, for example by specifying the full covariance matrix in the multivariate normal case. For Poisson demands, it is less clear how this could be done. Another area for future research is to study other demand distributions, such as exponential or gamma. The two distributions considered here are common for demand in the newsvendor problem, but there are many others that might be used. Finally, we would like to improve the process used to compute the ambiguity sets. Since CS has been shown to be very fast, the only real bottleneck of our approach is the time taken to compute ambiguity sets. It would also be interesting to study how this process can be further improved, potentially with parallelisation and/or faster programming languages. 
\Large{\textbf{Acknowledgements}}

\normalsize
We would like to acknowledge the support of the Engineering and Physical
Sciences Research Council funded (EP/L015692/1) STOR-i Centre for Doctoral
Training. 

\bibliographystyle{apalike}
\bibliography{bib}

\newpage
\begin{appendices}
    \setcounter{equation}{19}

\section{Proof of Theorem~3.1}\label{sec:KKT_proof}
\begin{proof}
To prove this theorem, we need to differentiate the objective function. To do so, firstly we evaluate the expected values. For ease of notation, let $Q_t = \sum_{k=1}^t q_k$ and $Y_t = \sum_{k=1}^t X_k$. Then, $\tilde{F}_t$ and $\tilde{f}_t$ are the CDF and PDF of $Y_t$. In addition, we have:
\begin{align*}
    \E_F\l[I_t^+\r] &= \E_F\l[\max(I_t, 0)\r] \\
    &= \P(I_t \ge 0)\E_F(I_t | I_t \ge 0) \\
    &= \P(Y_t \le Q_t)\E_F\l(Q_t - Y_t | Y_t \le Q_t\r)\\
    &= \tilde{F}_t\l(Q_t\r)\l[Q_t - \E_{\tilde{F}_t}\l(Y_t | Y_t \le Q_t\r)\r]
\end{align*}
where $\tilde{F}_t$ and $\tilde{f}_t$ are the CDF and PDF of $Y_t$. 
Now note that we have:
\begin{equation*}
    \pderiv{\tilde{F}_t\l(Q_t\r)}{q_j} = \begin{cases}
        \tilde{f}_t\l(Q_t\r) &\text{ if } j \le t,\\
        0 &\text{ otherwise}
    \end{cases}.
\end{equation*}
To differentiate the expected values we need the following. Firstly, the random variable $Y_t | Y_t \le Q_t$ is $Y_t$ truncated above at $Q_t$. Therefore, the PDF of this random variable is:
\begin{equation*}
    \bar{f}_t(u) = \begin{cases}
        \frac{\tilde{f}_t(u)}{\tilde{F}_t(Q_t)} &\text{ if } u \le Q_t\\
        0 &\text{ otherwise}
    \end{cases}
\end{equation*}
and therefore:
\begin{equation*}
    \E_{\tilde{F}_t}\l(Y_t | Y_t \le Q_t\r) = \int_{-\infty}^{\infty} u \bar{f}_t(u) \text{d}u = \frac{\int_{-\infty}^{Q_t} u\tilde{f}_t(u) \text{d}u}{\tilde{F}_t\l(Q_t\r)}.
\end{equation*}
This means that:
\begin{align*}
    \E_F[I^+_t] &= \tilde{F}_t\l(Q_t\r)\l[Q_t - \frac{\int_{-\infty}^{Q_t} u\tilde{f}_t(u) \text{d}u}{\tilde{F}_t\l(Q_t\r)}\r]\\
    &= \tilde{F}_t\l(Q_t\r)Q_t - \int_{-\infty}^{Q_t} u\tilde{f}_t(u) \text{d}u.
\end{align*}
Applying the second fundamental theorem of calculus and the chain rule we find that:
\begin{align*}
    \pderiv{}{q_j}\int_{-\infty}^{Q_t} u\tilde{f}_t(u) \text{d}u &= \pderiv{}{q_j}\l(\l[\int u\tilde{f}_t(u) \text{d}u\r]\Bigg|_{u=Q_t} -  \lim_{v\to-\infty}\l[\int u\tilde{f}_t(u) \text{d}u\r]\Bigg|_{u=v}\r) \\
    &= \pderiv{}{q_j}\l[\int u\tilde{f}_t(u) \text{d}u\r]\Bigg|_{u=Q_t}\\
    &=\l(Q_t \times \tilde{f}_t\l(Q_t\r)\r) \times \pderiv{}{q_j}Q_t \\
    &= Q_t \times \tilde{f}_t\l(Q_t\r)
\end{align*}
for $j \le t$ and 0 otherwise. Therefore:
\begin{align*}
    \pderiv{}{q_j}\E_F[I^+_t] &= \tilde{f}_t\l(Q_t\r) Q_t + \tilde{F}_t\l(Q_t\r) - \tilde{f}_t\l(Q_t\r) Q_t \\
    &= \tilde{F}_t\l(Q_t\r).
\end{align*}
for $j \le t$ and 0 otherwise.  Similarly, we can find that:
\begin{align*}
    \E_F\l[I_t^-\r] &= \E_F\l[\max(-I_t, 0)\r] \\
    &= -\P(I_t \le 0)\E_F(I_t | I_t \le 0) \\
    &= \l(\tilde{F}_t\l(Q_t\r) - 1\r)Q_t + \int_{Q_t}^\infty u\tilde{f}_t(u) \text{d}u.
\end{align*}
Applying the same calculus as before, we find:
\begin{align*}
    \pderiv{}{q_j}\int_{Q_t}^\infty u\tilde{f}_t(u) \text{d}u 
    &= - Q_t \times \tilde{f}_t\l(Q_t\r).
\end{align*}
Hence:
\begin{align*}
    \pderiv{}{q_j} \E_F[I^-_t] &= \tilde{f}_t\l(Q_t\r)Q_t + \l(\tilde{F}_t\l(Q_t\r) - 1\r) - Q_t\times \tilde{f}_t\l(Q_t\r)\\
    &= \tilde{F}_t\l(Q_t\r) - 1
\end{align*}
for $j \le t$, 0 otherwise. Hence, we have:
\begin{align*}
    \pderiv{}{q_j} C_F(\bm{q}) &= p\l(\tilde{F}_T\l(Q_T\r) - 1\r) + \sum_{t=j}^T\l[(h + b)\tilde{F}_t\l(Q_t\r) - b\r] + w_j \\
    &= \sum_{t=j}^T\l[(h + b + p\mathds{1}\{t=T\})\tilde{F}_t\l(Q_t\r) - b\r] + w_j - p.
\end{align*}
The Lagrangian of MPNVP without considering the non-negativity constraints is:
\begin{equation*}
    L(\bm{q}, \nu) = C_F(\bm{q}) + \nu \l(\sum_{t=1}^T w_t q_t - W\r).
\end{equation*}
The KKT conditions are therefore:
\begin{align*}
    \pderiv{}{q_j} C_F(\bm{q}^*) + \nu w_j &= 0 \fa j = 1,\dots,T,\\
    \sum_{t=1}^T w_t q^*_t &\le W, \\
    \nu\l(\sum_{t=1}^T w_t q^*_t - W\r) &= 0 , \\
    \nu &\ge 0 
\end{align*}
Note that this means that the stationarity condition is equivalent to:
\begin{equation*}
    \sum_{t=j}^T\l[(h + b + p\mathds{1}\{t=T\})\tilde{F}_t\l(Q^*_t\r) - b\r] = p - (1+\nu)w_j \fa j = 1,\dots,T.
\end{equation*}
Now, assume that $\bm{q}^*$ is as given in~(5). Then, by condition 1.\ and 2.\ in Theorem~3.1, we have:
\begin{align}
    \frac{b - (1+\nu)(w_t - w_{t+1})}{h+b} &\ge 0 \fa t \in \{1,\dots,T-1\},\label{eq:LHS1} \\
    \frac{b + p - (1+\nu)w_t}{h+p+b} \ge 0, \label{eq:LHS2} 
\end{align}
and since $w_t - w_{t+1} \ge 0$, we also have:
\begin{align*}
    &\frac{b - (1+\nu)(w_t - w_{t+1})}{h+b} \le \frac{b}{h+b} \le 1 \fa t \in \{1,\dots,T-1\}, \\
    &\frac{b + p - (1+\nu)w_T}{h+p+b} \le \frac{b+p}{h+b+p} \le 1. 
\end{align*}
Hence, the left hand sides of~\eqref{eq:LHS1} and~\eqref{eq:LHS2} lie between 0 and 1, meaning the inverse CDFs are defined at these points. Hence, $\bm{q}^*$ is well-defined. Now, note that we have:
\begin{equation*}\label{eq:sum_of_qs}
    \sum_{k=1}^t q^*_k = \begin{cases}
        \tilde{F}^{-1}_t\l(\frac{b - (1+\nu)(w_t - w_{t+1})}{h+b}\r) &\text{ if } t < T,\\
        \tilde{F}^{-1}_T\l(\frac{b - (1+\nu)w_T + p}{h+b+p}\r) &\text{ if } t= T.
    \end{cases}
\end{equation*}
Substituting~(5) in to the Lagrangian's derivative, we get:
\begin{align*}
    \pderiv{}{q_j} L(\bm{q}^*, \nu) =& \sum_{t=j}^{T-1} \l[b - (1+\nu)(w_t - w_{t+1})- b\r] + \l[b - (1+\nu)w_T + p - b\r] - p + (1+\nu)w_j\\
    =&\ \sum_{t=j}^{T-1}\l[(1+\nu)(w_{t+1} - w_t )\r]  - (1+\nu)w_T + (1+\nu)w_j  \\
    =& \  (1+\nu)(w_T - w_j) - (1+\nu)(w_T - w_j)\\
    =& \ 0.
\end{align*}
Hence, the stationarity condition is satisfied.
\end{proof}

\section{Details on FD and its benchmarks}\label{sec:FD_details}

\subsection{Deriving FD}\label{sec:derive_FD}

The general idea of our algorithm is similar to the multi-product algorithm of~\cite{ALFARES2005465}. A high-level outline of their algorithm is as follows:
\begin{enumerate}
    \item Set $\nu = 0$ (i.e.\ no budget constraint) and solve using KKT conditions, setting negative order quantities to zero.
    \item If budget constraint is met then go to step 5. Else go to step 3.
    \item Starting from $\nu = 0$ increase $\nu$ until the first occurrence of either of the following:
    \begin{enumerate}
        \item An order quantity becomes negative. Go to step 4.
        \item Budget constraint met. Go to step 5.
    \end{enumerate}
    \item Set the negative order quantity to zero, remove from the set of products and go to step 1.
    \item Return solution.
\end{enumerate}
In essence, we solve with no constraints, and increase the Lagrange multiplier $\nu$ (e.g., by line search) until either the vector of order quantities becomes feasible, or one order quantity  becomes negative. If an order quantity becomes negative, we remove the corresponding product from consideration and start again. The main difference between our model and that of~\cite{ALFARES2005465} is the structure of the optimal KKT solution that results from the inventory tracking. Their solution is of the form:
\begin{equation*}
    x_\tau = F^{-1}(P_\tau) \fa \tau \in \{1,\dots,N\},
\end{equation*}
for some CDF $F$ and probabilities $P_\tau$. This means that the optimal unconstrained orders for different items are independent. However, our solution is of the form:
\begin{equation}\label{eq:q_t_simple}
    q^*_t = \tilde{F}^{-1}_t(P_t) - \tilde{F}^{-1}_{t-1}(P_{t-1}),
\end{equation}
where the subtraction of $\tilde{F}^{-1}_{t-1}(P_{t-1})$ is done to ensure that 
\begin{equation}\label{eq:important_relation}
    \sum_{k=1}^t q^*_k = \tilde{F}^{-1}_t(P_t).
\end{equation}
This relation holding true is what ensures that our KKT conditions are met. In the algorithm of~\cite{ALFARES2005465}, when an item is removed and set to zero, this is because its marginal profit has become negative. Therefore, if we continue to increase $\nu$ then this value will become negative and stay negative. We can apply similar logic to our problem, with one slight difference. When Alfares and Elmorra set one variable to zero, the other variables can be found via their respective formulas and the corresponding derivatives will still be zero. However, if we set $q^*_{t-1} = 0$ then the value of $q^*_t$ from~\eqref{eq:q_t_simple} will not give~\eqref{eq:important_relation}. In fact, $q^*_{t-1}$ being zero means that: 
\begin{align*}
    \sum_{k=1}^t q^*_k &= \sum_{k=1}^{t-2} q^*_k + q^*_t\\
    &= \tilde{F}^{-1}_{t-2}(P_{t-2}) + ( \tilde{F}^{-1}_t(P_t) - \tilde{F}^{-1}_{t-1}(P_{t-1})).
\end{align*}
In order to satisfy~\eqref{eq:important_relation}, the optimal value of $q^*_t$ becomes $\tilde{F}^{-1}_t(P_t) - \tilde{F}^{-1}_{t-2}(P_{t-2})$. Hence, when we remove a day from the set of days, we need to adjust the remaining solutions. Suppose we have a set $\m{T}^0$ of days that have been set to zero permanently. Then, when re-calculating our next solution we use $\bm{q} = \tilde{\bm{q}}^*$, where $\tilde{\bm{q}}^*$ is defined by~(6). This solution subtracts the sum of all previous orders from each $\tilde{q}_t$ to ensure that~\eqref{eq:important_relation} holds. 

\subsection{FD's line search algorithm}\label{sec:FD_LS}

The line search algorithm inside FD works as follows:
\begin{enumerate}
    \item Initialise $\nu = 0$, select scaling parameter $\tau \in (0,1)$, stepsize $\delta$, minimum stepsize $\delta^{\min}$, budget tolerance $\varepsilon$.
    \item Set $\bm{q}' = \bm{q} = \tilde{\bm{q}}^*(\nu)$, if $\sum_{t=1}^T q_t \le W$ then set $\texttt{done} = True$, else set $\texttt{done} =  False$.
    \item While not \texttt{done}:
    \begin{enumerate}
        \item While $\nu + \delta > \nu^{\text{UB}}$:
        \begin{enumerate}
            \item $\delta = \tau \delta$.
            \item $\delta^{\min} = \min\{\delta, \delta^{\min}\}$.
        \end{enumerate}
        \item Set $\nu' = \nu + \delta$, $\bm{q}' =  \tilde{\bm{q}}^*(\nu')$:
        \begin{enumerate}
            \item If either $\sum_{t=1}^T w_t q'_t > W$ or $\delta = \delta^{\min}$, then set $\nu = \nu'$.
            \item Otherwise, define $\delta' = \tau \delta$. If $\delta' \le \delta^{\min}$ then set $\delta^{\min} = \delta$. Otherwise set $\delta = \delta'$ and $\bm{q}' = \bm{q}$.
        \end{enumerate}
        \item If either $\l(\big| \sum_{t=1}^T q'_t - W\big| \le \varepsilon \text{ and } q'_t \ge 0 \fa t \in \m{T}\r)$ or $\min{\bm{q}'} < 0$, set $\texttt{done} = True$.
    \end{enumerate}
\end{enumerate}
The idea of the algorithm is as follows. We first evaluate $\tilde{\bm{q}}^*(0)$ and if this solution is feasible, we return this solution. If not, we know that $\nu$ needs to be increased in order to either make the solution feasible or make an order quantity negative. Step 3 therefore increases $\nu$ until one of these situations occurs. 

Step 3(a) reduces $\delta$ until $\nu' = \nu + \delta < \nu^{\text{UB}}$ if this does not hold already. This is to ensure that we can get arbitrarily close to $\nu^{\text{UB}}$ in order to find the first negative or feasible order. Then, in step 3(b), we evaluate the order quantities for $\nu' = \nu + \delta$. If the new solution $\bm{q}'$ does not satisfy the budget constraint, then we move to 3(c) and check whether any orders are negative. If there is a negative order, we know that this occurs before the budget constraint is met, so the algorithm ends. If there are no negative orders, we know that $\nu$ has to be increased further and we go to 3(a). If $\sum_{t=1}^T w_t q'_t \le W$ in 3(b)i, then we have found a solution that satisfies the budget constraint. Therefore, we need to decide if this is the first such solution, or if a smaller $\nu$ can be found that also achieves this. If $\delta = \delta^{\min}$ (which is $10^{-100}$ in our experiments), then it is very unlikely that this is not the first solution that satisfies the budget constraint. Hence, we move to 3(c). If, on the other hand, $\delta > \delta^{\min}$, then it is possible that we have skipped over the first solution that satisfies the budget constraint. Hence, we reduce $\delta$ in step 3(b)ii. This means that, if we need to go back to 3(a), our next solution will be closer to the first solution that satisfies the budget constraint.

\subsection{FD's benchmark algorithms}\label{sec:scipy_algs}

The algorithms that we compare FD with are described in more detail as follows:
\begin{enumerate}
    \item SLSQP: Sequential Least Squares Quadratic Programming~\citep{SLSQP}. This is an algorithm that can be found in Python's baseline package for scientific computing, Scipy~\citep{2020SciPy-NMeth}. The algorithm is a line search algorithm, where at each iteration the next search direction is chosen by solving a quadratic programming relaxation of the problem's Lagrangian. The algorithm ends when the gradient of the true objective function is sufficiently small at the candidate solution.
    \item TC: Trust Constraint~\citep{trust_constr}. This is a trust region algorithm that can be found in Scipy. Trust region algorithms are similar to line search algorithms, but where the next solution is generated by solving an approximate version of the model over a region where the approximation is ``trusted''. In each iteration, the ratio of the improvement in the true objective function versus the approximate objective function is used to determine whether the region should be enlarged or reduced. The algorithm ends when the gradient of the true objective function is sufficiently small at the candidate solution.
    \item PLA: Piecewise Linear Approximation. This entails using a piecewise linear approximation to the objective function, allowing it to be solved by Gurobi~\citep{gurobi}. 
\end{enumerate}

\section{Appendix for normally distributed demands}\label{sec:normal_proofs}

\subsection{Proof of Lemma~4.1}\label{sec:normal_obj_proof}

\begin{proof}
Again, let $Q_t = \sum_{k=1}^t q_k$ and $Y_t = \sum_{k=1}^t X_k$. Since the $X_t$'s are independent normal, random variables, we have:
\begin{equation*}\label{eq:X_dist}
    Y_t \sim \m{N}\l(\sum_{k=1}^t \mu_k, \sum_{k=1}^t \sigma^2_k\r),
\end{equation*}
Furthermore, since $I_t = Q_t - Y_t$, we have:
\begin{equation*}
    I_t \sim \m{N}\l(Q_t - \sum_{k=1}^t \mu_k, \sum_{k=1}^t \sigma^2_k\r) \fa t = 1,\dots,T,
\end{equation*}
Furthermore, $I_t | I_t \ge 0$ is a normal random variable truncated below at 0. Therefore, by~\cite{truncated_meanvar}, we have:
\begin{align*}
    \E_{\bm{\omega}}(I_t | I_t \ge 0) &= \l(Q_t - \sum_{k=1}^t \mu_k\r)  + \frac{\phi(\beta_t)}{1-\Phi(\beta_t)}\sqrt{\sum_{k=1}^t \sigma^2_t},
\end{align*}
where 
$$\beta_t =\frac{0 - \E_{\bm{\omega}}(I_t)}{\sqrt{\text{Var}_{\bm{\omega}}(I_t)}} =\frac{ \sum_{k=1}^t \mu_k - Q_t}{\sqrt{\sum_{k=1}^t
\sigma_k^2}}.$$
Therefore, we have:
\begin{align*}
    \E_{\bm{\omega}}[I^+_t]  &= \P(I_t \ge 0)\E_{\bm{\omega}}[I_t | I_t \ge 0] \\
    &= \P(Y_t \le Q_t)\E_{\bm{\omega}}[I_t | I_t \ge 0]\\
    &= \Phi(-\beta_t) \l(\l(Q_t - \sum_{k=1}^t \mu_k\r)  + \frac{\phi(\beta_t)}{1-\Phi(\beta_t)}\sqrt{\sum_{k=1}^t \sigma^2_t}\r) \\
    &= \Phi(-\beta_t)\l(Q_t - \sum_{k=1}^t \mu_k\r)  + \phi(\beta_t)\sqrt{\sum_{k=1}^t \sigma_k^2},
\end{align*}
Similarly, we have:
\begin{align*}
    \E_{\bm{\omega}}[I^-_t] &= \E_{\bm{\omega}}[\max\{-I_t, 0\}]\\
    &= \P(-I_t \ge 0)\E_{\bm{\omega}}[-I_t | -I_t \ge 0] \\
    &= \Phi(\beta_t) \l(-\l(Q_t - \sum_{k=1}^t\mu_k\r) + \sqrt{\sum_{k=1}^t \sigma_k^2} \frac{\phi(-\beta_t)}{1-\Phi(-\beta_t)}\r) \\
    &= \Phi(\beta_t)\l(\sum_{k=1}^t\mu_k - Q_t\r) + \phi(\beta_t)\sqrt{\sum_{k=1}^t \sigma_k^2},
\end{align*}
This means we can write the objective function as:
\begin{align*}
    C_{\bm{\omega}}(\bm{q}) =& \sum_{t=1}^T\l\{
    h\l[\Phi(-\beta_t)\l(Q_t - \sum_{k=1}^t\mu_k \r) +
    \phi(\beta_t)\sqrt{\sum_{k=1}^t \sigma_k^2}\r] \r.\\ 
    &\l.+
    (b + p\mathds{1}\{t= T\})\l[\Phi(\beta_t)\l(\sum_{k=1}^t\mu_k - Q_t\r) +
    \phi(\beta_t)\sqrt{\sum_{k=1}^t \sigma_k^2}\r] + w_t q_t - p \mu_t\r\}\\
    =& \sum_{t=1}^T\l\{(h(\Phi(\beta_t) - 1) +  (b +  p\mathds{1}\{t=
    T\})\Phi(\beta_t))\l(\sum_{k=1}^t\mu_k - Q_t\r)\r. \\
    &\l.+ (h + b +
    p\mathds{1}\{t= T\})\phi(\beta_t)\sqrt{\sum_{k=1}^t\sigma^2_k} +
    w_t q_t - p \mu_t\r\} \\
    =& \sum_{t=1}^T\l\{((b +  p\mathds{1}\{t=
    T\} + h)\Phi(\beta_t) - h)\l(\sum_{k=1}^t\mu_k - Q_t\r)\r. \\
    & \l.+ (h + b +
    p\mathds{1}\{t= T\})\phi(\beta_t)\sqrt{\sum_{k=1}^t\sigma^2_k} +
    w_t q_t- p \mu_t\r\} 
\end{align*}
or
\begin{equation*}
    C_{\bm{\omega}}(\bm{q}) \hiderel{=}
     \sum_{t=1}^T\l\{(a_t \Phi(\beta_t) - h)\sum_{k=1}^t\mu_k - Q_t
    + a_t \phi(\beta_t)\sqrt{\sum_{k=1}^t \sigma_k^2} +
    q_t w_t- p \mu_t\r\}.
\end{equation*}
where $\mathds{1}\{t=T\}$ is 1 if $t=T$ and 0 otherwise, and $a_t = b +
p\mathds{1}\{t=T\} + h$.
\end{proof}

\subsection{Proof of Proposition~4.2}\label{sec:normal_conf_proof}
\begin{proof}
The likelihood function for the data is:
\begin{align*}
    \m{L}(\bm{x}|\bm{\omega}) &= \prod_{n=1}^N f_{\bm{X}}(\bm{x}| \bm{\omega}) \\
    &= \prod_{n=1}^N \prod_{t=1}^T \frac{1}{\sigma_t\sqrt{2\pi}}\exp\l(-\frac{1}{2\sigma^2_t}(x^n_t - \mu_t)^2\r)\\
    &\propto \l(\prod_{t=1}^T(\sigma_t)^{-1}\r)^N \times \exp\l(-\sum_{n=1}^N\sum_{t=1}^T\frac{1}{2\sigma^2_t}(x^n_t - \mu_t)^2\r).
\end{align*}
The log-likelihood is:
\begin{align*}
    \ell(\bm{\omega}) &= -N \sum_{t=1}^T \log(\sigma_t) - \sum_{n=1}^N\sum_{t=1}^T \frac{1}{2\sigma^2_t}(x^n_t - \mu_t)^2.
\end{align*}
Taking derivatives, we find:
\begin{align*}
    \pderiv{\ell(\bm{\omega})}{\mu_t} &= \frac{1}{\sigma^2_t}\sum_{n=1}^N (x^n_t - \mu_t) \\
    \pderiv{\ell(\bm{\omega})}{\sigma_t} &= -\frac{N}{\sigma_t} + \frac{1}{\sigma^3_t}\sum_{n=1}^N (x^n_t - \mu_t)^2 \\
    \frac{\partial^2 \ell(\bm{\omega})}{\partial \mu_t \partial \sigma_t} &= -\frac{2}{\sigma^3_t}\sum_{n=1}^N(x^n_t - \mu_t) \\
    \frac{\partial^2 \ell(\bm{\omega})}{\partial \mu^2_t} &= -\frac{N}{\sigma^2_t} \\
    \frac{\partial^2 \ell(\bm{\omega})}{\partial \sigma_t \partial \mu_t} &= -\frac{2}{\sigma^3_t}\sum_{n=1}^N(x^n_t - \mu_t) \\
    \frac{\partial^2 \ell(\bm{\omega})}{\partial \sigma^2_t} &= \frac{N}{\sigma^2_t} - \frac{3}{\sigma^4_t}\sum_{n=1}^N(x^n_t - \mu_t)^2,
\end{align*}
and all other second partial derivatives are zero. As expected, this yields the MLEs:
\begin{align*}
    \hat{\mu}_t &= \frac{1}{N}\sum_{n=1}^N x^n_t\\
    \hat{\sigma}_t &= \sqrt{\frac{1}{N} \sum_{n=1}^N (x^n_t - \hat{\mu}_t)^2}.
\end{align*}
The expected values of the non-zero second derivatives are given by:
\begin{align*}
    \E_{\bm{\omega}}\l(\frac{\partial^2 \ell(\bm{\omega})}{\partial \mu_t \partial \sigma_t}\r) &= 0 \\
    \E_{\bm{\omega}}\l(\frac{\partial^2 \ell(\bm{\omega})}{\partial \mu^2_t}\r) &= -\frac{N}{\sigma^2_t} \\
    \E_{\bm{\omega}}\l(\frac{\partial^2 \ell(\bm{\omega})}{\partial \sigma_t \partial \mu_t}\r) &= 0 \\
    \E_{\bm{\omega}}\l(\frac{\partial^2 \ell(\bm{\omega})}{\partial \sigma^2_t}\r) &= \frac{-2N}{\sigma^2_t}.
\end{align*}
Therefore, the Fisher information matrix is given by:
\begin{align*}
    I_{\E}(\bm{\omega}) &= \begin{pmatrix}
        \frac{N}{\sigma^2_1} & 0 & \dots & \dots & \dots &\dots & \dots & 0 \\
        0 & \frac{N}{\sigma^2_2} & 0 & \dots & \dots & \dots & \dots &  0 \\
        \vdots &  & \ddots & & & & & \vdots \\
        0 & \dots & 0 & \frac{N}{\sigma^2_T} & 0   & \dots& \dots & 0 \\
        0 & \dots & \dots & 0 & \frac{2N}{\sigma^2_1} & 0 & \dots & 0 \\
        0 & \dots & \dots & \dots & 0 & \frac{2N}{\sigma^2_2}  & \dots & 0 \\
        \vdots & & & & &  &\ddots & \vdots \\
         0 & \dots & \dots & \dots & \dots & \dots & 0 &   \frac{2N}{\sigma^2_T}\\
    \end{pmatrix}
\end{align*}
Hence, by standard results in maximum likelihood theory (see e.g.~\cite{MillarRussellB2011MLEa}), we have that (for large $N$):
\begin{equation*}
    \hat{\bm{\omega}} - \bm{\omega}^0 \sim \m{N}\l(0, I_{\E}^{-1}(\bm{\omega}^0)\r),
\end{equation*}
which is asymptotically equivalent to:
\begin{equation*}
    \hat{\bm{\omega}} - \bm{\omega}^0 \sim \m{N}\l(0, I_{\E}^{-1}(\bm{\hat{\omega}})\r).
\end{equation*}
Therefore, for large $N$ we have:
\begin{equation*}
    (\hat{\bm{\omega}} - \bm{\omega}^0)I_{\E}(\hat{\bm{\omega}})(\hat{\bm{\omega}} - \bm{\omega}^0)^T \sim \chi^2_{2T}.
\end{equation*}
In other words:
\begin{equation*}
    \sum_{t=1}^T\l(\frac{N}{\hat{\sigma}^2_t}\l(\hat{\mu}_t - \mu^0_t\r)^2 + \frac{2N}{\hat{\sigma}^2_t}\l(\hat{\sigma}_t - \sigma^0_t\r)^2\r) \sim \chi^2_{2T}.
\end{equation*}
Thus, we can create an approximate $100(1-\alpha)\%$ confidence set for $\bm{\omega}^0$ using:
\begin{equation*}
    \Omega_{\alpha} = \l\{(\bm{\mu}, \bm{\sigma}) \in \R^{2T}: \sum_{t=1}^T\l(\frac{N}{\hat{\sigma}^2_t}\l(\hat{\mu}_t - \mu_t\r)^2 + \frac{2N}{\hat{\sigma}^2_t}\l(\hat{\sigma}_t - \sigma_t\r)^2\r) \le \chi^2_{2T, 1-\alpha} \r\}.
\end{equation*}
\end{proof}

\subsection{Piecewise linear approximation of DRO model}\label{sec:normal_PWL}

Since the objective function~(10) is nonlinear, and also not quadratic, we
use a piecewise linear approximation (PLA) to the problem to solve it. In order to create a PLA of the
model, we must define a PLA of the entire non-linear part of the objective function. Firstly, for $i = 1,\dots,|\Omega^{\prime}_{\alpha}|$
define $\beta^i_t = \frac{ \sum_{k=1}^t (\mu^i_k - q_k)}{\sqrt{\sum_{k=1}^t
(\sigma^i_k)^2}}$, where $\bm{\omega}^i$ is the $i$th element of $\Omega^{\prime}_{\alpha}$. Let $g^i_t(\beta^i_t)$ be the non-linear part of the objective function for day $t$, under parameter $\bm{\omega}^i \in \Omega^{\prime}_{\alpha}$. Then, $g^i_t$ can be written as
\begin{equation*}
    g^i_t(\beta^i_t) = \sqrt{\sum_{k=1}^t (\sigma^i_k)^2}\l\{(a_t \Phi(\beta^i_t) - h) \beta^i_t  + a_t \phi(\beta^i_t) \r\}, \fa t = 1,\dots,T \fa i = 1,\dots,|\Omega^{\prime}_{\alpha}|.
\end{equation*}
In order to define a piecewise linear approximation of $g^i_t$, we need to define a set of points $Z^{i, t} = \l\{z^{i, t}_1,\dots,z^{i, t}_J\r\}$ such that $z^{i, t}_1 \le \beta^i_t \le z^{i, t}_J \fa \bm{q}$. In our experiments, we will assume that these points are equally spaced, i.e.\ $z^{i, t}_j = z^{i, t}_{j-1} + \varepsilon$,  where $\varepsilon$ is referred to as a \textit{gap}. Since $\beta^i_t = \frac{\sum_{k=1}^t(\mu^i_k - q_k)}{\sqrt{\sum_{k=1}^t (\sigma^i_t)^2}}$, we have:
\begin{equation*}
    (\beta^i_t)^{\min} := \frac{\sum_{k=1}^t\l(\mu^i_k - \frac{W}{w_k}\r)}{\sqrt{\sum_{k=1}^t (\sigma^i_t)^2}} \le \beta^i_t
    \le \frac{\sum_{k=1}^t\mu^i_k}{\sqrt{\sum_{k=1}^t (\sigma^i_t)^2}} =: (\beta^i_t)^{\max}
\end{equation*}
for each $i = 1,\dots,|\Omega^{\prime}_{\alpha}|$ and $t=1,\dots,T$. The first inequality is due to two reasons. Firstly, the maximum value $q_k$ can take is when all other orders are zero. Secondly, this means that the budget constraint becomes $w_k q_k \le W$ i.e.\ $q_k \le \frac{W}{w_k}$. The second inequality is true since $q_k \ge 0 \fa k$. We can therefore define $Z^{i, t}$ as follows:
\begin{equation*}
    Z^{i, t} = \l\{(\beta^i_t)^{\min} + j\varepsilon: j=0,\dots,\l\lfloor \frac{(\beta^i_t)^{\max} - (\beta^i_t)^{\min}}{\varepsilon}\r\rfloor\r\}.
\end{equation*}
The reason for having a different range for each $(i,t)$ pair is that if we were to use only one, this would be the largest range of all of the individual ranges. Having one for each pair allows us to reduce the number of points in total that are used to define the piecewise linear model.  Given the sets $Z^{i, t}$, a piecewise linear approximation of $g^i_t$ is given by:
\begin{equation*}\label{eq:PWL_CDF}
    G^i_t(\beta^i_t) = \frac{\beta^i_t-z^{i, t}_{j+1}}{z^{i, t}_j - z^{i, t}_{j+1}}g^i_t(z^{i, t}_j) + \frac{\beta^i_t-z^{i, t}_j}{z^{i, t}_{j+1} - z^{i, t}_j} g_t(z^{i, t}_{j+1}),
\end{equation*}
where $j \in \{1,\dots, J - 1\}$ is chosen such that $\beta^i_t \in [z_j, z_{j+1}]$. Given this, we can construct a piecewise linear approximation of the DRO model via the following. Firstly, create dual-indexed dummy decision variables
$\tilde{\bm{\beta}} = (\tilde{\beta}^i_t)_{i \in \{1,\dots,|\Omega^{\prime}_{\alpha}|\}, t \in \m{T}}$ and $\bm{\tilde{g}} = (\tilde{g}^i_t)_{i \in \{1,\dots,|\Omega^{\prime}_{\alpha}|\}, t \in \m{T}}$ such
that an index $i$ corresponds to the value of the variables when the distribution
is given by $\bm{\omega}^i \in \Omega^{\prime}_{\alpha}$. Secondly, define a dummy
variable $\theta$ to represent the worst-case cost. Then, our piecewise linear model is given by:
\begin{align}
    \min_{\bm{q}, \bm{\tilde{g}}, \bm{\tilde{\beta}}} \quad  &\theta \label{eq:DRO_PWL_obj} \\
    \text{s.t. } &\theta \ge \sum_{t=1}^T\l\{\tilde{g}^i_t + q_t w_t - p \mu^i_t\r\} \fa i = 1,\dots,|\Omega^{\prime}_{\alpha}|\label{eq:inner_obj}\\
    &\sum_{t=1}^T w_t q_t \le W,\\
    &\tilde{g}^i_t = G^i_t(\tilde{\beta}^i_t) \fa t \in \m{T}, \fa i = 1,\dots,|\Omega^{\prime}_{\alpha}|\label{eq:fix_g}\\
    &\tilde{\beta}^i_t = \frac{\sum_{k=1}^t (\mu^i_k -
    q_k)}{\sqrt{\sum_{k=1}^t(\sigma^i_{k})^2}} \fa t \in \m{T}, \fa i = 1,\dots,|\Omega^{\prime}_{\alpha}|\label{eq:fix_alpha}\\
    &q_t \ge 0 \fa t \in \m{T}. \label{eq:DRO_PWL_lastconstraint}
\end{align}
Our numerical experiments are conducted using Gurobi and~\eqref{eq:DRO_PWL_obj}-\eqref{eq:DRO_PWL_lastconstraint} is presented for reproducibility when using this solver. Here, \eqref{eq:fix_alpha} ensures that $\tilde{\beta}^i_t = \beta^i_t$. We do not technically need a dummy variable for $\beta^i_t$ as it is linear in the $q_k$, but Gurobi requires the PLA function $G_t$ to be evaluated at a decision variable. Then, \eqref{eq:fix_g} ensures that $\tilde{g}^i_t$ is equal to the PLA of $g_t$ evaluated at the dummy variable $\tilde{\beta}^i_t$. Constraint~\eqref{eq:inner_obj} ensures that the value of $\theta$ returned by the model is equal to the worst-case expected cost for the chosen $\bm{q}$.

\subsection{Proof of Theorem~4.3}\label{sec:normal_convex_proofs}

\begin{proof}
    In order to differentiate the objective function, we first differentiate the various terms. Let $\tau \in \{1,\dots,T\}$.
Then $\mu_\tau$ and $\sigma_\tau$ appear in the $t^{\text{th}}$ term of the sum
for each $t \ge \tau$. In fact, recalling that $\beta_t =
\frac{-\sum_{k=1}^t(q_k - \mu_k)}{\sqrt{\sum_{k=1}^t \sigma_k^2}}$, we have
\begin{align*}
    \pderiv{\beta_t}{\mu_\tau} &= \begin{cases}
        0 & \text{ if } \tau > t, \\
        \frac{1}{\sqrt{\sum_{k=1}^t\sigma^2_k}}, & \text{ if } \tau \le t.
    \end{cases}
    \\
    \pderiv{\beta_t}{\sigma_\tau} &= \begin{cases}
        0 & \text{ if } \tau > t,\\
        \frac{\sum_{k=1}^t(q_k - \mu_k)}{\sigma_\tau\l(\sum_{k=1}^t \sigma^2_k \r)^{\frac{3}{2}}} &\text{ if } \tau \le t.
    \end{cases}
\end{align*}
Using the chain rule, we then find that:
\begin{align*}
    \pderiv{\Phi(\beta_t)}{\mu_\tau} &= \frac{\phi(\beta_t)}{\sqrt{\sum_{k=1}^t\sigma^2_k}} \\
    \pderiv{\Phi(\beta_t)}{\sigma_\tau} &= \frac{\sum_{k=1}^t(q_k - \mu_k)}{\sigma_\tau\l(\sum_{k=1}^t \sigma^2_k \r)^{\frac{3}{2}}} \phi(\beta_t) \\
    &= \frac{-\beta_t}{\sigma_\tau \sum_{k=1}^t \sigma^2_k } \phi(\beta_t) 
\end{align*} 
for $\tau \le t$. Furthermore, we have:
\begin{align*}
    \pderiv{\phi(\beta_t)}{\mu_\tau} &= \frac{1}{\sqrt{\sum_{k=1}^t\sigma^2_k}} (-\beta_t\phi(\beta_t))\\
    \pderiv{\phi(\beta_t)}{\sigma_\tau} &= \frac{-\beta_t}{\sigma_\tau \sum_{k=1}^t \sigma^2_k} (-\beta_t\phi(\beta_t))  \\
    &= \frac{1}{\sigma_\tau \sum_{k=1}^t \sigma^2_k}\beta_t^2 \phi(\beta_t) 
\end{align*} 
Now, we can differentiate the objective as given in Lemma~4.1 to find:
\begin{equation*}
    \pderiv{C_{\bm{\omega}}(\bm{q})}{\mu_\tau} \hiderel{=} \sum_{t=\tau}^T\l\{
    (a_t \Phi(\beta_t) - h) + a_t \frac{\phi(\beta_t)}{\sqrt{\sum_{k=1}^t \sigma^2_k}}\sum_{k=1}^t(\mu_k -
    q_k) - a_t \frac{\beta_t \phi(\beta_t)}{\sqrt{\sum_{k=1}^t \sigma^2_k}}\sqrt{\sum_{k=1}^t \sigma^2_k} - p
    \r\}
\end{equation*}
and we can simplify this to:
\begin{align*}
    \pderiv{C_{\bm{\omega}}(\bm{q})}{\mu_\tau} &= \sum_{t=\tau}^{T}\l\{
    - h + a_t\l(\Phi(\beta_t)  + \frac{\phi(\beta_t)}{\sqrt{\sum_{k=1}^t
    \sigma^2_k}}\sum_{k=1}^t(\mu_k - q_k)\r) - a_t \beta_t \phi(\beta_t) - p
    \r\}\\
    &= \sum_{t=\tau}^T \l\{-h + a_t\l(\Phi(\beta_t)+\beta_t\phi(\beta_t)\r) - a_t\beta_t
    \phi(\beta_t) - p\r\}\\
    &= \sum_{t=\tau}^T \l\{a_t\Phi(\beta_t) - h - p\r\} \label{eq:MPNVP_mu_deriv}
\end{align*}
Furthermore, with $s_t = \sqrt{\sum_{k=1}^t\sigma^2_t}$ we have:
\begin{equation*}
    \frac{\partial^2 C_{\bm{\omega}}(\bm{q})}{\partial \mu^2_\tau} = \sum_{t=\tau}^T a_t \frac{\phi(\beta_t)}{s_t} \ge 0
\end{equation*}
i.e.\ the objective is convex in each $\mu_\tau$.
We also have:
\begin{align*}
    \pderiv{C_{\bm{\omega}}(\bm{q})}{\sigma_\tau} &= \sum_{t=\tau}^T\l\{
    a_t\l(\frac{-\beta_t\phi(\beta_t)}{\sigma_\tau \sum_{k=1}^t\sigma^2_k}\r)\sum_{k=1}^t(\mu_k - q_k) + \frac{a_t \phi(\beta_t)\sigma_\tau}{\sqrt{\sum_{k=1}^t\sigma^2_k}} + a_t\sqrt{\sum_{k=1}^t\sigma^2_k}\frac{\alpha^2_t\phi(\beta_t)}{\sigma_\tau \sum_{k=1}^t \sigma^2_k}     
    \r\}\\
    &= \sum_{t=\tau}^T\l\{a_t\l(\frac{-\alpha^2_t \phi(\beta_t)}{\sigma_\tau \sqrt{\sum_{k=1}^t\sigma^2_k}}\r) + \frac{a_t\phi(\beta_t)\sigma_\tau}{\sqrt{\sum_{k=1}^t\sigma^2_k}} + a_t \frac{\alpha^2_t\phi(\beta_t)}{\sigma_\tau \sqrt{\sum_{k=1}^t\sigma^2_k}}\r\}\\
    &= \sum_{t=\tau}^T \frac{a_t \phi(\beta_t)\sigma_\tau}{ \sqrt{\sum_{k=1}^t\sigma^2_k}}\label{eq:MPNVP_sig_deriv}\\
    & \ge 0 \fa \tau \in \m{T}.
\end{align*}
Hence, $C_{\bm{\omega}}(\bm{q})$ is increasing in $\sigma_\tau$ for all $\tau
\in \m{T}$. 
\end{proof}

\section{Appendix for Poisson demands}\label{sec:pois_proofs}

\subsection{Proof of Lemma~5.1}\label{sec:pois_obj_proof}

\begin{proof}
    In order to evaluate the objective function, we need to evaluate $\E_{\bm{\lambda}}[I^+_t]$ and $\E_{\bm{\lambda}}[I^-_t]$. Recall that:
\begin{align*}
    \E_{\bm{\lambda}}[I^+_t] &= \E_{\bm{\lambda}}[\max\{I_t, 0\}] \\
    \E_{\bm{\lambda}}[I^-_t] &= \E_{\bm{\lambda}}[\max\{-I_t, 0\}] 
\end{align*}
where $I_t =  \sum_{k=1}^t(q_t - X_t)$, $\P(I_t \ge 0) = \P\l(\sum_{k=1}^t X_k \le \sum_{k=1}^t q_k\r)$ and $\E_{\bm{\lambda}}[I_t|I_t \ge 0] = \sum_{k=1}^t q_k - \E_{\bm{\lambda}}[\sum_{k=1}^t X_k| \sum_{k=1}^t X_k \le \sum_{k=1}^t q_k]$. To simplify notation, let us define $Q_t = \sum_{k=1}^t q_k$ and $\Lambda_t = \sum_{k=1}^t \lambda_k$. Now, we can calculate the expected values as follows:
\begin{align*}
    \E_{\bm{\lambda}}[\max\{I_t, 0\}] &= \E_{\bm{\lambda}}\l[\max\l\{Q_t - \sum_{k=1}^t X_k, 0\r\}\r] \\
    &= \sum_{x = 0}^{\infty} \max\{Q_t - x, 0\} \tilde{f}_{t}(x) \\
    &= \sum_{x=0}^{Q_t} (Q_t - x) \tilde{f}_t(x)\\
    &= Q_t\tilde{F}_t(Q_t) - \sum_{x=0}^{Q_t} x \frac{(\Lambda_t)^{x}\exp{(-\Lambda_t)}}{(x)!} \\
    &= Q_t\tilde{F}_t(Q_t) - \sum_{x=1}^{Q_t}  \frac{(\Lambda_t)^{x}\exp{(-\Lambda_t)}}{(x-1)!} \\
    &= Q_t \tilde{F}_t(Q_t) - \Lambda_t\sum_{x=1}^{Q_t}  \frac{(\Lambda_t)^{x-1}\exp{(-\Lambda_t)}}{(x-1)!} \\
    &= Q_t \tilde{F}_t(Q_t) - \Lambda_t\tilde{F}_t(Q_t-1),
\end{align*}
where $\tilde{f}_t$ and $\tilde{F}_t$ are the PMF and CDF of $\sum_{k=1}^t X_k$. Similarly we have:
Similarly:
\begin{align*}
    \E_{\bm{\lambda}}[\max\{-I_t, 0\}] &= \sum_{x = 0}^{\infty} \max\{x - Q_t, 0\} \tilde{f}_{t}(x) \\
    &= \sum_{x=Q_t + 1}^{\infty}(x-Q_t)\tilde{f}_t(x) \\
    &= -Q_t(1-\tilde{F}_t(Q_t)) + \sum_{x=Q_t+1}^{\infty} x\frac{(\Lambda_t)^{x}\exp{(-\Lambda_t)}}{(x)!}\\
    &= -Q_t(1-\tilde{F}_t(Q_t)) + \Lambda_t - \sum_{x=0}^{Q_t} x\frac{(\Lambda_t)^{x}\exp{(-\Lambda_t)}}{(x)!}\\
    &= -Q_t(1-\tilde{F}_t(Q_t)) + \Lambda_t - \sum_{x=1}^{Q_t} \frac{(\Lambda_t)^{x}\exp{(-\Lambda_t)}}{(x-1)!}\\
    &= -Q_t(1-\tilde{F}_t(Q_t)) +  \Lambda_t(1-\tilde{F}_t(Q_t-1)). 
\end{align*}
This gives the following objective function (with $a_t = b +
p\mathds{1}\{t=T\} + h$):
\begin{align*}
    C_{\bm{\lambda}}(\bm{q}) &= p\E_{\bm{\lambda}}[I^-_T] +  \sum_{t=1}^T(h \E_{\bm{\lambda}}[I^+_t] +  b\E_{\bm{\lambda}}[I^-_t] + w_t q_t - p \E_{\bm{\lambda}}[X_t]) \\
    &= \sum_{t=1}^T(h \E_{\bm{\lambda}}[I^+_t] +  (b + p\mathds{1}\{t=T\})\E_{\bm{\lambda}}[I^-_t] + w_t q_t - p \E_{\bm{\lambda}}[X_t]) \\
    &= \sum_{t=1}^T\l(h (Q_t \tilde{F}_t(Q_t) - \Lambda_t\tilde{F}_t(Q_t-1)) +  (b + p\mathds{1}\{t=T\})(-Q_t(1-\tilde{F}_t(Q_t)) \r.\nonumber\\ &\l. \quad +\Lambda_t(1-\tilde{F}_t(Q_t-1))) + w_t q_t - p \lambda_t\r) \\
    &= \sum_{t=1}^T\l(a_t Q_t \tilde{F}_t(Q_t) - \Lambda_t a_t \tilde{F}_t(Q_t-1) + (b+p\mathds{1}\{t=T\})(\Lambda_t - Q_t) + w_t q_t - p \lambda_t\r),
\end{align*}
as required.
\end{proof}

\subsection{Proof of Proposition~5.2}\label{sec:pois_conf_proof}

\begin{proof}
By standard results in maximum likelihood estimation~\citep{MillarRussellB2011MLEa}, for large $N$ we have the approximate result that:
\begin{equation*}
    \hat{\lambda}_t \sim \m{N}\l(\lambda^0_t, \frac{\lambda^0_t}{N}\r),
\end{equation*}
or equivalently:
\begin{equation*}
    \hat{\lambda}_t \sim \m{N}\l(\lambda^0_t, \frac{\hat{\lambda}_t}{N}\r).
\end{equation*}
Now, by independence of the $T$ MLEs, we have the approximate result that:
\begin{equation*}
    \sum_{t=1}^T \frac{N}{\hat{\lambda_t}}(\hat{\lambda}_t - \lambda^0_t)^2 \le \chi^2_T.
\end{equation*}
Hence, we can calculate an approximate $100(1-\alpha)\%$ confidence set for $\bm{\lambda}^0$ using:
\begin{equation*}
    \Omega_{\alpha} = \l\{\bm{\lambda} \in \R^T_{+}: \sum_{t=1}^T \frac{N}{\hat{\lambda}_t}(\hat{\lambda}_t - \lambda^0_t)^2 \le \chi^2_{T, 1-\alpha}\r\}.
\end{equation*}
\end{proof}

\subsection{Piecewise linear DRO model}\label{sec:poisson_PWL}

Note that, for Poisson demands, the objective function is already piecewise linear. The Poisson CDF $\tilde{F}_t$ is piecewise constant and hence piecewise linear. In addition, the non-linear term $Q_t \tilde{F}_t(Q_t)$ is piecewise linear, since $\tilde{F}_t(Q_t)$ is piecewise constant. Since $\tilde{F}_t(Q_t)$ is constant on the interval $[Q_t, Q_t+1)$ for any $Q_t \in \mathbb{N}_0$, the objective function is linear in these intervals. Therefore, we only need to consider $\varepsilon = 1$ as the gap between points for this model, and in addition the piecewise linear objective function is not approximate. We now present the piecewise linear model formulation, for reproducibility when using Gurobi. Despite the fact that the cost function is already piecewise linear, we still use Gurobi's piecewise linear constraints in order to formulate the model. In order to do so, define $\bm{\lambda}^i$ as the $i$th element of $\Omega^{\prime}_{\alpha}$ and define:
\begin{equation*}
    D^i_t(Q) =  a_t Q \tilde{F}^i_t(Q) - \Lambda^i_t a_t \tilde{F}^i_t(Q-1) + (b+p\mathds{1}\{t=T\})(\Lambda^i_t - Q) + w_t Q - p \lambda^i_t
\end{equation*}
for each $i = 1,\dots|\Omega^{\prime}_{\alpha}|$ and $t = 1,\dots,T$,
where $\tilde{F}^i_t$ is the CDF of $\sum_{k=1}^t X_k$ under $\bm{\lambda}^i$. To define the piecewise linear constraint in Gurobi, we need a set $U^t = \l\{u^t_1,\dots,u^t_{J^{\prime}}\r\}$  for each $t \in \m{T}$ such that $u^t_1 \le Q_t \le u^{t}_{J^{\prime}}$ holds for any $\bm{q}$. Since $Q_t = \sum_{k=1}^t q_k$, it is easy to see that:
$$0 \le Q_t \le \sum_{k=1}^t \frac{W}{w_k} \fa t = 1,\dots,T.$$
Hence, we take $U^t = \l\{0,1,\dots, \sum_{k=1}^t \frac{W}{w_k}\r\}$ as our points for each $t$. Then, the model can be formulated using:
\begin{align*}
    \min_{\bm{q}, \bm{Q}, \theta, h} \ &\theta \\
    \text{s.t. } & \theta \ge \sum_{t=1}^T d^i_t \fa i =1,\dots,|\Omega^{\prime}_{\alpha}|,\\
    &Q_t = \sum_{k=1}^t q_t \fa t =1,\dots,T,\\
    &d^i_t = D^i_t(Q_t), \fa i =1,\dots,|\Omega^{\prime}_{\alpha}| \fa t=1,\dots,T\label{eq:pois_PWL_constr} \\
    &\sum_{t=1}^T w_t q_t \le W \\
    & q_t \in \mathbb{N}_0\fa t=1,\dots,T.
\end{align*}

\subsection{Proof of Theorem~5.3}\label{sec:pois_convex_proof}

\begin{proof}
    In order to differentiate $C_{\bm{\lambda}}(\bm{q})$, we first differentiate the Poisson CDF to find that:
\begin{align*}
    \pderiv{\tilde{F}_t(Q_t)}{\lambda_k} &=  \pderiv{}{\lambda_k}\l(\sum_{x=0}^{Q_t} \frac{\exp\l(-\Lambda_t\r)\l(\Lambda_t\r)^x}{x!}\r)\\
    &= \sum_{x=0}^{Q_t} \frac{-\exp\l(-\Lambda_t\r)\l(\Lambda_t\r)^x + \exp\l(-\Lambda_t\r)x\l(\Lambda_t\r)^{x-1}}{x!}\\
    &= \sum_{x=1}^{Q_t}\frac{\exp\l(-\Lambda_t\r)\l(\Lambda_t\r)^{x-1}}{(x-1)!} - \sum_{x=0}^{Q_t} \frac{\exp\l(-\Lambda_t\r)\l(\Lambda_t\r)^x}{x!} \\
    &= \sum_{x=0}^{Q_t - 1}\frac{\exp\l(-\Lambda_t\r)\l(\Lambda_t\r)^{x}}{x!} - \sum_{x=0}^{Q_t} \frac{\exp\l(-\Lambda_t\r)\l(\Lambda_t\r)^x}{x!} \\
    &= \tilde{F}_t(Q_t - 1) - \tilde{F}_t(Q_t)\\
    &= -\tilde{f}_t(Q_t)
\end{align*}
for $k \le t$. Similarly,
\begin{equation*}
    \pderiv{\tilde{F}_t(Q_t - 1)}{\lambda_k} =  -\tilde{f}_t(Q_t-1) \quad (k \le t).
\end{equation*}
Hence, using the expression given in Lemma~5.1, we have:
\begin{align}
    \pderiv{C_{\bm{\lambda}}(\bm{q})}{\lambda_k} &= \sum_{t=k}^T\l(
        - a_t Q_t \tilde{f}_t(Q_t) + \Lambda_t a_t \tilde{f}_t(Q_t - 1) - a_t \tilde{F}_t(Q_t - 1) + (a_t - h) - p 
    -p\r) \\
    &= \sum_{t=k}^T\l(- a_t \tilde{F}_t(Q_t - 1) + (a_t - h) - p\r) \label{eq:cancel_poisson}
\end{align}
where~\eqref{eq:cancel_poisson} is true because:
\begin{align*}
    Q_t \tilde{f}_t(Q_t) &= Q_t \frac{\exp(-\Lambda_t) (\Lambda_t)^{Q_t}}{Q_t!}\\
    &= \frac{\exp(-\Lambda_t) (\Lambda_t)^{Q_t}}{(Q_t - 1)!}\\
    &= \Lambda_t \tilde{f}_t(Q_t-1).
\end{align*}
Differentiating again, we find:
\begin{align*}
    \frac{\partial^2 C_{\bm{\lambda}}(\bm{q})}{\partial \lambda^2_k} &= \sum_{t=k}^T a_t \tilde{f}_t(Q_t - 1)  \\
    &\ge 0 \fa k = 1,\dots,T.
\end{align*}
Therefore, the objective function is convex in $\lambda_k$ for each $k$.
\end{proof}
\end{appendices}

\end{document}